\documentclass[reqno,11pt]{amsart}
\usepackage{amssymb, amsmath,latexsym,amsfonts,amsbsy, amsthm}
\usepackage{txfonts}
\allowdisplaybreaks[4]

\newcommand{\beq}{\begin{equation}}
\newcommand{\eeq}{\end{equation}}
\newcommand{\ben}{\begin{eqnarray}}
\newcommand{\een}{\end{eqnarray}}
\newcommand{\beno}{\begin{eqnarray*}}
\newcommand{\eeno}{\end{eqnarray*}}

\renewcommand{\theequation}{\thesection.\arabic{equation}}

\newtheorem{theorem}{Theorem}[section]

\newtheorem{lemma}[theorem]{Lemma}

\newtheorem{Theorem}{Theorem}[section]
\newtheorem{Definition}[Theorem]{Definition}
\newtheorem{Proposition}[Theorem]{Proposition}
\newtheorem{Lemma}[Theorem]{Lemma}

\newtheorem{Remark}[Theorem]{Remark}

\begin{document}
\title[Prandtl--Batchelor flows with a point vortex]
{Prandtl--Batchelor flows with a point vortex on a disk}

\author{Zhi Chen}
\address{School of  Mathematics and Statistics, Anhui Normal University, Wuhu 241002, China}
\email{zhichenmath@ahnu.edu.cn}

\author{Mingwen Fei}
\address{School of  Mathematics and Statistics, Anhui Normal University, Wuhu 241002, China}
\email{mwfei@ahnu.edu.cn}

\author{Zhiwu Lin}
\address{School of Mathematical Sciences, Fudan University,
Shanghai 200433, China}
\email{zwlin@fudan.edu.cn}

\author{Jianfeng Zhao}
\address{School of Mathematics and Information Science, Guangxi University, Nanning 530004, China}
\email{zhaojianfeng@amss.ac.cn}

\date{\today}

\renewcommand{\theequation}{\thesection.\arabic{equation}}
\setcounter{equation}{0}
\begin{abstract}
We study steady two-dimensional incompressible Navier--Stokes flow in the unit disk in a high Reynolds number regime, driven by a localized central forcing and a nearly rigid rotating boundary. The forcing models a compact source of circulation that sustains a vortex core in a viscous flow. For sufficiently small viscosity and boundary perturbation, we construct a very weak steady solution on the whole disk and identify its inviscid limit: a point vortex embedded in a constant vorticity background and joined to the boundary motion by a thin boundary layer. The background vorticity is selected by the Batchelor--Wood formula. The singular core introduces non-coercive interactions in the linearized problem that are absent for regular Prandtl--Batchelor flows. To overcome this difficulty, we map the punctured disk to a semi-infinite cylinder by a logarithmic radial transformation and develop a Fourier-mode analysis that isolates the delicate first Fourier mode, yielding uniform stability and vorticity estimates in the presence of the point-vortex singularity.
\end{abstract}
\smallbreak
      \smallskip

\maketitle

\numberwithin{equation}{section}
\setcounter{secnumdepth}{3}

\indent

\numberwithin{equation}{section}

\indent

\section{Introduction}

\subsection{Prandtl--Batchelor selection with a singular vortex core}

The problem considered in this paper lies at the intersection of two classical
themes in fluid mechanics: the vanishing-viscosity limit near solid
boundaries and the persistence of concentrated vortical structures.  Since
Prandtl's boundary-layer theory \cite{prandtl}, a basic expectation has been
that, at large Reynolds number, a steady Navier--Stokes flow in a bounded
domain is described in the bulk by an Euler flow, with the mismatch at the
no-slip wall resolved inside a thin viscous layer.  In a recirculating region
filled by nested closed streamlines, however, the limiting Euler flow is not
arbitrary.  Batchelor \cite{B} and Feynman--Lagerstrom
\cite{feymann-lagerstrom} observed that the weak viscous effect, accumulated
along closed streamlines, imposes a solvability condition that forces the
limiting core vorticity to be constant.  For circular configurations, Batchelor \cite{B} and Wood
\cite{W} identified the relation between this constant and the tangential
velocity prescribed at the wall.  This Prandtl--Batchelor selection
principle is one of the basic mechanisms by which a small viscosity selects
a distinguished steady Euler state.

Rigorous work on this selection problem has so far focused principally on
regular eddy cores or on geometries that avoid a singular center.  Kim
\cite{K1998,K2000,kim-formula} initiated the mathematical study of
Prandtl--Batchelor flows in a disk.  For wall velocities close to rigid
rotation, the convergence of the Prandtl expansion was proved on the disk in
\cite{FGLT}, where the selected outer flow has a smooth solid-body core, and
on an annulus in \cite{FGLT2}, where Couette-type profiles are permitted but
the inner boundary removes the vortex center.  The present paper treats the
singular-core disk problem left between these two regimes: the limiting
Euler flow contains a point-vortex component whose velocity behaves as
\(b/r\) all the way to the origin.

This singular core is natural from the viewpoint of coherent vortex
modeling.  Beginning with Helmholtz \cite{H}, point vortices have served as
idealized descriptions of sharply concentrated circulation: away from the
small vortex core, a localized vortex induces a tangential velocity whose
leading behavior is proportional to \(1/r\).  Such idealizations are used to
describe transport and interaction by coherent vortices in two-dimensional
and geophysical fluid dynamics; see, for instance, \cite{D2022,Gallay,
GMS,Kizner,Lugt,Provenzale,wmz06}.  In the setting of the present work, the
point-vortex component models a localized core, while the regular
solid-body component represents the surrounding eddy whose vorticity is
selected in the high-Reynolds-number limit.

There is, nevertheless, an essential distinction between prescribing a point
vortex in an inviscid Euler model and constructing a \emph{steady viscous}
flow containing a point-vortex-like core.  In the Euler description, a point vortex may be prescribed directly as a concentrated circulation carried by the flow.  In a viscous fluid confined by a no-slip outer boundary, however, a persistent concentrated vortex cannot be regarded merely as an initial structure: viscous transport and wall dissipation tend to diffuse and dissipate its circulation.  A steady solution with a fixed central vortex strength should therefore be understood as the response to a localized source of angular
momentum or circulation.  In our model this source is represented by a
distributional force supported at the center, obtained conceptually as the
vanishing-radius limit of a compact forcing region with fixed integrated
circulation input.

This interpretation has concrete laboratory and engineering analogues.  A
small mechanical or magnetic stirrer placed near the center of a cylindrical
vessel injects angular momentum locally and produces a coherent rotating
core together with a surrounding circulation.  Experiments with a magnetic
stirrer show that outside a small central core the measured tangential
velocity is well described by an ideal-vortex profile \cite{Halasz}.  Thus,
after the details of the actuator are removed, a localized forcing producing
a \(1/r\) bulk swirl is a natural reduced model for a steadily driven vortex
in a mixing or actuation device. There is also a geophysical motivation: in rotating or stratified geophysical flows, localized topographic features, such as seamounts or isolated obstacles, can organize long-lived vortical motions and exchange momentum with the ambient flow \cite{SDR,DN}. Although the physical forcing in such problems is distributed over the obstacle and depends on additional effects not included here, its leading influence on the surrounding circulation motivates an idealized localized torque or circulation input. Accordingly, the point-vortex force in this paper is not intended as a literal model of a particular device or topographic configuration; rather, it isolates the analytical question of how a steadily maintained singular vortex core interacts with a viscous boundary layer and with the Prandtl--Batchelor selection of the regular background vorticity.

Our goal is therefore a singular-core version of the Prandtl--Batchelor
selection problem: we prescribe a centered localized circulation input,
construct the resulting steady viscous flow, and prove that its inviscid
limit consists of a point vortex embedded in a constant-vorticity background
selected by the rotating wall.  The central question is not whether the
singular circulation persists---it is imposed by the forcing---but whether
the surrounding flow is regular and continues to obey the Batchelor--Wood selection
law in the presence of the $1/r$ core singularity.

\subsection{The localized forcing model}

Let $B_1=B_1(0)\subset\mathbb R^2$ be the unit disk.  We consider the steady
forced Navier--Stokes equations
\begin{eqnarray}
\left \{
\begin {array}{ll}
\mathbf{u}^\varepsilon\cdot\nabla\mathbf{u}^\varepsilon+\nabla p^\varepsilon-\varepsilon^2\Delta\mathbf{u}^\varepsilon=\mathbf{f},\\[5pt]
\nabla\cdot \mathbf{u}^\varepsilon=0,\\[5pt]
\mathbf{u}^\varepsilon\big|_{\partial B_1}=\mathbf{g}.
\end{array}
\right.\label{ns}
\end{eqnarray}
Here $\varepsilon^2>0$ is the reciprocal Reynolds number.  The external force
is chosen at the viscous scale and is supported at the origin:
\begin{align}\label{force-condition}
\mathbf{f}=\varepsilon^2\mathbf{F}=\varepsilon^2F\mathbf{t}.
\end{align}
It sustains the singular circulation associated with the point vortex
\begin{align}\label{point-vortex}
\omega_{\rm pv}(x,y)=2\pi b\delta_0,\qquad b>0.
\end{align}
More precisely, $\mathbf{F}\in (W_0^{1,q'}\cap W^{2,q'})'$ for $q'>2$,
is supported at $\{0\}$, and is defined by
\begin{align}\label{point-vortex-0}
\iint_{B_1}\mathbf{F}\cdot {\boldsymbol{\varphi}}\,dxdy
=-4b\pi(\partial_r\overline{\varphi_1})(0)
=-2b\big(\nabla{\boldsymbol{\varphi}}(0,0):\mathbf{I_0}\big)
\end{align}
for divergence-free test fields
${\boldsymbol{\varphi}}=\varphi_1 \mathbf{t}+\varphi_2\mathbf{n}
\in W_0^{1,q'}\cap W^{2,q'}$ satisfying
$\boldsymbol{\varphi}(0,0)=0$, where
\[
 \overline{\varphi_1}(r)=\frac{1}{2\pi}
 \int_0^{2\pi}\varphi_1(\theta,r)\,d\theta,
 \qquad
 \mathbf{I_0}=
 \begin{pmatrix}0&\pi\\-\pi&0\end{pmatrix}
 =\int_0^{2\pi}\mathbf{n}\otimes\mathbf{t}\,d\theta .
\]
The force $\mathbf F$ should therefore be read as the effective localized
source that sustains the prescribed singular circulation; it is not a force
exerted by a point vortex regarded as a material particle.  The scaling
$\mathbf f=\varepsilon^2\mathbf F$ places this sustaining input at the same
order as viscous effects and permits a nontrivial Euler limit outside the
boundary layer.

\subsection{Very weak solvability}

The classical theory of steady Navier--Stokes boundary-value problems in
bounded $C^{1,1}$ domains (see, e.g., \cite{La,Leray,T}) yields weak solutions
in $W^{1,q}$ for data
$\mathbf f\in W^{-1,q}$ and $\mathbf g\in W^{1-1/q,q}$ when $q\geq2$.
Using Cattabriga's $L^q$ theory for the Stokes system, Serre \cite{Se}
extended this existence theory to $1<q<2$.  For more singular data, the
notion of very weak solutions in $L^q$ was introduced by Giga \cite{Giga},
Conca \cite{C}, and Galdi et al.\ \cite{GSS}; see also
\cite{AR,FGS,kim-h,MP}.  In particular, Kim \cite{kim-h} obtained a
solvability theory for very weak steady solutions on smooth bounded domains.
Pukhnachev \cite{Pukhnachev} studied a planar steady viscous problem
involving a point vortex at low Reynolds number.

The present problem has two features not covered by these general existence
theories.  First, the point-supported forcing in \eqref{force-condition} is
chosen precisely to produce the velocity singularity $b/r$, so that the
solution lies naturally in $L^q(B_1)$ for $q<2$ rather than in a standard
energy class.  Second, we seek a quantitative high-Reynolds-number
construction that identifies both the boundary layer and the selected
inviscid limit.  Thus the problem combines very weak solvability with a
singular perturbation analysis.

\subsection{Prandtl--Batchelor flows and steady boundary-layer expansions}

The rigorous analysis of Prandtl--Batchelor flows in circular geometries was initiated by Kim \cite{K1998,K2000,kim-formula}. For boundary rotations close to rigid rotation, the validity of the Prandtl expansion was subsequently established on the disk in \cite{FGLT} and on an annulus in \cite{FGLT2}.

These works are closely related to the development of steady Prandtl boundary-layer expansions for Navier--Stokes flows. The validity of such expansions was first obtained by Guo and Nguyen \cite{GN} for steady flows over a moving plate with a shear Euler background. Iyer then extended this circle of results to a rotating disk and developed global-in-\(x\) steady Prandtl layers over moving boundaries \cite{Iyer1,Iyer2}. For the classical no-slip boundary condition, Guo and Iyer proved the validity for shear Euler flows on a short \(x\)-interval and for related forced problems \cite{GI1,GI2}. Gao and Zhang further treated large \(L\) under a concavity assumption on the Prandtl profile and non-shear Euler flows for small \(L\) \cite{GZ,GZ2}. In a periodic setting, G\'erard-Varet and Maekawa established Sobolev stability of steady Prandtl expansions \cite{GMM}. Iyer and Masmoudi proved boundary-layer expansions and global-in-\(x\) stability for stationary Navier--Stokes flows \cite{IM1,IM2}. More recently, Gao and Xin \cite{GX} studied Prandtl boundary layers in an infinitely long convergent channel, giving another natural global model for the steady boundary-layer expansion; see also \cite{DIN2024,FPZ}. The Prandtl--Batchelor constructions in \cite{FGLT,FGLT2} use the same matched Euler--Prandtl expansion and remainder-stability strategy, with the additional feature that the leading Euler state is selected by the Batchelor--Wood relation rather than prescribed in advance.

\subsection{The inviscid profile and main results}

We prescribe a nearly rigid rotating velocity on the outer boundary:
\begin{align}\label{rotating-boundary-condition}
\mathbf{g}=(\alpha+\eta \varpi(\theta))\mathbf{t},
\qquad \alpha>0,
\end{align}
where $\eta$ is small and $\varpi$ is smooth and $2\pi$-periodic.  We first
construct the solution on the punctured disk
\[
 \Omega=[0,2\pi]\times(0,1]=B_1\setminus\{0\},
\]
and then identify its extension to $B_1$ in the very weak sense.  Writing
$\mathbf u^\varepsilon=u^\varepsilon\mathbf t+v^\varepsilon\mathbf n$, the
equations on $\Omega$ take the form
\begin{eqnarray}
\left \{
\begin {array}{lll}
u^\varepsilon u^\varepsilon_\theta+rv^\varepsilon u^\varepsilon_r+u^\varepsilon v^\varepsilon+p^\varepsilon_\theta-\varepsilon^2\big(ru^\varepsilon_{rr}
 + u^\varepsilon_{r}+\frac{u^\varepsilon_{\theta\theta}}{r}+\frac{2}{r} v^\varepsilon_\theta-\frac{ u^\varepsilon}{r}\big)=0,\\[5pt]
u^\varepsilon v^\varepsilon_\theta+rv^\varepsilon v^\varepsilon_r-(u^\varepsilon)^2+rp^\varepsilon_r-\varepsilon^2\big(rv^\varepsilon_{rr}
 +v^\varepsilon_{r}+\frac{ v^\varepsilon_{\theta\theta}}{r}-\frac{2}{r} u^\varepsilon_\theta-\frac{ v^\varepsilon}{r}\big)=0,\\[5pt]
 u^\varepsilon_\theta+rv^\varepsilon_r+ v^\varepsilon=0,
\end{array}
\right.\label{NS-curvilnear}
\end{eqnarray}
with
\begin{align}\label{rotating-boundary-condition-0}
u^\varepsilon(\theta,1)=\alpha+\eta \varpi(\theta),\qquad
v^\varepsilon(\theta,1)=0.
\end{align}
Formally, the inviscid outer equations are
\begin{align}\label{equ:Euler}
\left\{
\begin{aligned}
& u_e\partial_\theta u_e+v_er\partial_ru_e+u_ev_e+\partial_\theta p_e=0,\\
&u_e\partial_\theta v_e+v_er\partial_rv_e-(u_e)^2+r\partial_rp_e =0,\\
&\partial_\theta u_e+\partial_r(rv_e)=0.
\end{aligned}
\right.
\end{align}
The appropriate outer profile is the Couette--point-vortex flow
\begin{align}
(u_e,v_e)=\left(ar+\frac{b}{r},0\right),\qquad b>0,\quad a+b>0.
\label{Couette-flow}
\end{align}
Away from the origin, the term $b/r$ is the velocity induced by the point
vortex \eqref{point-vortex}; the regular part $ar$ has constant vorticity
$2a$.  In accordance with the annular Prandtl--Batchelor mechanism discussed
in \cite{FGLT2}, the background parameter is selected by the
Batchelor--Wood formula
\begin{align}
u_e^2(1)=(a+b)^2
=\alpha^2+\frac{\eta^2}{2\pi}
\int_0^{2\pi}\varpi^2(\theta)\,d\theta .
\label{Batchelor-Wood-formula}
\end{align}
Thus
\begin{align}
a=\left(\alpha^2+\frac{\eta^2}{2\pi}
\int_0^{2\pi}\varpi^2(\theta)\,d\theta\right)^{1/2}-b .
\label{re-Couette-flow}
\end{align}

The boundary mismatch is corrected by the leading Prandtl layer near $r=1$.
With $Y=(r-1)/\varepsilon$, let $(u_p^{(0)},v_p^{(1)})$ solve
\begin{align}
\left \{
\begin {array}{ll}
\big(u_e(1)+u_p^{(0)}\big)\partial_\theta u_p^{(0)}
+\big(v_p^{(1)}-v_p^{(1)}(\theta,0)\big)\partial_Yu_p^{(0)}
-\partial_{YY}u_p^{(0)}=0,\\[5pt]
\partial_\theta u_p^{(0)}+\partial_Yv_p^{(1)}=0,\\[5pt]
u_p^{(0)}(\theta,Y)=u_p^{(0)}(\theta+2\pi,Y),\quad
v_p^{(1)}(\theta,Y)=v_p^{(1)}(\theta+2\pi,Y),\\[5pt]
u_p^{(0)}\big|_{Y=0}=\alpha+\eta\varpi(\theta)-u_e(1),\quad
\displaystyle\lim_{Y\to-\infty}(u_p^{(0)},v_p^{(1)})=(0,0).
\label{eq:leading-prandtl}
\end{array}
\right.
\end{align}
This system is derived and solved in Section~2 as part of the matched
asymptotic construction.

\begin{Theorem}\label{main-theorem}
Assume that $\varpi$ is a smooth $2\pi$-periodic function satisfying
$\int_0^{2\pi}\varpi(\theta)\,d\theta=0$.  There exist
$\varepsilon_0>0$ and $\eta_0>0$ such that, for
$0<\varepsilon<\varepsilon_0$ and $0<\eta<\eta_0$, the system
\eqref{NS-curvilnear}--\eqref{rotating-boundary-condition-0} admits a solution
$(u^\varepsilon,v^\varepsilon)$ satisfying
\begin{align}
\left\|u^\varepsilon(\theta,r)-u_e(r)
-u_p^{(0)}\left(\theta,\frac{r-1}{\varepsilon}\right)
\right\|_{L^\infty(\Omega)}
&\leq C\varepsilon,\label{velocity-behavior-01}\\
\|v^\varepsilon\|_{L^\infty(\Omega)}
&\leq C\varepsilon.\label{velocity-behavior-02}
\end{align}
Moreover, if $\omega^\varepsilon$ denotes the vorticity on $\Omega$, then for
each $0<r_0<1$,
\begin{align}
\lim_{\varepsilon\to0}
\|\omega^\varepsilon-2a\|_{L^\infty(B_{r_0}\setminus\{0\})}=0.
\label{vorticity-behavior-0}
\end{align}
Here $a$ is given by \eqref{re-Couette-flow}.
\end{Theorem}

The singular part of the velocity is explicit.  Subtracting it leaves a
regular field at the origin, and hence the construction extends to the full
disk in the very weak formulation.

\begin{Definition}[Very weak formulation on the whole disk]\label{def:very-weak-whole-disk}
Let $1<q<2$ and let $q'=q/(q-1)>2$.  We say that
$\mathbf u^\varepsilon\in L^q(B_1)$ is a very weak solution of
\eqref{ns}, \eqref{force-condition}, and \eqref{rotating-boundary-condition}
if, for every
$\boldsymbol{\varphi}\in W^{2,q'}(B_1)\cap W^{1,q'}_0(B_1)$ satisfying
$\nabla\cdot\boldsymbol{\varphi}=0$ and
$\boldsymbol{\varphi}(0)=0$, the punctured-domain limit
\[
\lim_{\delta\to0^+}
\iint_{B_1\setminus B_\delta}
\mathbf u^\varepsilon\otimes\mathbf u^\varepsilon:
\nabla\boldsymbol{\varphi}\,dxdy
\]
exists and
\begin{align}\label{def:very-weak-identity}
-\varepsilon^2\iint_{B_1}
\mathbf u^\varepsilon\cdot\Delta\boldsymbol{\varphi}\,dxdy
-\lim_{\delta\to0^+}
\iint_{B_1\setminus B_\delta}
\mathbf u^\varepsilon\otimes\mathbf u^\varepsilon:
\nabla\boldsymbol{\varphi}\,dxdy
&=\varepsilon^2
\langle\mathbf F,\boldsymbol{\varphi}\rangle
-\varepsilon^2
\int_{\partial B_1}
\mathbf g\cdot\frac{\partial\boldsymbol{\varphi}}
{\partial\mathbf n}\,dl .
\end{align}
Moreover,
\[
\iint_{B_1}\mathbf u^\varepsilon\cdot\nabla\psi\,dxdy=0
\]
for every admissible scalar test function $\psi\in W^{1,q'}(B_1)$.
The nonlinear term in \eqref{def:very-weak-identity} is understood only
through the above punctured limit; this is essential because the singular
component $(b/r)\mathbf t$ is not square integrable at the origin.
\end{Definition}

\begin{Theorem}\label{main-theorem-0}
Under the hypotheses of Theorem~\ref{main-theorem}, the forced system
\eqref{ns}, \eqref{force-condition}, and \eqref{rotating-boundary-condition}
admits a velocity $\mathbf u^\varepsilon\in L^q(B_1)$ for every
$1\leq q<2$.  For every $1<q<2$, this velocity is a very weak solution
in the sense of Definition~\ref{def:very-weak-whole-disk}.  Set
\[
 \mathbf u_{\rm reg}^\varepsilon
 :=\mathbf u^\varepsilon-\frac{b}{r}\mathbf t,
 \qquad
 \omega_{\rm reg}^\varepsilon
 :=\operatorname{curl}\mathbf u_{\rm reg}^\varepsilon.
\]
Then $\mathbf u_{\rm reg}^\varepsilon\in H^1(B_1)$.  The tangential component in the following estimate is understood on $B_1\setminus\{0\}$ and extended to $r=0$ by the bounded regular part:
\begin{align}
\left\|\mathbf u_{\rm reg}^\varepsilon\cdot\mathbf{t}
-ar
-u_p^{(0)}\left(\theta,\frac{r-1}{\varepsilon}\right)
\right\|_{L^\infty(B_1)}
&\leq C\varepsilon, \label{velocity-behavior-1}\\
\|\mathbf{u}^\varepsilon\cdot\mathbf{n}\|_{L^\infty(B_1)}
&\leq C\varepsilon.\label{velocity-behavior-2}
\end{align}
Moreover,
\[
 \operatorname{curl}\mathbf u^\varepsilon
 =2\pi b\,\delta_0+\omega_{\rm reg}^\varepsilon
 \qquad\hbox{in }\mathcal D'(B_1).
\]
For every $0<r_0<1$,
\begin{align}
\lim_{\varepsilon\to0}
\|\omega_{\rm reg}^\varepsilon-2a\|_{L^\infty(B_{r_0})}=0.
\label{vorticity-behavior}
\end{align}
\end{Theorem}

\begin{Remark}
The zero-mean assumption on $\varpi$ is only a normalization: one may absorb
its mean into $\alpha$ by writing
\[
 \alpha+\eta\varpi(\theta)
 =\left(\alpha+\frac{\eta}{2\pi}
 \int_0^{2\pi}\varpi(\theta)\,d\theta\right)
 +\eta\widetilde{\varpi}(\theta),
 \qquad
 \int_0^{2\pi}\widetilde{\varpi}(\theta)\,d\theta=0.
\]
We impose smallness of $\eta$ in order to work in the perturbative regime
required for the boundary-layer construction.
\end{Remark}

\begin{Remark}[Shielding of the singular core]\label{rmk:shielding}
The estimates in Theorem~\ref{main-theorem-0} exhibit a shielding effect
in the following analytical sense.  Although the total velocity contains the
prescribed singular component $(b/r)\mathbf t$, the regular correction
obtained after separating this point-vortex part and subtracting the
Euler--Prandtl approximation remains uniformly bounded up to the origin.
Correspondingly, the vorticity splits distributionally as
\[
 \operatorname{curl}\mathbf u^\varepsilon
 =2\pi b\,\delta_0+\omega_{\rm reg}^\varepsilon,
 \qquad
 \omega_{\rm reg}^\varepsilon\to 2a
\]
locally uniformly in the interior.  Thus the point singularity does not
create an additional singular regular vorticity and does not destroy the
constant-vorticity selection in the surrounding eddy.  Here ``shielding''
refers to boundedness of the regular correction and to the persistence of
the selected regular background vorticity; the imposed point-vortex velocity
$b/r$ itself remains singular at the origin.
\end{Remark}

\subsection{Mathematical difficulties and the main ideas}

\begingroup
We briefly explain why the singular-core problem is not a routine
perturbation of the regular-disk construction in \cite{FGLT}.  In the
regular problem, the leading Euler flow is the smooth rigid rotation
\((ar,0)\), so the velocity field, the vorticity, and the linearized error
operator remain regular at the origin.  Here the outer profile contains the
point-vortex component \(b/r\), whose distributional vorticity is
\(2\pi b\delta_0\).  Although the Prandtl layer is still localized near the
outer boundary, the stability estimates must remain uniform up to the vortex
center.  This also explains why the result cannot be obtained by simply
letting the inner radius in the annular construction of \cite{FGLT2} tend to
zero: the annular estimates are not uniform in that limit and do not yield a
whole-disk very weak solution with the point-supported forcing
\eqref{force-condition}.

The proof has five main components.

\smallskip
\noindent
{\bf (i) Separation of the prescribed singular core.}
The point-vortex velocity \((b/r)\mathbf t\) belongs to \(L^q(B_1)\) for
\(q<2\), but is not bounded at the origin.  Thus the full solution cannot be
constructed in the regular Sobolev framework used for nonsingular disk
flows.  We write
\[
 \mathbf u^\varepsilon=\frac{b}{r}\mathbf t+
 \mathbf u_{\rm reg}^\varepsilon,
\]
and prove that the regular part remains controlled across \(r=0\).  This
singular/regular decomposition is also what makes the whole-disk very weak
formulation in Definition~\ref{def:very-weak-whole-disk} natural.

\smallskip
\noindent
{\bf (ii) High-order Euler--Prandtl approximation.}
The approximate velocity has the schematic form
\[
 u^a(\theta,r)
 = ar+\frac{b}{r}
   +\chi(r)u_p^{(0)}\left(\theta,\frac{r-1}{\varepsilon}\right)
   +\hbox{higher-order terms},
 \qquad
 v^a(\theta,r)=\hbox{regular terms}.
\]
The singular part is prescribed by the forcing, while the regular component
\(ar\) is selected by the Batchelor--Wood relation
\eqref{Batchelor-Wood-formula}.  The expansion is carried to high order so
that the residuals are sufficiently small for the nonlinear stability
argument.  The full expansion and the residual estimates are given in
Section~2.

\smallskip
\noindent
{\bf (iii) Failure of the regular and annular positivity mechanisms.}
For the regular disk problem, the stability estimate can be closed with
multipliers adapted to the smooth core.  In the present problem, the
linearized transport also contains the singular component \(b/r\).  Formally,
testing by angular derivatives produces the non-coercive singular expression
\[
 \int_0^1\int_0^{2\pi}
 \frac{b}{r}\bigl(u_\theta^2+v_\theta^2+2u_\theta v\bigr)
 \,d\theta\,dr,
\]
which cannot be controlled uniformly near \(r=0\) by the regular-disk
argument.  Weighted estimates of the type used on an annulus give, at best,
bounds schematically like
\[
 \left\|r^{1/2+\delta}
 \left(u^\varepsilon-u_e-u_p^{(0)}\right)\right\|_{L^\infty}
 +\left\|r^{1/2+\delta}v^\varepsilon\right\|_{L^\infty}
 \lesssim \varepsilon,
\]
which allow growth as \(r\to0\).  The central task is therefore to obtain an
unweighted estimate that is uniform at the vortex center.

\smallskip
\noindent
{\bf (iv) Logarithmic geometry and mode-sensitive stability.}
The key device is the logarithmic radial coordinate
\[
 s=-\log r.
\]
It maps the punctured disk to the half-cylinder
\([0,2\pi]\times[0,\infty)\), sends the outer boundary to \(s=0\), and sends
the vortex center to \(s=\infty\).  In this geometry, uniform control near
the singular core becomes a far-field stability problem in \(s\).  The
solution is then decomposed into the angular mean, the first Fourier mode,
and the higher-mode component:
\[
 u=u_0+c(s)\sin\theta+d(s)\cos\theta+\widetilde u,
 \qquad
 v=e(s)\sin\theta+f(s)\cos\theta+\widetilde v.
\]
The first Fourier mode is the delicate part.  The positivity estimate gives
direct control of \(e'\) and \(f'\), but not of all first derivatives of
\(c\) and \(d\).  A key cancellation, using the divergence-free relation
\(u_\theta-v_s+v=0\), allows the singular interaction represented
schematically by
\[
 \iint e^{-s}v^a u_{1s}u_{1\theta}\,d\theta\,ds
\]
to be controlled by the available first-mode quantities.  This
mode-specific cancellation is the main new stability mechanism for the
singular-core problem.  The detailed estimates are proved in Section~3.

\smallskip
\noindent
{\bf (v) Nonlinear closure and vorticity selection.}
The linear stability bounds are combined with second-order estimates and a
Sobolev embedding on the half-cylinder to obtain the \(L^\infty\) control
needed for the nonlinear error system.  The first Fourier mode again
requires special treatment, because second-order integrations by parts
produce boundary terms at \(s=0\); these are handled through local Stokes and
trace estimates near the outer boundary.  This yields a contraction argument
for the nonlinear remainder in Section~4.

Finally, after the prescribed Dirac vorticity is subtracted, the regular
vorticity is governed on the half-cylinder by an equation whose leading
transport-diffusion part has the form
\[
 \varepsilon^2(\partial_s^2+\partial_\theta^2)\omega
 -b\,\partial_\theta\omega=\Lambda_b .
\]
The angular mean is controlled by an ordinary differential equation, while
the nonzero modes are estimated through the vorticity equation together with
the velocity bounds.  Transforming back to the disk gives
\[
 \operatorname{curl}\mathbf u^\varepsilon
 =2\pi b\,\delta_0+\omega_{\rm reg}^\varepsilon,
 \qquad
 \omega_{\rm reg}^\varepsilon\to 2a
\]
locally uniformly in the interior.  Thus the imposed point-vortex mass
remains isolated, while the surrounding regular vorticity is still selected
by the Batchelor--Wood relation.

The construction is first carried out on the punctured disk.  After
subtracting the prescribed point-vortex velocity, the regular part
\[
 \left(u^\varepsilon-\frac{b}{r},\,v^\varepsilon\right)
\]
extends across \(r=0\) as an \(H^1(B_1)\) vector field.  Reintroducing the
singular component \((b/r)\mathbf t\) and testing against admissible
divergence-free test fields gives the whole-disk very weak solution asserted
in Theorem~\ref{main-theorem-0}.
\endgroup
\subsection{Organization of the paper}

Section~2 constructs the high-order approximate solution by matched
asymptotic expansion and establishes its residual bounds.  Section~3 derives
the error equations in logarithmic radial coordinates and proves the basic
energy and positivity estimates, including the separate first-mode analysis.
Section~4 derives the higher-order estimates needed for $L^\infty$ control
and proves the existence of a solution to the nonlinear error system.
Section~5 establishes the vorticity estimates and completes the proofs of
Theorems~\ref{main-theorem} and \ref{main-theorem-0}.

\medskip
\noindent
{\bf Notation.}
Throughout the paper, $\|\cdot\|_p$ denotes the $L^p$ norm for
$2\leq p\leq\infty$, and $C,C_1,\ldots$ denote positive constants independent
of sufficiently small $\varepsilon$ and $\eta$.  We write
$A_1\lesssim A_2$ when $A_1\leq CA_2$, and similarly for $\gtrsim$.
Finally,
$\iint f$ denotes
$\int_0^\infty\int_0^{2\pi}f(\theta,s)\,d\theta\,ds$.

\section{Construction of approximate solutions}
\indent

In this section, we use the method of multi-scale matched asymptotic expansion to construct an approximate solution of the Navier--Stokes equations (\ref{NS-curvilnear}).

Before the logarithmic change of variables, all estimates in this section are written in the physical polar variables $(\theta,r)$.  Thus, when an $L^2$ estimate for a function of $(\theta,r)$ is used, the integration is with respect to $d\theta dr$; in the residual estimates below we write the relevant integrals explicitly.

\subsection{Euler expansions for the leading-order and first-order terms}

\indent

Away from the boundary, we make the following formal expansions
\begin{align*}
&u^{\varepsilon}(\theta,r)=u_e^{(0)}(\theta,r)+\varepsilon u_e^{(1)}(\theta,r)+\cdots,\\[5pt]
&v^{\varepsilon}(\theta,r)=v_e^{(0)}(\theta,r)+\varepsilon v_e^{(1)}(\theta,r)+\cdots,\\[5pt]
&p^{\varepsilon}(\theta,r)=p_e^{(0)}(\theta,r)+\varepsilon p_e^{(1)}(\theta,r)+\cdots.
\end{align*}

\subsubsection{Equations for $(u_e^{(0)},v_e^{(0)},p_e^{(0)})$}

\indent

By substituting the above expansions into (\ref{NS-curvilnear}) and collecting the zeroth-order (leading-order) terms, we deduce that $(u_e^{(0)},v_e^{(0)},p_e^{(0)})$ satisfies the following steady nonlinear Euler equations
\begin{eqnarray}
\left \{
\begin {array}{ll}
u_e^{(0)} \partial_\theta u_e^{(0)}+rv_e^{(0)} \partial_ru_e^{(0)}+u_e^{(0)} v_e^{(0)}+\partial_\theta p_e^{(0)}=0,\\[7pt]
 u_e^{(0)} \partial_\theta v_e^{(0)}+rv_e^{(0)} \partial_rv_e^{(0)}-(u_e^{(0)})^2+r\partial_rp_e^{(0)}=0,\\[7pt]
 \partial_\theta u_e^{(0)}+r\partial_rv_e^{(0)}+v_e^{(0)}=0.
\end{array}
\right.\label{outer-leading-order-equation}
\end{eqnarray}

We choose the leading-order Euler flow $(u_e^{(0)},v_e^{(0)})$ as in \eqref{Couette-flow}; then
(\ref{outer-leading-order-equation}) reduces to
 \begin{align}
\partial_\theta p_e^{(0)}(\theta,r)=0,\quad  \partial_rp_e^{(0)}(\theta,r)=\frac{1}{r}u_e^2.\nonumber
 \end{align}
 Thus we deduce that
 \begin{align}
 p_e^{(0)}(\theta,r)=p_e(r),\quad  p'_e(r)=\frac{(ar+b/r)^2}{r}.\label{outer-leading-order-pressure}
 \end{align}

\subsubsection{Equations for $(u_e^{(1)},v_e^{(1)},p_e^{(1)})$}

\indent

By collecting the $\varepsilon$-order terms, we deduce that $(u_e^{(1)},v_e^{(1)},p_e^{(1)})$ satisfies the following linearized Euler equations in $\Omega$
\begin{eqnarray}
\left \{
\begin {array}{ll}
u_e(r) \partial_\theta u_e^{(1)}+ru'_e(r)v_e^{(1)}+u_e(r)v^{(1)}_e +\partial_\theta p_e^{(1)}=0,\\[5pt]
u_e(r) \partial_\theta v_e^{(1)}-2u_e(r)u_e^{(1)}+r\partial_rp_e^{(1)}=0,\\[5pt]
\partial_\theta u_e^{(1)}+r\partial_rv_e^{(1)}+ v_e^{(1)}=0,
\end{array}
\right.\label{outer-1-order-equation}
\end{eqnarray}
which are equipped with the boundary conditions
\begin{align}
 v_e^{(1)}|_{r=1}=-v_p^{(1)}|_{Y=0},\quad v_e^{(1)}(\theta,r)=v_e^{(1)}(\theta+2\pi,r) ,\label{outer-1-order-bc}
\end{align}
where $v_p^{(1)}$ is the solution of Prandtl equations which will be derived in the next subsection.

Before we perform higher-order expansions we need to modify $(u_e^{(1)},v_e^{(1)},p_e^{(1)})$, e.g., \eqref{modify-Euler}. The details will be given later.
\subsection{Prandtl expansions near $r=1$ for the $\varepsilon^{-1}$-, leading-, and first-order terms}
\indent

We  introduce the scaled variable  $Y=\frac{r-1}{\varepsilon}\in (-\infty,0]$ and make the following Prandtl expansions near $r=1$
\begin{align}\label{first-order-expansion}
\begin{aligned}
&u^\varepsilon=u_e(r)+u_p^{(0)}(\theta,Y)+\varepsilon\big[u_e^{(1)}(\theta,r)+u_p^{(1)}(\theta,Y)\big]
+\cdots,\\[5pt]
&v^\varepsilon=v_p^{(0)}(\theta,Y)+\varepsilon\big[v_e^{(1)}(\theta,r)+v_p^{(1)}(\theta,Y)\big]
+\varepsilon^2\big[v_e^{(2)}(\theta,r)+v_p^{(2)}(\theta,Y)\big]+\cdots,\\[5pt]
&p^\varepsilon=p_e(r)+p_p^{(0)}(\theta,Y)+\varepsilon\big[p_e^{(1)}(\theta,r)+p_p^{(1)}(\theta,Y)\big]
+\varepsilon\big[p_e^{(2)}(\theta,r)+p_p^{(2)}(\theta,Y)\big]+\cdots,
\end{aligned}
\end{align}
where as $Y\rightarrow -\infty$
\begin{align}
  \partial_\theta^l\partial_Y^mv_p^{(i)}(\theta,Y)\rightarrow 0, \ \partial_\theta^l\partial_Y^mp_p^{(i)}(\theta,Y)\rightarrow 0,\label{matching-condition}
\end{align}
here $l,m\geq 0, i=0,1,\cdots,$
which satisfy the following boundary conditions
\begin{align*}
&u_e^{(0)}(\theta,1)+u_p^{(0)}(\theta,0)=\alpha+\eta \varpi(\theta), \quad u_e^{(i)}(\theta,1)+u_p^{(i)}(\theta,0)=0,\ i\geq 1,\\
&v_e^{(i)}(\theta,1)+v_p^{(i)}(\theta,0)=0, \ i\geq 0.
\end{align*}
The boundary conditions of $u_p^{(i)}(\theta,Y)$ as $Y\rightarrow-\infty $ will be given later.

\subsubsection{Equations for $(v_p^{(0)}, p_p^{(0)})$}

\indent

By substituting the above expansions into (\ref{NS-curvilnear}) and collecting the $\frac{1}{\varepsilon}$-order terms, we get
\begin{align*}
\partial_Yv_p^{(0)}(\theta,Y)=0,\ \ \partial_Yp_p^{(0)}(\theta,Y)=0,
\end{align*}
which together with (\ref{matching-condition}) imply
\begin{align*}
v_p^{(0)}=0,\ \ p_p^{(0)}=0.
\end{align*}

\subsubsection{Equations for $(u_p^{(0)},v_p^{(1)},p_p^{(1)})$}

\indent

By substituting the above expansions into (\ref{NS-curvilnear}) and collecting the zeroth-order (leading-order) terms, we obtain
 the following steady Prandtl equations for $(u_p^{(0)},v_p^{(1)})$
\begin{eqnarray}
\left \{
\begin {array}{ll}
\big(u_e(1)+u_p^{(0)}\big)\partial_\theta u_p^{(0)}+\big(v_e^{(1)}(\theta,1)+ v_p^{(1)}\big)\partial_Yu_p^{(0)}-\partial_{YY}u_p^{(0)}=0,\\[5pt]
\partial_\theta u_p^{(0)}+\partial_Yv_p^{(1)}=0,\\[5pt]
u_p^{(0)}(\theta,Y)=u_p^{(0)}(\theta+2\pi,Y),\quad v_p^{(1)}(\theta,Y)=v_p^{(1)}(\theta+2\pi,Y),\\[5pt]
u_p^{(0)}\big|_{Y=0}=\alpha+\eta \varpi(\theta)-u_e(1),\ \ \lim\limits_{Y\rightarrow -\infty}(u_p^{(0)},v_p^{(1)})=(0,0)
\label{prandtl-problem-1}
\end{array}
\right.
\end{eqnarray}
and the pressure $p_p^{(1)}$ satisfies
\begin{align}\label{equation-of-first-pressure}
\partial_Yp_p^{(1)}(\theta, Y)=(u_p^{(0)})^2(\theta,Y)+2u_e(1)u_p^{(0)}(\theta,Y), \ \lim_{Y\rightarrow -\infty}p_p^{(1)}(\theta,Y)=0.
\end{align}

\subsubsection{Equations for $(u_p^{(1)},v_p^{(2)})$}

\indent

By substituting the above expansions into the first and third equation in (\ref{NS-curvilnear}) and collecting the $\varepsilon$-order terms, we obtain
 the following linearized steady Prandtl equations for $(u_p^{(1)},v_p^{(2)})$
\begin{eqnarray}
\left \{
\begin {array}{ll}
\big(u_e(1)+u_p^{(0)}\big)\partial_\theta u_p^{(1)}+\big(v_e^{(1)}(\theta,1)+ v_p^{(1)}\big)\partial_Yu_p^{(1)}+(v_e^{(2)}(\theta,1)+v_p^{(2)})\partial_{Y}u_p^{(0)}\\[5pt]
\quad \quad \quad \quad +(u_p^{(1)}+u_e^{(1)}(\theta,1))\partial_\theta u_p^{(0)}-\partial_{YY}u_p^{(1)}=f_1(\theta,Y),\\[5pt]
\partial_\theta u_p^{(1)}+\partial_Yv_p^{(2)}+\partial_Y(Yv_p^{(1)})=0,\\[5pt]
u_p^{(1)}(\theta,Y)=u_p^{(1)}(\theta+2\pi,Y),\quad v_p^{(2)}(\theta,Y)=v_p^{(2)}(\theta+2\pi,Y),\\[5pt]
u_p^{(1)}\big|_{Y=0}=-u_e^{(1)}\big|_{r=1},\ \ \lim\limits_{Y\rightarrow -\infty}(\partial_Yu_p^{(1)},v_p^{(2)})=(0,0),
\label{first-linearized-prandtl-problem-near-1}
\end{array}
\right.
\end{eqnarray}
where
\begin{align*}
f_1(\theta,Y)=&-\partial_\theta p_p^{(1)}+Y\partial_{YY}u_p^{(0)}+\partial_Yu_p^{(0)}\\[5pt]
&-u_p^{(0)}\big(\partial_\theta u_e^{(1)}(\theta,1)+v_e^{(1)}(\theta,1)+v_p^{(1)}\big)-u'_e(1)Y\partial_\theta u_p^{(0)}\\[5pt]
&-(\partial_rv_e^{(1)}(\theta,1)+v_e^{(1)}(\theta,1))Y\partial_Y u_p^{(0)}-(u_e'(1)+Y\partial_Yu_p^{(0)}+u_e(1))v_p^{(1)}.
\end{align*}

Since the Prandtl expansion equations depend on the Euler expansions, we will present the higher-order Prandtl equations later.

\subsection{Solvability of the leading-order and first-order Prandtl equations and the first-order Euler equations}
\indent

 We have the following result which comes from  Corollary 2.4 in \cite{FGLT}.
\begin{Proposition} There exists $\eta_0>0$ such that for any $\eta\in(0,\eta_0)$, the equations (\ref{prandtl-problem-1}) have a unique solution $(u_p^{(0)},v_p^{(1)})$ which satisfies
\begin{align}\label{decay-behavior-prandtl}
\sum_{j+k\leq m}\int_{-\infty}^0\int_0^{2\pi}\Big|\partial_\theta^j\partial_Y^k (u_p^{(0)},v_p^{(1)})\Big|^2\big<Y\big>^{2l}d\theta dY\leq C(m,l)\eta^2, \ m,l \geq 0.
\end{align}
\end{Proposition}

Notice that
$$v^{(1)}_p(\theta,Y)=\int_{-\infty}^Y\partial_Yv^{(1)}_{p}(\theta,Y')dY'=-\int_{-\infty}^Y\partial_\theta u^{(0)}_{p}(\theta,Y')dY',$$
we have
\begin{align*}
\int_{0}^{2\pi}v^{(1)}_p(\theta, Y)d\theta=0, \ \forall\ Y\leq 0.
\end{align*}
Finally, solving (\ref{equation-of-first-pressure}), we obtain $p_p^{(1)}(\theta,Y)$  which decays very fast as $Y\rightarrow -\infty$.

Completely similar to Proposition 2.6 in \cite{FGLT} we get
\begin{Proposition}\label{solvability-of-first-Euler}
The linearized Euler equations (\ref{outer-1-order-equation}) have a solution $(u_e^{(1)}, v_e^{(1)}, p_e^{(1)}))$ which satisfies
\begin{align}
|\partial_\theta u_e^{(1)}+v_e^{(1)}|(\theta,r)\leq& C\eta r, \ |\partial_\theta v_e^{(1)}-u_e^{(1)}|(\theta,r)\leq C\eta r,  \ \forall (\theta, r)\in \Omega, \label{Estimate-of-first-combined-linearized-Euler-equation}\\[5pt]
 \left(\int_0^1\int_0^{2\pi}\left|\partial^k_\theta\partial^j_r(u_e^{(1)},v_e^{(1)})\right|^2d\theta dr\right)^{\frac12}\leq& C(k,j)\eta, \quad \forall j,k\geq 0,\label{Estimate-of-first-linearized-Euler-equation}\\[5pt]
 r^2\triangle u_e^{(1)}-u_e^{(1)}+2\partial_\theta& v_e^{(1)}=0,\quad \int_{0}^{2\pi}v_e^{(1)}(\theta,r)d\theta=0,\label{identity-for-first-Euler}
\end{align}
here and below, $\Delta=\partial_{rr}+\frac{\partial_r}{r}+\frac{\partial_{\theta\theta}}{r^2}.$ More concretely, if  we write
\begin{align*}
v_p^{(1)}(\theta,0)=\sum_{n=1}^{+\infty}\big(a_{n1} \cos(n\theta)+b_{n1}\sin(n\theta)\big),
\end{align*}
then
\begin{align}
u_e^{(1)}(\theta, r)&=\sum_{n=1}^{+\infty}\big(a_{n1} r^{n-1} \sin(n\theta)-b_{n1}r^{n-1}\cos(n\theta)\big),\label{1orderEuler-1}\\
v_e^{(1)}(\theta,r)&=-\sum_{n=1}^{+\infty}\big(a_{n1} r^{n-1} \cos(n\theta)+b_{n1}r^{n-1}\sin(n\theta)\big).\label{1orderEuler-2}
\end{align}

\end{Proposition}

The proof of solving the linearized Prandtl equations (\ref{first-linearized-prandtl-problem-near-1}) is similar to the proof of Proposition 2.7 in \cite{FGLT}, so we only show the following result and omit its proof here.

\begin{Proposition}\label{decay-estimates-of-linearized-Prandtl} There exists $\eta_0>0$ such that for any $\eta\in(0,\eta_0)$, the equations (\ref{first-linearized-prandtl-problem-near-1}) have a unique solution $(u_p^{(1)},v_p^{(2)})$ which satisfies
\begin{align}
&\sum_{j+k\leq m}\int_{-\infty}^0\int_0^{2\pi}\big|\partial_\theta^j\partial_Y^k \big(u_p^{(1)}-A_{1\infty},v_p^{(2)}\big)\big|^2\big<Y\big>^{2l}d\theta dY\leq C(m,l)\eta^2, \ \ m, l\geq 0; \nonumber\\
& \int_0^{2\pi}v_p^{(2)}(\theta,Y)d\theta =0, \ \forall \ Y\leq 0,  \label{decay-behavior-prandtl-1}
\end{align}
where
$A_{1\infty}:=\lim\limits_{Y\rightarrow -\infty}u_p^{(1)}(\theta,Y)$ is a constant which satisfies $|A_{1\infty}|\leq C\eta.$
\end{Proposition}

Next, we construct the pressure $p_p^{(2)}(\theta,Y)$. For that we define $\tilde{u}_p^{(1)}=u_p^{(1)}-A_{1\infty}$  and  consider the equation
\begin{align}\label{equation-for-second-pressure}
\partial_Yp_p^{(2)}(\theta, Y)=g_1(\theta,Y), \quad \lim_{Y\rightarrow -\infty}p_p^{(2)}(\theta,Y)=0,
\end{align}
where
\begin{align*}
g_1(\theta,Y)=&-Y\partial_Yp_p^{(1)}+\partial_{YY}v_p^{(1)}-u_e(1)\partial_\theta v_p^{(1)}-u_p^{(0)} (\partial_\theta v_e^{(1)}(\theta,1)+\partial_\theta v_p^{(1)})\\[5pt]
&-\partial_Yv_p^{(1)}(v_e^{(1)}(\theta,1)+v_p^{(1)})-2(Yu'_e(1)u_p^{(0)}+u_e(1)\tilde{u}_p^{(1)}+[u_e^{(1)}(\theta,1)+A_1]u_p^{(0)}+u_p^{(0)}\tilde{u}_p^{(1)}
\end{align*}
can be obtained by replacing $u_p^{(1)}$ by $\tilde{u}_p^{(1)}$ in the expansion (\ref{first-order-expansion}) and putting the new expansion into the second equation of (\ref{NS-curvilnear}), then collecting the $\varepsilon$-order terms together.
 Notice that $g_1(\theta,Y)$ decays fast as $Y\rightarrow -\infty$, we can get $p_p^{(2)}(\theta,Y)$ by solving (\ref{equation-for-second-pressure}) and deduce that $p_p^{(2)}(\theta,Y)$ decays fast as $Y\rightarrow -\infty$.

\subsection{Linearized Euler equations for $(u_e^{(2)}, v_e^{(2)}, p_e^{(2)})$ and their solvabilities}
\indent

We first modify $(u_e^{(1)}, v_e^{(1)}, p_e^{(1)})$. Let $\chi(r)\in C^\infty([0,1])$ be an increasing smooth function such that
\begin{align}\label{cut-off-function}
\chi(r)=
\left\{
\begin{array}{lll}
0, \quad r\in [0,\frac12], \\[5pt]
1, \quad r\in [\frac34,1].
\end{array}
\right.
\end{align}
Then, let $\phi_1(r)=-A_{1\infty}\big(r\chi''(r)+\chi'(r)-\frac{\chi(r)}{r}\big)$ and
\begin{align*}
A_1(r):=a_1r+r\int_0^r\frac{\phi_1(s)}{2s}-\frac{1}{r}\int_0^r\frac{s\phi_1(s)}{2}ds,
\end{align*}
where $a_1$ is a constant such that $A_1(1)=0$. Obviously, $|a_1|\leq C\eta$ and
\begin{align}\label{corrector-of-first-order-Euler-equation}
\left\{
\begin{array}{ll}
rA''_1(r)+A'_1(r)-\frac{A_1(r)}{r}=\phi_1(r),\ 0<r\leq 1 \\[5pt]
A_1(1)=0.
 \end{array}
 \right.
\end{align}
Direct computation gives $\|\partial_r^kA_1(r)\|_\infty\leq C(k)\eta.$  Moreover, notice that $\chi(r)=0$ for $r\leq \frac12$, we deduce that $A_1(r)=a_1r$ for $r\leq \frac12$.

Set
\begin{align}\label{modify-Euler}
\begin{aligned}
\tilde{u}_e^{(1)}(\theta,r):&=u_e^{(1)}(\theta,r)+\chi(r)A_{1\infty}+A_1(r), \\[5pt]
\tilde{v}_e^{(1)}(\theta,r):&=v_e^{(1)}(\theta,r),\\
\tilde{p}_e^{(1)}(\theta,r):&=p_e^{(1)}(\theta,r)+\int_{0}^r\frac{2u_e(s)}{s}\big(\chi(s)A_{1\infty}+A_1(s)\big)ds,
\end{aligned}
\end{align}
then $(\tilde{u}_e^{(1)},\tilde{v}_e^{(1)},\tilde{p}_e^{(1)})$ also satisfies the linearized Euler equations (\ref{outer-1-order-equation})
with the boundary conditions (\ref{outer-1-order-bc}). Moreover, there holds
\begin{align}\label{Estimate-of-modified-first-linearized-Euler-equation}
\begin{aligned}
|\partial_\theta \tilde{u}_e^{(1)}+\tilde{v}_e^{(1)}|(\theta,r)\leq& C\eta r, \ |\partial_\theta \tilde{v}_e^{(1)}-\tilde{u}_e^{(1)}|(\theta,r)\leq C\eta r, \ \forall (\theta, r)\in \Omega, \\[5pt]
 \left(\int_0^1\int_0^{2\pi}\left|\partial^k_\theta\partial^j_r(\tilde{u}_e^{(1)},\tilde{v}_e^{(1)})\right|^2d\theta dr\right)^{\frac12}\leq& C(k,j)\eta, \quad \forall j,k\geq 0;\\
  r^2\triangle \tilde{u}_e^{(1)}-\tilde{u}_e^{(1)}+2\partial_\theta& \tilde{v}_e^{(1)}=0, \quad
  \int_{0}^{2\pi}\tilde{v}_e^{(1)}d\theta=0.
\end{aligned}
\end{align}

Putting
\begin{align*}
&u^{\varepsilon}(\theta,r)=u_e(r)+\varepsilon \tilde{u}_e^{(1)}(\theta,r)+\varepsilon^2 u_e^{(2)}(\theta,r)+\cdots,\\[5pt]
&v^{\varepsilon}(\theta,r)=\varepsilon \tilde{v}_e^{(1)}(\theta,r)+\varepsilon^2 v_e^{(2)}(\theta,r)+\cdots,\\[5pt]
&p^{\varepsilon}(\theta,r)=p_e(r)+\varepsilon \tilde{p}_e^{(1)}(\theta,r)+\varepsilon^2 p_e^{(2)}(\theta,r)+\cdots
\end{align*}
into the Navier--Stokes equations (\ref{NS-curvilnear}), we obtain the following linearized Euler equations for $(u_e^{(2)},v_e^{(2)}, p_e^{(2)})$
\begin{eqnarray}
\left \{
\begin {array}{ll}
u_e(r) \partial_\theta u_e^{(2)}+rv_e^{(2)}u'_e(r)+u_e(r)v^{(2)}_e+\partial_\theta p_e^{(2)}+\tilde{u}_e^{(1)}\partial_\theta \tilde{u}_e^{(1)}+\tilde{v}_e^{(1)}r\partial_r\tilde{u}_e^{(1)}+\tilde{u}_e^{(1)}\tilde{v}_e^{(1)}=0,\\[5pt]
u_e(r)\partial_\theta v_e^{(2)}-2u_e(r)u_e^{(2)}+r\partial_rp_e^{(2)}+\tilde{u}_e^{(1)}\partial_\theta \tilde{v}_e^{(1)}+\tilde{v}_e^{(1)}r\partial_r \tilde{v}_e^{(1)}-(\tilde{u}_e^{(1)})^2=0,\\[7pt]
\partial_\theta u_e^{(2)}+r\partial_rv_e^{(2)}+ v_e^{(2)}=0,
\end{array}
\right.\label{outer-2-order-equation}
\end{eqnarray}
 with the boundary conditions
\begin{align}\label{boundary-condition-of-second-Euler}
 v_e^{(2)}|_{r=1}=-v_p^{(2)}|_{Y=0}, \quad v_e^{(2)}(\theta,r)=v_e^{(2)}(\theta+2\pi,r).
\end{align}

\begin{Proposition}\label{solvability-of-second-Euler-equation}
The linearized Euler equations  (\ref{outer-2-order-equation}) have a solution $(u_e^{(2)}, v_e^{(2)}, p_e^{(2)})$ which satisfies
\begin{align*}
|\partial_\theta u_e^{(2)}+v_e^{(2)}|(\theta,r)\leq& C\eta r, \ |\partial_\theta v_e^{(2)}-u_e^{(2)}|(\theta,r)\leq C\eta r,  \ \forall (\theta, r)\in \Omega, \\[5pt]
 \left(\int_0^1\int_0^{2\pi}\left|\partial^k_\theta\partial^j_r(u_e^{(2)},v_e^{(2)})\right|^2d\theta dr\right)^{\frac12}\leq& C(k,j)\eta, \quad \forall j,k\geq 0,\\
  r^2\triangle u_e^{(2)}-u_e^{(2)}+2\partial_\theta &v_e^{(2)}=0, \quad \int_{0}^{2\pi}v_e^{(2)}d\theta=0.
\end{align*}
\end{Proposition}
\begin{proof}
 Eliminating the pressure $p_e^{(2)}$ in the equation (\ref{outer-2-order-equation}), we obtain
\begin{align*}
-ru_e(r)\triangle(rv_e^{(2)})-\tilde{u}_e^{(1)}r\triangle(r\tilde{v}_e^{(1)})+\tilde{v}_e^{(1)}(r^2\triangle \tilde{u}_e^{(1)}-\tilde{u}_e^{(1)}+2\partial_\theta \tilde{v}_e^{(1)})=0.
\end{align*}
Recall that $\triangle(rv_e^{(1)})=0$ and using (\ref{Estimate-of-modified-first-linearized-Euler-equation}), we obtain the following equation for $rv_e^{(2)}$ in $\Omega$
\begin{align*}
\left\{\begin{array}{lll}
-\triangle(rv_e^{(2)})=0,\\[5pt]
rv_e^{(2)}|_{r=1}=-v_p^{(2)}(\theta,0).
\end{array}
\right.
\end{align*}
Then, we can complete the proof by following the argument of Proposition \ref{solvability-of-first-Euler} line by line; we omit the repetitive details.
\end{proof}

\subsection{Linearized Prandtl equations for  $(u_p^{(2)},v_p^{(3)})$ and their solvabilities}
\indent

Putting the expansion
\begin{align*}
&u^\varepsilon(\theta,r)=u_e(r)+u_p^{(0)}(\theta,Y)+\varepsilon\big[\tilde{u}_e^{(1)}(\theta,r)+\tilde{u}_p^{(1)}(\theta,Y)\big]
+\varepsilon^2\big[u_e^{(2)}(\theta,r)+u_p^{(2)}(\theta,Y)\big]+\cdots,\\[5pt]
&v^\varepsilon(\theta,r)=\varepsilon\big[\tilde{v}_e^{(1)}(\theta,r)+v_p^{(1)}(\theta,Y)\big]
+\varepsilon^2\big[v_e^{(2)}(\theta,r)+v_p^{(2)}(\theta,Y)\big]+\varepsilon^3[v_e^{(3)}(\theta,r)+v_p^{(3)}(\theta,Y)]+\cdots,\\[5pt]
&p^\varepsilon(\theta,r)=p_e(r)+\varepsilon\big[\tilde{p}_e^{(1)}(\theta,r)+p_p^{(1)}(\theta,Y)\big]
+\varepsilon^2\big[p_e^{(2)}(\theta,r)+p_p^{(2)}(\theta,Y)\big]+\varepsilon^3p_p^{(3)}(\theta,Y)+\cdots
\end{align*}
with the boundary conditions
\begin{align*}
 u_e^{(2)}(\theta,1)+u_p^{(2)}(\theta,0)=0,\ v_e^{(3)}(\theta,1)+v_p^{(3)}(\theta,0)=0, \ \lim_{Y\rightarrow -\infty}(\partial_Yu_p^{(2)},v_p^{(3)})=(0,0)
\end{align*}
into the first and third equation of (\ref{NS-curvilnear}), collecting $\varepsilon^2$-order terms together, we obtain the following linearized steady Prandtl equations for $(u_p^{(2)},v_p^{(3)})$
\begin{eqnarray}
\left \{
\begin {array}{ll}
\big(u_e(1)+u_p^{(0)}\big)\partial_\theta u_p^{(2)}+\big(v_e^{(1)}(\theta,1)+ v_p^{(1)}\big)\partial_Yu_p^{(2)}+(v_e^{(3)}(\theta,1)+v_p^{(3)})\partial_{Y}u_p^{(0)}\\[5pt]
\quad \quad \quad \quad +(u_p^{(2)}+u_e^{(2)}(\theta,1))\partial_\theta u_p^{(0)}-\partial_{YY}u_p^{(2)}=f_2(\theta,Y),\\[5pt]
\partial_\theta u_p^{(2)}+\partial_Yv_p^{(3)}+\partial_Y(Yv_p^{(2)})=0,\\[5pt]
u_p^{(2)}(\theta,Y)=u_p^{(2)}(\theta+2\pi,Y),\quad v_p^{(3)}(\theta,Y)=v_p^{(3)}(\theta+2\pi,Y),\\[5pt]
u_p^{(2)}\big|_{Y=0}=-u_e^{(2)}\big|_{r=1},\ \ \lim\limits_{Y\rightarrow -\infty}(\partial_Yu_p^{(2)},v_p^{(3)})=(0,0),
\label{second-linearized-prandtl-problem-near-1}
\end{array}
\right.
\end{eqnarray}
where
{\small\begin{align*}
f_2(\theta,Y)=&-\partial_\theta p_p^{(2)}+Y\partial_{YY}u_p^{(1)}+\partial_Yu_p^{(1)}+\partial_{\theta\theta}u_p^{(0)}-u_p^{(0)}
-\tilde{u}_p^{(1)}\partial_\theta u_p^{(1)}-v_p^{(2)}\partial_Y u_p^{(1)}-\sum_{i+j=2}v_p^{(i)}Y\partial_Y u_p^{(j)}\\[5pt]
&-\sum_{k=0}^2\sum_{i+j=2-k, (k,j)\neq (0,2)}\Big(\frac{\partial_r^k\tilde{u}_e^{(i)}(\theta,1)}{k!}Y^k \partial_\theta \tilde{u}_p^{(j)}+\tilde{u}_p^{(j)}\frac{\partial_r^k\partial_\theta \tilde{u}_e^{(i)}(\theta,1)}{k!}Y^k\Big)+u_e^{(2)}(\theta,1)\partial_\theta u_p^{(0)}\\[5pt]
&-\sum_{k=0}^1\sum_{i+j=2-k}\Big(\frac{\partial_r^k\tilde{v}_e^{(i)}(\theta,1)}{k!}Y^{k+1}\partial_Y \tilde{u}_p^{(j)}+v_p^{(i)}\frac{\partial_r^k(r\partial_r \tilde{u}_e^{(j)})(\theta,1)}{k!}Y^k\Big)\\[5pt]
&-\sum_{k=0}^2\sum_{i+j=3-k, (k,j)\neq (0,2),(0,0)}\frac{\partial_r^k\tilde{v}_e^{(i)}(\theta,1)}{k!}Y^k \partial_Y \tilde{u}_p^{(j)}
\end{align*}}
with $\tilde{u}_p^{(0)}=u_p^{(0)}, \ \tilde{u}_e^{(0)}=u_e(r),\ \tilde{u}_e^{(2)}=u_e^{(2)}, \ \tilde{v}_e^{(2)}=v_e^{(2)}$ in the above expression.

\begin{Proposition}\label{decay-estimates-of-higher-order-linearized-Prandtl} There exists $\eta_0>0$ such that for any $\eta\in(0,\eta_0)$, equations (\ref{second-linearized-prandtl-problem-near-1}) have a unique solution $(u_p^{(2)},v_p^{(3)})$ which satisfies
\begin{align}\label{decay-behavior-prandtl-2}
&\sum_{j+k\leq m}\int_{-\infty}^0\int_0^{2\pi}\Big|\partial_\theta^j\partial_Y^k (u_p^{(2)}-A_{2\infty},v_p^{(3)} )\Big|^2\big<Y\big>^{2l}d\theta dY\leq C(m,l)\eta^2, \ \ m, l\geq 0, \nonumber\\
&\int_0^{2\pi}v_p^{(3)}(\theta, Y)d\theta=0, \ \forall\ Y\leq 0,
\end{align}
where $A_{2\infty}:=\lim\limits_{Y\rightarrow -\infty}u_p^{(2)}(\theta,Y)$ is a constant which satisfies $|A_{2\infty}|\leq C\eta$.

\end{Proposition}
The proof follows the argument of Proposition~\ref{decay-estimates-of-linearized-Prandtl}. Since $f_2(\theta,Y)$ decays rapidly as $Y\rightarrow -\infty$, we omit the repetitive details.

We construct the pressure $p_p^{(3)}(\theta,Y)$ by considering the equation
\begin{align}\label{equation-of-second-pressure}
\partial_Yp_p^{(3)}(\theta, Y)=g_2(\theta,Y), \quad \lim_{Y\rightarrow -\infty}p_p^{(3)}(\theta,Y)=0,
\end{align}
where
{\small\begin{align*}
g_2(\theta,Y)=&\partial_{YY}v_p^{(2)}+Y\partial_{YY}v_p^{(1)}+\partial_Yv_p^{(1)}-2\partial_\theta u_p^{(0)}-Y\partial_Yp_p^{(2)}-\sum_{i+j=3} v_p^{(i)}\partial_Yv_p^{(j)}\\[5pt]
&-\sum_{i+j=2}\Big(\tilde{u}_p^{(i)}\partial_\theta v_p^{(j)}+v_p^{(i)}Y\partial_Yv_p^{(j)}-\tilde{u}_p^{(i)}\tilde{u}_p^{(j)}+\tilde{v}_e^{(i)}(\theta,1)Y\partial_Yv_p^{(j)}
+v_p^{(i)}\partial_r\tilde{v}_e^{(j)}(\theta,1)\Big)\\[5pt]
&-\sum_{k=0}^1\sum_{i+j=2-k}\Big(\frac{\partial_r^k\tilde{u}_e^{(i)}(\theta,1)}{k!}Y^k\partial_\theta v_p^{(j)}+\frac{\partial_r^k\partial_\theta \tilde{v}_e^{(j)}(\theta,1)}{k!}Y^k\partial_\theta \tilde{u}_p^{(i)}\Big)\\[5pt]
&-\sum_{k=0}^2\sum_{i+j=2-k}\Big(\frac{\partial_r^k\tilde{u}_e^{(i)}(\theta,1)}{k!}Y^k\tilde{u}_p^{(j)}+\frac{\partial_r^k \tilde{u}_e^{(j)}(\theta,1)}{k!}Y^k \tilde{u}_p^{(i)}\Big),
\end{align*}}
in which $\tilde{u}_p^{(0)}=u_p^{(0)},  \ \tilde{u}_e^{(0)}=u_e(r),\ \tilde{u}_e^{(2)}=u_e^{(2)}+A_{2\infty}, \ \tilde{v}_e^{(2)}=v_e^{(2)}$, and here and below $\tilde{u}_p^{(2)}=u_p^{(2)}-A_{2\infty}$.  $g_2(\theta,Y)$ can be derived by the same argument as $g_1(\theta,Y)$.
Moreover, since $g_2(\theta,Y)$ decays fast as $Y\rightarrow -\infty$, we obtain $p_p^{(3)}$ by solving \eqref{equation-of-second-pressure}, and $p_p^{(3)}$ also decays fast as $Y\rightarrow -\infty$.

By induction, we can perform higher-order Euler and Prandtl expansions.
The order in which we solve the equations is as follows
\begin{align*}
&(u_e(r),0)\rightarrow (u_p^{(0)},v_p^{(1)},p_p^{(1)})\rightarrow (u_e^{(1)},v_e^{(1)},p_e^{(1)})\rightarrow(u_p^{(1)},v_p^{(2)},p_p^{(2)})\\&\rightarrow (\tilde{u}_e^{(1)},\tilde{v}_e^{(1)},\tilde{p}_e^{(1)})\rightarrow (u_e^{(2)},v_e^{(2)},p_e^{(2)})\rightarrow(u_p^{(2)},v_p^{(3)},p_p^{(3)})\rightarrow\cdot\cdot\cdot\\&\rightarrow (\tilde{u}_e^{(i)},\tilde{v}_e^{(i)},\tilde{p}_e^{(i)})\rightarrow (u_e^{(i+1)},v_e^{(i+1)},p_e^{(i+1)})\rightarrow(u_p^{(i+1)},v_p^{(i+2)},p_p^{(i+2)})\\&\rightarrow\cdot\cdot\cdot,
\end{align*}
where $(\tilde{u}_e^{(1)},\tilde{v}_e^{(1)},\tilde{p}_e^{(1)})$ and  $(\tilde{u}_e^{(i)},\tilde{v}_e^{(i)},\tilde{p}_e^{(i)}) $ are the modifications of $(u_e^{(1)},v_e^{(1)},p_e^{(1)})$ and  $(u_e^{(i)},v_e^{(i)},p_e^{(i)})$, respectively. We omit the repetitive algebraic details and list only the properties used below. One can refer to  Section 2.3 in \cite{FGLT}.

\begin{Proposition}\label{prop:uniform-higher-euler}
There exists $\eta_0>0$ such that, for any $\eta\in(0,\eta_0)$, the higher-order Euler and Prandtl correctors may be chosen so that, for $i\geq2$,

\begin{align}\label{Estimate-of-modified-third-linearized-Euler-equation-i}
|\partial_\theta \tilde{u}_e^{(i)}+\tilde{v}_e^{(i)}|(\theta,r)\leq& C\eta r, \ |\partial_\theta \tilde{v}_e^{(i)}-\tilde{u}_e^{(i)}|(\theta,r)\leq C\eta r,  \ \forall (\theta, r)\in \Omega,\nonumber\\[5pt]
 \left(\int_0^1\int_0^{2\pi}\left|\partial^k_\theta\partial^j_r(\tilde{u}_e^{(i)},\tilde{v}_e^{(i)})\right|^2d\theta dr\right)^{\frac12}\leq& C(k,j)\eta, \quad \forall j,k\geq 0;\\
 r^2\triangle \tilde{u}_e^{(i)} -\tilde{u}_e^{(i)}+2\partial_\theta &\tilde{v}_e^{(i)}=0,\quad \int_{0}^{2\pi}\tilde{v}^{(i)}_ed\theta=0\nonumber
\end{align}and for $i\geq3$
\begin{align}\label{decay-behavior-prandtl-i}
&\sum_{j+k\leq m}\int_{-\infty}^0\int_0^{2\pi}\big|\partial_\theta^j\partial_Y^k (u_p^{(i)}-A_{i\infty},v_p^{(i+1)} )\big|^2\big<Y\big>^{2l}d\theta dY\leq C(m,l)\eta^2, \ \ m, l\geq 0, \nonumber\\
& \int_0^{2\pi}v_p^{(i+1)}(\theta,Y)d\theta=0, \ \forall \ Y\leq 0,
\end{align}
where $A_{i\infty}:=\lim\limits_{Y\rightarrow -\infty} u_p^{(i)}$ is a constant which satisfies $|A_{i\infty}|\leq C\eta$.
In particular, no higher Euler correction introduces an additional
$r^{-1}$ singularity at the center; the only singular velocity in the
approximate solution is the prescribed leading term $b/r$.
\end{Proposition}

\subsection{Approximate solutions}

\indent

In this subsection, we construct an approximate solution of the Navier--Stokes equations (\ref{NS-curvilnear}).
Set
\begin{align*}
\tilde{u}_p^a(\theta,r)&=\chi(r)\Big(u_p^{(0)}(\theta,Y)+\sum_{i=1}^{11}\varepsilon^i\tilde{u}_p^{(i)}(\theta,Y)+
\varepsilon^{12}u_p^{(12)}(\theta,Y)\Big)\\&
=:\chi(r)u_p^a,\ \ \tilde{u}_p^{(i)}=u_p^{(i)}-A_{i\infty},\\[5pt]
\tilde{v}_p^a(\theta,r)&=\chi(r)\Big(\sum_{i=1}^{13}\varepsilon^i v_p^{(i)}(\theta,Y)\Big)
=:\chi(r)v_p^a,\\
\tilde{p}_p^a(\theta,r)&=\chi^2(r)\Big(\sum_{i=1}^{13}\varepsilon^i p_p^{(i)}(\theta,Y)\Big)
=:\chi^2(r)p_p^a,
\end{align*}
and
\begin{align*}
u_e^a(\theta,r)=&u_e(r)+\sum_{i=1}^{12}\varepsilon^i \tilde{u}_e^{(i)}(\theta,r), \\[5pt]
 v_e^a(\theta,r)=&\sum_{i=1}^{12}\varepsilon^i \tilde{v}_e^{(i)}(\theta,r),\\[5pt]
 p_e^a(\theta,r)=&p_e(r)+\sum_{i=1}^{12}\varepsilon^i \tilde{p}_e^{(i)}(\theta,r).
\end{align*}

We construct an approximate solution
\begin{align}\label{approximate-solution}
\begin{aligned}
u^a(\theta,r)&=u_e^a(\theta,r)+\tilde{u}_p^a(\theta,r)+\varepsilon^{13}h(\theta,r),\\[5pt]
v^a(\theta,r)&=v_e^a(\theta,r)+\tilde{v}_p^a(\theta,r),\\[5pt]
p^a(\theta,r)&=p_e^a(\theta,r)+\tilde{p}_p^a(\theta,r),
\end{aligned}
\end{align}
where the corrector $h(\theta,r)$ will be given in Appendix~A which satisfies
\beno
h(\theta, 1)=0, \ \left(\int_0^1\int_0^{2\pi}\left|\partial_\theta^j\partial_r^kh\right|^2d\theta dr\right)^{\frac12}\leq C(j,k)\varepsilon^{-k}
\eeno
and makes $(u^a, v^a)$ be divergence-free
\beno
u^a_\theta+rv^a_r+v^a=0.
\eeno
Moreover, $(u^a, v^a)$  satisfies the following boundary conditions
\begin{align*}
 u^a(\theta+2\pi,r)&=u^a(\theta,r), \ v^a(\theta+2\pi,r)=v^a(\theta,r),\\[5pt]
u^a(\theta, 1)&=\alpha+\eta \varpi(\theta), \quad \ \ v^a(\theta,1)=0.
\end{align*}
And $v^a$ satisfies
\begin{align}\label{va-null-0}
\int_{0}^{2\pi}v^ad\theta=0.
\end{align}

Notice that $\chi(r)=0$ and $A_i(r)=a_ir$ for $r\leq \frac12$, by collecting the estimates (\ref{decay-behavior-prandtl-1}), (\ref{Estimate-of-modified-first-linearized-Euler-equation}),(\ref{decay-behavior-prandtl-2}), \eqref{Estimate-of-modified-third-linearized-Euler-equation-i} and \eqref{decay-behavior-prandtl-i},  we deduce that
\begin{align}\label{estimate-on-Euler-parts}
\begin{aligned}
&\|u^a_e-u_e(r)\|_{\infty}+\|(\partial_r(u^a_e-u_e(r)),\partial_\theta u^a_e)\|_{\infty}\leq C\varepsilon\eta, \ \ \|v^a_e\|_{\infty}+\|(\partial_rv^a_e, \partial_\theta v^a_e)\|_{\infty}\leq C\varepsilon\eta, \\[5pt]
&|\partial_\theta u_e^a(\theta,r)+v_e^a(\theta,r)|\leq C\varepsilon \eta r, \ \ |\partial_\theta v_e^a(\theta,r)-u_e^a(\theta,r)|\leq C\varepsilon \eta r, \ \forall (\theta, r)\in \Omega,
\end{aligned}
\end{align}
and
\begin{align}\label{estimate-on-Prandtl-parts}
\|(Y^j\partial_Y^ku_p^a, Y^j\partial_\theta^ku_p^a)\|_{\infty}\leq C \eta, \ \ \|(Y^j\partial_Y^kv_p^a,Y^j\partial_\theta^kv_p^a)\|_{\infty}\leq C \varepsilon\eta,  \ \ \forall j\leq 2, k\leq 1.
\end{align}

 Finally, set
 \begin{align*}
 R_u^a&=u^au^a_\theta+v^aru^a_r+u^av^a+p^a_\theta-\varepsilon^2\Big( ru^a_{rr}
+\frac{u^a_{\theta\theta}}{r}+2\frac{v^a_{\theta}}{r}+u^a_r-\frac{u^a}{r}\Big),\\[5pt]
R_v^a&=u^av^a_\theta+v^arv^a_r-(u^a)^2+rp^a_r-\varepsilon^2 \Big(rv^a_{rr}+\frac{v^a_{\theta\theta}}{r}-2\frac{u^a_{\theta}}{r}+v^a_r-\frac{v^a}{r}\Big),
 \end{align*}
 then there hold
\begin{align}\label{remainder-estimate}
\Big(\int_0^1\int_0^{2\pi}(R_u^a)^2(\theta,r)d\theta dr\Big)^{\frac12}
+\Big(\int_0^1\int_0^{2\pi}(R_v^a)^2(\theta,r)d\theta dr\Big)^{\frac12}
\leq C\varepsilon^{13}
 \end{align}

and
\begin{align}\label{remainder-r-derivative-estimate}
\Big(\int_0^1\int_0^{2\pi}(\partial_r R_u^a)^2(\theta,r)d\theta dr\Big)^{\frac12}
+\Big(\int_0^1\int_0^{2\pi}(\partial_r R_v^a)^2(\theta,r)d\theta dr\Big)^{\frac12}
\leq C\varepsilon^{12}.
\end{align}

 and
{\small\begin{align*}
\left\{
\begin{array}{lll}
u^au^a_\theta+v^aru^a_r+u^av^a+p^a_{\theta}-\varepsilon^2\big( ru^a_{rr}
+\frac{u^a_{\theta\theta}}{r}+2\frac{v^a_{\theta}}{r}+u^a_r-\frac{u^a}{r}\big)=R_u^a,\ (\theta,r)\in \Omega,\\[5pt]
u^av^a_\theta+v^arv^a_r-(u^a)^2+rp^a_r-\varepsilon^2\big( rv^a_{rr}+\frac{v^a_{\theta\theta}}{r}-2\frac{u^a_{\theta}}{r}+v^a_r-\frac{v^a}{r}\big)=R_v^a,\ (\theta,r)\in \Omega,\\[5pt]
 u^a_\theta+(rv^a)_r=0, \ (\theta,r)\in \Omega, \\[5pt]
 u^a(\theta+2\pi,r)=u^a(\theta,r), \ v^a(\theta+2\pi,r)=v^a(\theta,r), \ (\theta,r)\in \Omega,\\[5pt]
u^a(\theta,1)=\alpha+\varpi(\theta)\eta,\ v^a(\theta,1)=0, \ \theta\in [0,2\pi].
\end{array}
\right.
\end{align*}}

In fact, when $r\leq \frac12$, there hold
\begin{align*}
R_u^a=&u_e^a\partial_\theta u_e^a+v_e^ar\partial_ru_e^a+u^a_ev^a_e+\partial_\theta p_e^a-\varepsilon^2\Big( r\partial_{rr}u^a_e
+\frac{\partial_{\theta\theta}u^a_e}{r}+2\frac{\partial_\theta v^a_e}{r}+\partial_ru^a_e-\frac{u_e^a}{r}\Big)\\
=&\sum_{m=13}^{24}\varepsilon^m\Big( \sum_{i+j=m, i,j\leq 12} \tilde{u}_e^{(i)}[\partial_\theta \tilde{u}_e^{(j)}+\tilde{v}_e^{(j)}]+\sum_{i+j=m, i,j\leq 12} \tilde{v}_e^{(i)}r\partial_r \tilde{u}_e^{(j)}\Big),\\
R_v^a=&u_e^a\partial_\theta v_e^a+v_e^ar\partial_rv_e^a-(u^a_e)^2+r\partial_r p_e^a-\varepsilon^2\Big( r\partial_{rr}v^a_e
+\frac{\partial_{\theta\theta}v^a_e}{r}+2\frac{\partial_\theta u^a_e}{r}+\partial_rv^a_e-\frac{v_e^a}{r}\Big)\\
=&\sum_{m=13}^{24}\varepsilon^m\Big( \sum_{i+j=m, i,j\leq 12} \tilde{u}_e^{(i)}[\partial_\theta \tilde{v}_e^{(j)}-\tilde{u}_e^{(j)}]+\sum_{i+j=m, i,j\leq 12} \tilde{v}_e^{(i)}r\partial_r \tilde{v}_e^{(j)}\Big).
\end{align*}
Using \eqref{Estimate-of-modified-first-linearized-Euler-equation} and \eqref{Estimate-of-modified-third-linearized-Euler-equation-i}, we deduce that
\begin{align}
|R_u^a(\theta,r)|+|R_v^a(\theta,r)|\leq C\varepsilon^{13} r, \ |\partial_rR_u^a(\theta,r)|+|\partial_rR_v^a(\theta,r)|\leq C\varepsilon^{12}, \  r<\frac12.\label{remainder-decay}
\end{align}
\begingroup
Moreover, on the annular region $\frac12\leq r<1$, there holds
\begin{align*}
&\left(\int_{1/2}^1\int_0^{2\pi}(R_u^a)^2(\theta,r)d\theta dr\right)^{\frac12}
+\left(\int_{1/2}^1\int_0^{2\pi}(R_v^a)^2(\theta,r)d\theta dr\right)^{\frac12}
\leq C\varepsilon^{13},\\
&\left(\int_{1/2}^1\int_0^{2\pi}(\partial_rR_u^a)^2(\theta,r)d\theta dr\right)^{\frac12}
+\left(\int_{1/2}^1\int_0^{2\pi}(\partial_rR_v^a)^2(\theta,r)d\theta dr\right)^{\frac12}
\leq C\varepsilon^{12}.
\end{align*}
\endgroup
Together with \eqref{remainder-decay}, these estimates yield
\eqref{remainder-estimate} and \eqref{remainder-r-derivative-estimate}.

We also record the corresponding residual estimates in the logarithmic coordinate
$s=-\log r$, which are the bounds used in the fixed-point and vorticity
estimates below.  With the $L^2$ norm taken on
$[0,2\pi]\times[0,+\infty)$ in $(\theta,s)$,
\begingroup
\begin{align}\label{eq:log-residual-bounds}
\|e^{-s/2}R_u^a\|_2+
\|e^{-s/2}R_v^a\|_2
&\le C\varepsilon^{13},\\
\|\partial_s R_u^a\|_2+
\|\partial_s R_v^a\|_2+
\|e^{-s}\partial_s R_u^a\|_2+
\|e^{-s}\partial_s R_v^a\|_2
&\le C\varepsilon^{12},\\
\|\partial_\theta R_u^a\|_2+
\|\partial_\theta R_v^a\|_2
&\le C\varepsilon^{13}.
\end{align}
\endgroup
Indeed, $\partial_s=-r\partial_r$ and $dr=e^{-s}ds$; hence the estimates
above follow directly from \eqref{remainder-estimate},
\eqref{remainder-r-derivative-estimate}, and the same construction applied
after one angular differentiation.

\section{Linear stability estimates of error equations}

\indent
In this section, we derive the error equations and establish linear stability estimates.

\subsection{Error equations}
\indent

Set the error
\begin{align*}
u=u^\varepsilon-u^a,\ v=v^\varepsilon-v^a,\ p=p^\varepsilon-p^a,
\end{align*}
then there hold
{\small \begin{align}
\left\{
\begin{array}{lll}
u^au_\theta+uu^a_\theta+uu_\theta+v^aru_r+vru^a_r+vru_r+v^au+vu^a+vu+p_\theta-\varepsilon^2\big( ru_{rr}
+\frac{u_{\theta\theta}}{r}+2\frac{v_{\theta}}{r}+u_r-\frac{u}{r}\big)=R_u^a,\\[5pt]
u^av_\theta+uv^a_\theta+uv_\theta+v^arv_r+vrv^a_r+vrv_r-(u^2+2u u^a)+rp_r-\varepsilon^2 \big( rv_{rr}+\frac{v_{\theta\theta}}{r}-2\frac{u_{\theta}}{r}+v_r-\frac{v}{r}\big)=R_v^a,\\[5pt]
u_\theta+(rv)_r=0,  \\[5pt]
u(\theta+2\pi,r)=u(\theta,r), \ v(\theta+2\pi,r)=v(\theta,r), \\[5pt]
u(\theta,1)=0,\ v(\theta,1)=0.\nonumber
\end{array}
\right.
\end{align}}

Let
\begin{align*}
S_u&=u^au_\theta+v^aru_r+uu^a_\theta+vru^a_r+v^au+vu^a,\\[5pt]
S_v&=u^av_\theta+v^arv_r+uv^a_\theta+vrv^a_r-2uu^a,
\end{align*}
then the error equations become
\begin{align}\label{e:error-equation}
\left\{
\begin{array}{lll}
-\varepsilon^2\big(r u_{rr}
+\frac{u_{\theta\theta}}{r}+2\frac{v_{\theta}}{r}+u_r-\frac{u}{r}\big)+p_\theta+S_u=R_u,\\[5pt]
-\varepsilon^2\big( rv_{rr}+\frac{v_{\theta\theta}}{r}-2\frac{u_{\theta}}{r}+v_r-\frac{v}{r}\big)+rp_r+S_v=R_v,\\[5pt]
u_\theta+(rv)_r=0,  \\[5pt]
u(\theta,r)=u(\theta+2\pi,r),\  v(\theta,r)=v(\theta+2\pi,r), \\[5pt]
u(\theta,1)=0,\ v(\theta,1)=0,
 \end{array}
\right.
\end{align}
where
\begin{align*}
R_u=R_u^a-uu_\theta-vru_r-vu,\quad R_v=R_v^a-uv_\theta-vrv_r+u^2.
\end{align*}

\subsection{Linear stability estimates}
\indent

In this subsection, we consider the linear equations
\begin{align}\label{linear-equation-of-error}
\left\{
\begin{array}{lll}
-\varepsilon^2\big(r u_{rr}
+\frac{u_{\theta\theta}}{r}+2\frac{v_{\theta}}{r}+u_r-\frac{u}{r}\big)+p_\theta+S_u=F_u,\\[5pt]
-\varepsilon^2\big( rv_{rr}+\frac{v_{\theta\theta}}{r}-2\frac{u_{\theta}}{r}+v_r-\frac{v}{r}\big)+rp_r+S_v=F_v,\\[5pt]
u_\theta+(rv)_r=0,  \\[5pt]
u(\theta,r)=u(\theta+2\pi,r),\  v(\theta,r)=v(\theta+2\pi,r), \\[5pt]
u(\theta,1)=0,\ v(\theta,1)=0,
 \end{array}
\right.
\end{align}
 and establish its stability estimates. Due to the divergence-free condition and the boundary condition of $v$, we deduce that
\begin{align}
\int_0^{2\pi}v(\theta,r)d\theta=0\label{v-ava}
\end{align}
which implies that  the Poincar\'{e} inequality is valid:
\begin{align}\label{poin}
\int_0^{2\pi}v^2d\theta\leq\int_0^{2\pi}v_\theta^2d\theta.
\end{align}

We also give the following Hardy-type inequalities which will be used frequently.
\begin{Lemma}\label{H1}
\indent

(i)If $f$ is a function defined in $[\frac{1}{2},1]$ and $f(\frac{1}{2})=f(1)=0$,
 then there holds
\begin{align}
\int_{\frac{1}{2}}^1\frac{f^2(r)}{(1-r)^2}dr\leq C\int_{\frac{1}{2}}^1\big(f'(r)\big)^2dr.\label{weight-Hardy}
\end{align}

(ii)If $f$ is a function defined in $[0,+\infty)$ and $f(0)=0$,  then there holds
\begin{align}
&\int_0^{+\infty} t^\beta f^2(t)dt\leq C(\beta)\int_0^{+\infty} t^{\beta+2}\big(f'(t)\big)^2dt,\ \forall \beta <-1.\label{Hardy-3}
\end{align}

\end{Lemma}

\begin{proof}

(i) We only prove \eqref{weight-Hardy} in the case that  $f\in C^1$ and $f(1)=0, f(\frac{1}{2})=0$. By integrating by parts, one has
\begin{align*}
\int_{\frac{1}{2}}^1\frac{f^2}{(1-r)^2}dr=\int_{\frac{1}{2}}^1f^2d\frac{1}{1-r}=-\int_{\frac{1}{2}}^1\frac{1}{1-r}d(f^2)=-\int_{\frac{1}{2}}^12\frac{ff'}{1-r}dr.
\end{align*}
So we deduce that
\begin{align*}
\int_{\frac{1}{2}}^1\frac{f^2}{(1-r)^2}dr\leq 2\Big(\int_{\frac{1}{2}}^1\frac{f^2}{(1-r)^2} dr\Big)^{\frac12}\Big(\int_{\frac{1}{2}}^1(f')^2dr\Big)^{\frac12}.
\end{align*}
Thus \eqref{weight-Hardy} holds.

(ii) Recalling the Hardy's inequality in \cite{Stein}, for  $1\leq p<\infty,$ $0<q<\infty$ we have
\begin{equation*}
  \bigg(\int^\infty_0\bigg(\int^x_0|g(t)|dt\bigg)^px^{-q-1}dx\bigg)^{\frac{1}{p}}\leq\frac{p}{q}\bigg(\int^\infty_0|g(t)|^pt^{p-q-1}dt\bigg)^{\frac{1}{p}}.
\end{equation*}
Choosing the $g=f'$ and $p=2$, one obtains that
\begin{equation*}
  \bigg(\int^{\infty}_0f^2x^{-q-1}dx\bigg)^{\frac{1}{2}}\leq\bigg(\int^{\infty}_0|f'|^2t^{1-q}dt\bigg)^{\frac{1}{2}}.
\end{equation*}
Thus we get \eqref{Hardy-3}.

Consequently we complete the proof of this lemma.
\end{proof}

\subsubsection{Basic energy estimate}
\indent
\par
Introduce the new variable
\begin{equation*}
  s=-\ln r\in[0,+\infty).
\end{equation*}
Then $r=e^{-s}, dr=-e^{-s}ds$ and $(\theta,s)\in \mathcal C:=[0,2\pi]\times[0,
+\infty)$. Due to the divergence-free condition, we have
\begin{equation*}
  \partial_\theta u(\theta,s)=\partial_sv(\theta,s)-v,
\end{equation*}
where $u(\theta,s)=u(\theta,r(s))$ and $v(\theta,s)=v(\theta,r(s))$.
\par
By a direct calculation, we can rewrite the system \eqref{linear-equation-of-error} as follows
\begin{align}\label{rewrite-linear-equation-of-error}
\left\{
\begin{array}{lll}
-\varepsilon^2\big(u_{ss}
+u_{\theta\theta}+2v_{\theta}-u\big)+e^{-s}p_\theta+e^{-s}S_u=e^{-s}F_u,\\[5pt]
-\varepsilon^2\big(v_{ss}+v_{\theta\theta}-2u_{\theta}-v\big)-e^{-s}p_s+e^{-s}S_v=e^{-s}F_v,\\[5pt]
u_\theta-v_s+v=0,  \\[5pt]
u(\theta,s)=u(\theta+2\pi,s),\  v(\theta,s)=v(\theta+2\pi,s), \\[5pt]
u(\theta,0)=0,\ v(\theta,0)=0,
 \end{array}
\right.
\end{align}
where $(\theta,s)\in \mathcal{C}$.
\par
We decompose $(u,v)$ as follows
\begin{align}
u(\theta,s)&=u_0+c(s)\sin\theta+d(s)\cos\theta+\tilde{u}\nonumber
\\&=:u_0+u_1+\tilde{u},\label{decompose-u}\\
v(\theta,s)&=e(s)\sin\theta+f(s)\cos\theta+\tilde{v}\nonumber
\\&=:v_1+\tilde{v},\label{decompose-u2}
\end{align}
where $u_0=\frac{1}{2\pi}\int^{2\pi}_0u(\theta,s)d\theta$ and
 \begin{align*}
 c(s)&=\frac{1}{\pi}\int^{2\pi}_0u(\theta,s)\sin \theta d\theta,\ d(s)=\frac{1}{\pi}\int^{2\pi}_0u(\theta,s)\cos \theta d\theta,\\
e(s)&=\frac{1}{\pi}\int^{2\pi}_0v(\theta,s)\sin \theta d\theta, \ f(s)=\frac{1}{\pi}\int^{2\pi}_0v(\theta,s)\cos \theta d\theta.
 \end{align*}

It is easy to check that
\begin{align*}
\int_0^{2\pi}\tilde{u}d\theta=\int_0^{2\pi}\tilde{v}d\theta=0
\end{align*}
and
\begin{align*}
&\int_0^{2\pi}\partial_\theta^{m_1}\partial_s^{l_1}u_j\partial_\theta^{m_2}\partial_s^{l_2}\tilde{u}d\theta =\int_0^{2\pi}\partial_\theta^{m_1}\partial_s^{l_1}v_1\partial_\theta^{m_2}\partial_s^{l_2}\tilde{u}d\theta  \\&=\int_0^{2\pi}\partial_\theta^{m_1}\partial_s^{l_1}u_j\partial_\theta^{m_2}\partial_s^{l_2}\tilde{v}d\theta  =\int_0^{2\pi}\partial_\theta^{m_1}\partial_s^{l_1}v_1\partial_\theta^{m_2}\partial_s^{l_2}\tilde{v}d\theta =0,
\end{align*}
where $j=0,1, 0\leq m_1,l_1,m_2,l_2\leq2.$

Invoked by the behaviors   in \eqref{1orderEuler-1}, \eqref{1orderEuler-2} and \eqref{modify-Euler} as $r\rightarrow 0$,  we assume that
\begin{equation*}
  \lim_{s\rightarrow\infty}c(s)=c_{\infty},~~~ \lim_{s\rightarrow\infty}d(s)=d_{\infty},~~~\lim_{s\rightarrow\infty}e(s)=e_{\infty},~~~ \lim_{s\rightarrow\infty}f(s)=f_{\infty}
\end{equation*}
and define
\begin{equation*}
  u_{1\infty}=c_{\infty}\sin\theta+d_{\infty}\cos\theta,~~~~v_{1\infty}=e_{\infty}\sin\theta+f_{\infty}\cos\theta.
\end{equation*}
Thanks to divergence-free condition in \eqref{rewrite-linear-equation-of-error}, we obtain
\begin{align*}
&c(s)\cos\theta-d\sin\theta+\tilde{u}_\theta=(f'-f)\cos\theta+(e'-e)\sin\theta+\tilde{v}_s-\tilde{v}.
\end{align*}
Multiplying the above equality by $\cos \theta$ and $\sin\theta$ respectively, and integrating in $\theta\in [0,2\pi]$ imply
\begin{align}\label{divergence-free-condition-1}
\begin{aligned}
&c(s)=f'(s)-f(s), \\[5pt]
&d(s)=e(s)-e'(s).
\end{aligned}
\end{align}
We then arrive at
\begin{align}
c_{\infty}=-f_{\infty},~~~~d_{\infty}=e_{\infty},~~~\partial_\theta u_{1\infty}=-v_{1\infty},~~~\partial_\theta v_{1\infty}=u_{1\infty}.\label{relation-inf}
\end{align}

The following argument will need the following hypothesis

\[
{\bf (H)}:\quad
u_0\in H^1\big((0,+\infty)\big),\quad
u_1-u_{1\infty},\ v_1-v_{1\infty}\in H^2(\mathcal C),\quad
\tilde{u},\tilde{v},\tilde{u}_{\theta},\tilde{v}_{\theta}\in H^1(\mathcal C).
\]

\begin{Lemma}\label{hardy-for-c-d-e-f}
If $(u,v)$ satisfies  the hypothesis {\bf (H)}, then there holds
\begin{align*}
	&\int^{+\infty}_0\big(u_1^2+v_1^2+u_{1\theta}^2+v_{1\theta}^2+v_{1s}^2\big)e^{-s}ds\leq C\int^{+\infty}_0(e'^2+f'^2)ds.
\end{align*}
\end{Lemma}
\begin{proof} We only need to verify
\begin{align}
\int^{+\infty}_0(c^2+d^2+e^2+f^2)e^{-s}ds\leq C\int^{+\infty}_0(e'^2+f'^2)ds.\label{hardy-for-c-d-e-f-0}
\end{align}

In fact,
by using \eqref{divergence-free-condition-1} and Lemma \ref{H1}, we obtain
\begin{align*}
&\int^{+\infty}_0(c^2+d^2+e^2+f^2)e^{-s}ds\\
&\lesssim\int^{+\infty}_0(e^2+f^2)e^{-s}ds+\int^{+\infty}_0(e'^2+f'^2)e^{-s}ds\\
&\lesssim\int^{+\infty}_0\big({\frac{e^2}{s^2}}+{\frac{f^2}{s^2}}\big)e^{-s}s^2ds+\int^{+\infty}_0(e'^2+f'^2)e^{-s}ds\\
&\lesssim\int^{+\infty}_0(e'^2+f'^2)ds.
\end{align*}

Thus we complete the proof of this lemma.
\end{proof}

Next we establish the  basic energy estimate of (\ref{rewrite-linear-equation-of-error}).
\begin{Lemma}\label{basic-Energy-estimate}
 If $(u,v)$ satisfies \eqref{rewrite-linear-equation-of-error} and the hypothesis {\bf (H)}, then there holds
\begin{align}\label{e:basic-Energy-estimate}
&\varepsilon^2\int_{0}^{+\infty}\int_0^{2\pi}\big({\tilde{u}}_s^2+{\tilde{v}}_s^2+\tilde{u}^2+\tilde{v}^2\big)d\theta ds\nonumber\\
&\leq \varepsilon^2\eta^2\|u_0\|^2_{2}+\varepsilon^2\eta^2\|u_{0s}\|^2_{2}
+C\|e'\|^2_{2}+C\|f'\|^2_{2}\nonumber\\
&\quad+C\|\tilde{u}_\theta\|^2_{2}+C\|\tilde{v}_\theta\|^2_{2}
+\left|\int_{0}^{+\infty}\int_0^{2\pi}
\big(e^{-s}F_u\tilde{u}+ e^{-s}F_v\tilde{v}\big)d\theta ds\right|.
\end{align}
\end{Lemma}
\begin{proof}

Multiplying the first equation in \eqref{rewrite-linear-equation-of-error}  by $\tilde{u}$,
the second equation in \eqref{rewrite-linear-equation-of-error} by $\tilde{v}$, adding them together and integrating in $\mathcal C$, we obtain that
\begin{align*}
&\underbrace{\iint-\varepsilon^2\big(u_{ss}
+u_{\theta\theta}+2v_{\theta}-u\big)\tilde{u}-\varepsilon^2\big(v_{ss}+v_{\theta\theta}-2u_{\theta}-v\big)\tilde{v}}_{diffusion \ term}\\
&+\underbrace{\iint e^{-s}S_u\tilde{u}+e^{-s}S_v\tilde{v}}_{convective \ term}=\iint e^{-s}F_u\tilde{u}+ e^{-s}F_v\tilde{v},
\end{align*}
where we have used  the divergence-free condition
$\tilde{v}_s-\tilde{v}-\tilde{u}_\theta=0$.
We next deal with them term by term. \\
{\bf Diffusion term:} Integrating by parts implies
\begin{align}\label{e:diffusion-estimate-in-basic-energy-estimate}
&\iint-\varepsilon^2\big(u_{ss}
+u_{\theta\theta}+2v_{\theta}-u\big)\tilde{u}-\varepsilon^2\big(v_{ss}+v_{\theta\theta}-2u_{\theta}-v\big)\tilde{v}\nonumber\\
&=\varepsilon^2\iint \tilde{u}^2_s+\tilde{u}^2_\theta+\tilde{u}^2+
\tilde{v}^2_s+\tilde{v}^2_\theta+\tilde{v}^2-2\tilde{v}_\theta\tilde{u}+2\tilde{u}_\theta\tilde{v}.
\end{align}
{\bf Convective term:}
It is noted that
\begin{align*}
&\iint\big(e^{-s}S_u\tilde{u}+e^{-s}S_v\tilde{v}\big)\\
&=\underbrace{\iint e^{-s}(u^au_\theta+v^aru_r+uu^a_\theta+vru^a_r+v^au+vu^a)\tilde{u}}_{I_1}
\\&\quad+\underbrace{\iint e^{-s}(u^av_\theta+v^arv_r+uv^a_\theta+vrv^a_r-2uu^a)\tilde{v}}_{I_2}.
\end{align*}
{\bf 1)Estimate of $I_1$:} We estimate $I_1$ term by term and prove
\begin{align}
	\label{I1}
	I_1&\leq \varepsilon^2\eta^2\|u_0\|^2_{2}+\varepsilon^2\eta^2\|u_{0s}\|^2_{2}+\varepsilon^2\eta^2\|\tilde{u}_{s}\|^2_{2}\nonumber\\&\quad+C\|e'\|^2_{2}+C\|f'\|^2_{2}+C\|\tilde{u}_\theta\|^2_{2}+C\|\tilde{v}_\theta\|^2_{2}.
\end{align}

  We first observe that
\begin{align*}
\iint e^{-s}(ar+\frac{b}{r})u_\theta\tilde{u} =0
\end{align*}
and in terms of (\ref{estimate-on-Euler-parts}), the Poincar\'{e} inequality, \eqref{decompose-u} and Lemma \ref{hardy-for-c-d-e-f}, we obtain that
\begin{align*}
&\Big|\iint e^{-s}\big(u^a-u_e(r)\big)u_\theta\tilde{u} \Big|\\
&\leq\Big|\iint e^{-s}\big(u^a-u_e(r)\big)\partial_\theta u_1\tilde{u} \Big|+\Big|\iint e^{-s}\big(u^a-u_e(r)\big)\partial_\theta \tilde{u}\tilde{u} \Big|\\
&\leq\Big(\int_{0}^{+\infty} e^{-s}(c^2+d^2)ds\Big)^{\frac{1}{2}}\Big(\iint\tilde{u}^2_\theta \Big)^{\frac{1}{2}}+\iint\tilde{u}^2_\theta \\
&\leq C\big(\|e'\|^2_{2}+\|f'\|^2_{2}+\|\tilde{u}_\theta\|^2_{2}\big).
\end{align*}
Therefore we have
\begin{align}
\Big|\iint e^{-s}u^au_\theta\tilde{u} \Big|\lesssim\|e'\|^2_{2}+\|f'\|^2_{2}+\|\tilde{u}_\theta\|^2_{2}.\label{I1-1}
\end{align}

Noting  $u_s=-u_rr$, thanks to (\ref{estimate-on-Euler-parts}), (\ref{estimate-on-Prandtl-parts}) and \eqref{decompose-u}, one has
\begin{align*}
\Big|\iint e^{-s}v^aru_r \tilde{u}\Big|&=\Big|\iint e^{-s}v^a
(u_0+u_1+\tilde{u})_s \tilde{u}\Big|\\
&\lesssim\varepsilon\eta\Big|\iint
e^{-s}(u_0+\tilde{u})_s\tilde{u}\Big|+\Big|\underbrace{\iint e^{-s}v^a
(u_1)_s \tilde{u}}_{I_{11}}\Big|\\
&\lesssim\varepsilon^2\eta^2\|{u_0}_s\|^2_{2}+\varepsilon^2\eta^2\|\tilde{{u}}_s\|^2_{2}+\|\tilde{u}_\theta\|^2_{2}+ I_{11}.
\end{align*}
By integration by parts we have
\begin{align*}
I_{11}&=-\iint e^{-s}v^a_su_1\tilde{u}+\iint e^{-s}v^au_1(\tilde{u}-\tilde{u}_s)\\&\lesssim\iint e^{-2s}u_1^2+\iint\tilde{u}^2+\varepsilon^2\eta^2\iint \tilde{u}_s^2\\
&\lesssim\|e'\|^2_{2}+\|f'\|^2_{2}+\iint\tilde{u}^2_\theta+\varepsilon^2\eta^2\iint \tilde{u}^2_s,
\end{align*}
where we have used the fact that $\partial_sv^a_p=O(1)$.
Accordingly one has
\begin{align}
	\Big|\iint e^{-s}v^aru_r \tilde{u}\Big|\lesssim\|e'\|^2_{2}+\|f'\|^2_{2}+\iint\tilde{u}^2_\theta+\varepsilon^2\eta^2\iint \tilde{u}^2_s.\label{I1-2}
\end{align}

Moreover, we decompose $u$ as in \eqref{decompose-u} and get
\begin{align}
\Big|\iint e^{-s}uu^a_\theta\tilde{u}\Big|&=
\Big|\iint e^{-s}(u_0+u_1+\tilde{u})u^a_\theta\tilde{u}\Big|\nonumber\\
&\leq\Big|\underbrace{\iint e^{-s}u_0u^a_\theta\tilde{u}}_{I_{12}}\Big|+ C\big(\|e'\|^2_{2}+\|f'\|^2_{2}+\|\tilde{u}_\theta\|^2_{2}\big).\label{I1-3-1}
\end{align}
For $I_{12}$, one decomposes $u^a$ into the Euler part and Prandtl part, and then has
\begin{align}
	I_{12}&\lesssim\varepsilon\eta\iint |u_0\tilde{u}|+\Big|\iint e^{-s}u_0\chi(r)\partial_{\theta} u_p^{(0)}\tilde{u}\Big|\nonumber
	\\&
	\lesssim\varepsilon\eta\iint |u_0\tilde{u}|+\varepsilon\Big|\iint e^{-s}\frac{u_0\chi(r)}{r-1}Y\partial_{\theta} u_p^{(0)}\tilde{u}\Big|\nonumber
	\\&\lesssim
	\varepsilon^2\eta^2\|u_0\|^2_{2}+\|\tilde{u}_\theta\|^2_{2}+\varepsilon^2\eta^2
	\int_0^{2\pi}	\int_{\frac{1}{2}}^1\frac{(u_0\chi(r))^2}{(1-r)^2}drd\theta\nonumber
	\\&\lesssim
	\varepsilon^2\eta^2\|u_0\|^2_{2}+\|\tilde{u}_\theta\|^2_{2}+\varepsilon^2\eta^2
	\int_0^{2\pi}	\int_{\frac{1}{2}}^1\big(\partial_r(u_0\chi(r))\big)^2drd\theta\nonumber
	\\&\lesssim
	\varepsilon^2\eta^2\|u_0\|^2_{2}+\varepsilon^2\eta^2\|u_{0s}\|^2_{2}+\|\tilde{u}_\theta\|^2_{2},\label{I1-3-2}
\end{align}
where we have used  \eqref{weight-Hardy}.
Substituting \eqref{I1-3-2}  into \eqref{I1-3-1} leads to
\begin{align}
	\Big|\iint e^{-s}uu^a_\theta\tilde{u}\Big|\leqslant 	\varepsilon^2\eta^2\|u_0\|^2_{2}+\varepsilon^2\eta^2\|u_{0s}\|^2_{2}+ C\big(\|e'\|^2_{2}+\|f'\|^2_{2}+\|\tilde{u}_\theta\|^2_{2}\big).\label{I1-3}
	\end{align}

It follows from the definition of $u^a$ that
\begin{align*}
\iint e^{-s}vru^a_r\tilde{u}&=\iint e^{-s}vr(u_e(r))_r\tilde{u}+\iint e^{-s}vr(u^a_e-u_e(r))_r\tilde{u}+\iint e^{-s}vr(\chi(r)u^a_p)_r\tilde{u}\\
&=\iint e^{-s}\tilde{v}r(u_e(r))_r\tilde{u}+\iint e^{-s}(v_1+\tilde{v})r(u^a_e-u_e(r))_r\tilde{u}\\&\qquad+\iint e^{-s}(v_1+\tilde{v})r\chi'(r)u^a_p\tilde{u}+\iint e^{-s}vr\chi(r)(u^a_p)_r\tilde{u}\\&\leqslant C\|e'\|^2_{2}+C\|f'\|^2_{2}+C\|\tilde{u}_\theta\|^2_{2}+C\|\tilde{v}_\theta\|^2_{2}
+\underbrace{\iint e^{-s}vr\chi(r)(u^a_p)_r\tilde{u}}_{I_{13}}.
\end{align*}
By the divergence-free condition, \eqref{weight-Hardy}, the Poincar\'{e} inequality, \eqref{decompose-u} and \eqref{decompose-u2}, we obtain
\begin{align}
I_{13}&=-\int_{0}^{1}\int_0^{2\pi} vr\chi(r)\partial_Yu^a_p\frac{Y}{r-1}\tilde{u}(s(r))d\theta dr\nonumber\\
&\lesssim\Big(\int_{0}^{1}\int_0^{2\pi}\frac{(rv\chi)^2}{(r-1)^2}d\theta dr\Big)^{\frac{1}{2}}\Big(\int_{0}^{1}\int_0^{2\pi}\tilde{u}^2(r)d\theta dr\Big)^{\frac{1}{2}}\nonumber\\
&\lesssim\Big(\int_{0}^{1}\int_0^{2\pi}{(u_\theta(r)\chi+\chi'rv(r))^2}d\theta dr\Big)^{\frac{1}{2}}\Big(\int_{0}^{1}\int_0^{2\pi}\tilde{u}_\theta^2(r)d\theta dr\Big)^{\frac{1}{2}}\nonumber\\
&\lesssim C\|e'\|^2_{2}+C\|f'\|^2_{2}+C\|\tilde{u}_\theta\|^2_{2}+C\|\tilde{v}_\theta\|^2_{2}.\label{I13}
\end{align}
Thus we have
\begin{align}
	\iint e^{-s}vru^a_r\tilde{u}\leqslant C\|e'\|^2_{2}+C\|f'\|^2_{2}+C\|\tilde{u}_\theta\|^2_{2}+C\|\tilde{v}_\theta\|^2_{2}.\label{I1-4}
\end{align}

Furthermore,
\begin{align}
\iint e^{-s}v^au\tilde{u}\lesssim\varepsilon\eta\bigg|\iint e^{-s}(u_0+u_1+\tilde{u})\tilde{u}\bigg|
\lesssim \varepsilon^2\eta^2\|u_0\|^2_{2}+\|e'\|^2_{2}+\|f'\|^2_{2}+\|\tilde{u}_\theta\|^2_{2}\label{I1-5}
\end{align}
and
\begin{align}
&\iint e^{-s}vu^a\tilde{u}=\iint e^{-s}(v_1+\tilde{v})u^a\tilde{u}\lesssim \|e'\|^2_{2}+\|f'\|^2_{2}+\|\tilde{u}_\theta\|^2_{2}+\|\tilde{v}_\theta\|^2_{2}.\label{I1-6}
\end{align}

We then obtain the desired conclusion \eqref{I1} by collecting \eqref{I1-1}, \eqref{I1-2}, \eqref{I1-3}, \eqref{I1-4}, \eqref{I1-5} and \eqref{I1-6}.

{\bf 2)Estimate of $I_2$:} We estimate $I_2$ term by term and prove
\begin{align}
	\label{I2}
	I_2&\leq \varepsilon^2\eta^2\|u_0\|^2_{2}+\varepsilon^2\eta^2\|u_{0s}\|^2_{2}+\varepsilon^2\eta^2\|\tilde{v}_{s}\|^2_{2}\nonumber\\&\quad+C\|e'\|^2_{2}+C\|f'\|^2_{2}+C\|\tilde{u}_\theta\|^2_{2}+C\|\tilde{v}_\theta\|^2_{2}.
\end{align}

First, it follows from  \eqref{decompose-u} and \eqref{decompose-u2} that
\begin{align}
\iint e^{-s}u^av_\theta\tilde{v}=\iint e^{-s}u^a(v_1+\tilde{v})_\theta\tilde{v}\leq C\|e'\|^2_{2}+C\|f'\|^2_{2}+C\|\tilde{v}_\theta\|^2_{2}\label{I2-1}
\end{align}
and
\begin{align}
\iint e^{-s}v^arv_r\tilde{v}&=-\iint e^{-s}v^av_s\tilde{v}\nonumber\\&=-\iint e^{-s}v^a(v_1+\tilde{v})_s\tilde{v}\nonumber\\
&\lesssim \|e'\|^2_{2}+\|f'\|^2_{2}+\|\tilde{v}_\theta\|^2_{2}+\varepsilon^2\eta^2\|\tilde{v}_s\|^2_{2}\label{I2-2}
\end{align}
and
\begin{align}
\iint e^{-s}uv^a_\theta\tilde{u}&=\iint e^{-s}v^a_\theta(u_0+u_1+\tilde{u})\tilde{u}\nonumber\\
&\lesssim \varepsilon^2\eta^2\|u_0\|^2_{2}+\|e'\|^2_{2}+\|f'\|^2_{2}+\|\tilde{u}_\theta\|^2_{2}.\label{I2-3}
\end{align}

Secondly, thanks to $v^a_r=O(1)$, one gets
\begin{align}
\iint e^{-s}vrv^a_r\tilde{v}&=\iint e^{-s}(v_1+\tilde{v})rv^a_r\tilde{v}\lesssim \|e'\|^2_{2}+\|f'\|^2_{2}+\|\tilde{v}_\theta\|^2_{2}.\label{I2-4}
\end{align}

Thirdly,
noting that
\begin{align*}
\iint-e^{-s}2u(ar+\frac{b}{r})\tilde{v}&=-2\iint (b+ar^2)\tilde{u}\tilde{v}\lesssim \|\tilde{v}_\theta\|^2_{2}+\|\tilde{u}_\theta\|^2_{2}\\
	-\iint e^{-s}2u(u^a_e-u_e(r))\tilde{v}&=\iint e^{-s}(u^a_e-u_e(r))(u_0+u_1+\tilde{u})\tilde{v}\\
	&\lesssim\varepsilon^2\eta^2\|u_0\|^2_{2}+\|e'\|^2_{2}+\|f'\|^2_{2}+\|\tilde{u}_\theta\|^2_{2}+\|\tilde{v}_\theta\|^2_{2}
\end{align*}
and
\begin{align*}
-\iint e^{-s}2u(\chi(r)u^a_p)\tilde{v}&=-\iint e^{-s}2(u_0+u_1+u_2)(\chi(r)u^a_p)\tilde{v}\\
&\lesssim\varepsilon^2\eta^2\|u_0\|^2_{2}+\varepsilon^2\eta^2\|u_{0s}\|^2_{2}+ \|e'\|^2_{2}+\|f'\|^2_{2}+\|\tilde{u}_\theta\|^2_{2}+\|\tilde{v}_\theta\|^2_{2},
\end{align*}
where we have used
\begin{align*}-\iint e^{-s}2u_0(\chi(r)u^a_p)\tilde{v}=
&-2\varepsilon\int_0^1\int_0^{2\pi}\frac{\chi u_0}{r-1}Yu^a_p\tilde{v}d\theta dr\\
&\lesssim\varepsilon\eta\Big(\int_0^1\int_0^{2\pi} \frac{r(\chi u_0)^2}{(r-1)^2}d\theta dr\Big)^{\frac{1}{2}}\Big(\int_0^1\int_0^{2\pi}\tilde{v}^2_\theta d\theta dr \Big)^{\frac{1}{2}}\\
&\lesssim\varepsilon^2\eta^2\|u_0\|^2_{2}+\varepsilon^2\eta^2\|u_{0s}\|^2_{2}+\|\tilde{v}_\theta\|^2_{2}.
\end{align*}
Thus there holds
\begin{align}
	-\iint e^{-s}2uu^a\tilde{v}\lesssim\varepsilon^2\eta^2\|u_0\|^2_{2}+\varepsilon^2\eta^2\|u_{0s}\|^2_{2}+ \|e'\|^2_{2}+\|f'\|^2_{2}+\|\tilde{u}_\theta\|^2_{2}+\|\tilde{v}_\theta\|^2_{2}.\label{I2-5}
\end{align}

With the help of \eqref{I2-1},\eqref{I2-2}, \eqref{I2-3}, \eqref{I2-4} and \eqref{I2-5} we immediately arrive at \eqref{I2}.

Collecting the estimates (\ref{e:diffusion-estimate-in-basic-energy-estimate}), \eqref{I1} and \eqref{I2},  and choosing small $\eta$, we  obtain (\ref{e:basic-Energy-estimate}) and hence complete the proof of this lemma.
\end{proof}

Next, we give the energy estimate of the zero-mode $u_0$ of the solution $u$.

 \begin{Lemma}\label{basic-Energy-estimate-for-zero-frequency}
If $(u,v)$ satisfies  \eqref{rewrite-linear-equation-of-error} and the hypothesis {\bf (H)}, then there holds
\begin{align}
\varepsilon^2\int_{0}^{+\infty}\int_0^{2\pi} u^2_0+u^2_{0s}d\theta ds\leq &C\int_{0}^{+\infty}\int_0^{2\pi}(\tilde{u}^2_\theta+\tilde{v}^2_\theta)d\theta ds+C\|e'\|^2_{2}+C\|f'\|^2_{2}\nonumber\\
\ \ \ &+C\left|\int_{0}^{+\infty}\int_0^{2\pi}e^{-s}F_uu_0d\theta ds\right|.\label{e:basic-Energy-estimate-2}
\end{align}
\end{Lemma}
\begin{proof}
Multiplying the first equation of  \eqref{rewrite-linear-equation-of-error} by $u_0$  and integrating by parts in $\mathcal C$, we arrive at
\begin{align}\label{energy-estimate-for-zero-frequency}
\varepsilon^2\iint
u^2_0+u^2_{0s}+\iint e^{-s}S_uu_0 =\iint e^{-s}F_uu_0.
\end{align}

Now we deal with the convective term:
\begin{align}
\iint e^{-s}S_uu_0=\iint e^{-s}(u^au_\theta+v^aru_r+uu^a_\theta+vru^a_r+v^au+vu^a)u_0.\label{con}
\end{align}

Noting that
\begin{align}
&\iint e^{-s}(u^au_\theta+uu^a_\theta)u_0=\iint e^{-s}\partial_\theta(u^au)u_0=-\iint e^{-s}(u^au)\partial_\theta u_0=0.\label{con-1}
\end{align}

Using the integration by parts, we have
\begin{align}
\iint e^{-s}v^aru_ru_0
&=-\iint e^{-s}v^au_su_0\nonumber\\
&=-\iint e^{-s}v^a(u_1+\tilde{u})_su_0\nonumber\\
&=\underbrace{-\iint e^{-s}v^a(u_1+\tilde{u})(u_0-u_{0s})}_{I_{31}}+\underbrace{\iint e^{-s}v^a_s(u_1+\tilde{u})u_0}_{I_{32}},\label{I312}
\end{align}
where we have used $\iint e^{-s}v^au_{0s}u_0=0$ which comes from \eqref{va-null-0}.

By Lemma \ref{hardy-for-c-d-e-f}, we have
\begin{align}
I_{31}\lesssim \varepsilon^2\eta^2\|u_0\|^2_{2}+\varepsilon^2\eta^2\|u_{0s}\|^2_{2}+\|e'\|^2_{2}+\|f'\|^2_{2}+\|\tilde{u}_\theta\|^2_{2}.\label{I31}
\end{align}
And direct computations imply that
\begin{align}
I_{32}
&=\iint e^{-s}[(v^a_e)_r+\chi'v^a_p]r(u_1+\tilde{u})u_0
+\iint e^{-2s} Y\partial_Yv^a_p\frac{\chi u_0}{r-1}(u_1+\tilde{u})\nonumber\\
&\lesssim \varepsilon^2\eta^2\|u_{0}\|^2_{2}+\varepsilon^2\eta^2\|u_{0s}\|^2_{2}+\|e'\|^2_{2}+\|f'\|^2_{2}+\|\tilde{u}_\theta\|^2_{2},\label{I32}
\end{align}
where we have used
\begin{align*}
&\iint e^{-2s} Y\partial_Yv^a_p\frac{\chi u_0}{r-1}(u_1+\tilde{u})\\
&\lesssim\varepsilon^2\eta^2\iint e^{-2s} \bigg(\frac{\chi u_0}{r-1}\bigg)^2+\iint e^{-s} (u_1+\tilde{u})^2
\\
&\lesssim\varepsilon^2\eta^2\iint e^{-2s} \big(\big(\chi'u_0\big)^2+\big(\chi u_{0r}\big)^2\big)+\iint e^{-s} (u_1+\tilde{u})^2\\
&\lesssim \varepsilon^2\eta\|u_{0s}\|^2_{2}+\varepsilon^2\eta^2\|u_{0}\|^2_{2}+\|e'\|^2_{2}+\|f'\|^2_{2}+\|\tilde{u}_\theta\|^2_{2}.
\end{align*}
Substituting \eqref{I31}  and  \eqref{I32}  into \eqref{I312}  leads  to
\begin{align}
\iint e^{-s}v^aru_ru_0\lesssim \varepsilon^2\eta^2\|u_0\|^2_{2}+\varepsilon^2\eta^2\|u_{0s}\|^2_{2}+\|e'\|^2_{2}+\|f'\|^2_{2}+\|\tilde{u}_\theta\|^2_{2}.\label{con-2}
\end{align}

Due to  \eqref{v-ava}, we have
\begin{align*}
	\iint e^{-s}vru_e(r)u_0=0
\end{align*}
and then
\begin{align*}
\iint e^{-s}vru^a_ru_0&=\iint e^{-s}vr(u^a_e-u_e(r))_ru_0+\iint e^{-s}vr\chi'u^a_pu_0+\iint e^{-s}vr\chi\partial_Yu^a_p\frac{Y}{r-1}u_0\\
&=\underbrace{\iint e^{-s}(v_1+\tilde{v})r(u^a_e-u_e(r))_ru_0}_{I_{34}}+\underbrace{\varepsilon\iint e^{-2s}(v_1+\tilde{v})Yu^a_p\frac{\chi'u_0}{r-1}}_{I_{35}}\\
&\qquad+\underbrace{\varepsilon\iint e^{-s}\frac{\sqrt{\chi}vr}{r-1}Y^2\partial_Yu^a_p\frac{\sqrt{\chi} u_0}{r-1}}_{I_{36}}.
\end{align*}

Thanks to Lemma \ref{hardy-for-c-d-e-f} one has
\begin{align*}
&I_{34}+I_{35}\lesssim\varepsilon^2\eta^2\|u_{0}\|^2_{2}+\|e'\|^2_{2}+\|f'\|^2_{2}+\|\tilde{v}_\theta\|^2_{2}.
\end{align*}
It follows from  Lemma \ref{H1} and Lemma \ref{hardy-for-c-d-e-f} that
\begin{align*}
I_{36}&\lesssim\varepsilon\eta\Big(\int_{0}^{1}\int_0^{2\pi} \frac{r(\sqrt{\chi} u_0)^2}{(r-1)^2}d\theta dr\Big)^{\frac{1}{2}}\Big(\int_{0}^{1}\int_0^{2\pi}\frac{r(\sqrt{\chi}rv)^2}{(r-1)^2}d\theta dr\Big)^{\frac{1}{2}}
\\&\lesssim\varepsilon^2\eta^2\|u_{0}\|^2_{2}+\varepsilon^2\eta^2\|u_{0s}\|^2_{2}+\|e'\|^2_{2}+\|f'\|^2_{2}+\|\tilde{u}_\theta\|^2_{2},
\end{align*}
where we have used the divergence-free condition: $\partial_r(rv)=-u_\theta$.

Consequently one has
\begin{align}
\iint e^{-s}vru^a_ru_0\lesssim\varepsilon^2\eta^2\|u_{0}\|^2_{2}+\varepsilon^2\eta^2\|u_{0s}\|^2_{2}+\|e'\|^2_{2}+\|f'\|^2_{2}+\|\tilde{u}_\theta\|^2_{2}.\label{con-3}
\end{align}

Moreover, we can find that
\begin{align}
\iint e^{-s}v^auu_0&=\iint e^{-s}v^a(u_1+\tilde{u})u_0\nonumber
\\&\lesssim\varepsilon^2\eta^2\|u_{0}\|^2_{2}+\|e'\|^2_{2}+\|f'\|^2_{2}+\|\tilde{u}_\theta\|^2_{2},\label{con-4}
\end{align}
and
\begin{align}
\iint e^{-s}vu^au_0&=\iint e^{-s}v(u^a_e-u_e(r))u_0+\iint e^{-s}v\chi(r)u^a_pu_0\nonumber\\&=\iint e^{-s}(v_1+\tilde{v})(u^a_e-u_e(r))u_0+\varepsilon\int_{0}^{1}\int_0^{2\pi} e^{-s}(v_1+\tilde{v})Yu^a_p\frac{\chi(r)u_0}{r-1}d\theta dr\nonumber
\\
&\lesssim\varepsilon^2\eta^2\|u_{0}\|^2_{2}+\varepsilon^2\eta^2\|u_{0s}\|^2_{2}+C\|e'\|^2_{2}+C\|f'\|^2_{2}+C\|\tilde{v}_\theta\|^2_{2}.\label{con-5}
\end{align}

Substituting \eqref{con-1},  \eqref{con-2},  \eqref{con-3} , \eqref{con-4}  and \eqref{con-5} into  \eqref{con} we arrive at
\begin{align}
	\bigg|\iint e^{-s}S_uu_0\bigg|\lesssim\varepsilon^2\eta^2\|u_{0}\|^2_{2}+\varepsilon^2\eta^2\|u_{0s}\|^2_{2}+\|e'\|^2_{2}+\|f'\|^2_{2}+\|\tilde{u}_\theta\|^2_{2}+\|\tilde{v}_\theta\|^2_{2}\label{con-0}
\end{align}
which, together with \eqref{energy-estimate-for-zero-frequency}, implies \eqref{e:basic-Energy-estimate-2} for small $\eta$.

 Thus we complete the proof of this lemma.
\end{proof}

\subsubsection{Positivity estimate}
\indent

In this subsection, we aim to establish the positivity estimates Lemma \ref{positivity-estimate} for the higher-mode component  $(\tilde{u},\tilde{v})$ and Lemma \ref{positivity-estimate-for-first-Fourier-mode} for the first Fourier mode $(u_1,v_1)$.
To proceed we give a modified version of the classical Wirtinger type inequality.
\begin{Lemma}\label{Wirtinger-inequality}
Let $H^1_{per}(0,T)$ be the closure of a smooth function with period $T$ in $H^1(0,T)$. For $u\in H^1_{per}(0,T),$
if $\int^T_0u(t)dt=0$ and $\int^T_0u(t)\cos\frac{2\pi t}{T}dt=\int^T_0u(t)\sin\frac{2\pi t}{T}dt=0$, one has
\begin{align}\label{Wirtinger-inequality-1}
\int^T_0|u|^2dt\leq\frac{T^2}{16\pi^2}\int^T_0|u'(t)|^2dt.
\end{align}
\end{Lemma}
\begin{proof}
Due to $\int^T_0u(t)dt=0$ and $\int^T_0u(t)\cos\frac{2\pi}{T}dt=\int^T_0u(t)\sin\frac{2\pi}{T}dt=0$, there holds
\begin{equation*}
  u(t)=\sum^{\infty}_{k=2}(a_k\cos\frac{2k\pi}{T}t+b_k\sin\frac{2k\pi}{T}t)
\end{equation*}
and then
\begin{equation*}
  u'(t)=\frac{2\pi}{T}\sum^{\infty}_{k=2}(-ka_k\sin\frac{2k\pi}{T}t+kb_k\cos\frac{2k\pi}{T}t).
\end{equation*}

It follows from the Parseval's inequality that
\begin{equation*}
  \int^T_0|u|^2 dt =\frac{T}{2}\sum^{\infty}_{k=2}(|a_k|^2+|b_k|^2)\leq\frac{1}{4}\cdot\frac{T}{2}\sum^{\infty}_{k=2}k^2(|a_k|^2+|b_k|^2)=\frac{T^2}{16\pi^2}\int^T_0|u'(t)|^2 dt.
\end{equation*}

Thus, we complete the proof of the lemma.
\end{proof}

\begin{Lemma}\label{positivity-estimate}
If $(u,v)$ satisfies  \eqref{rewrite-linear-equation-of-error} and the hypothesis {\bf (H)}, then there exists $\eta_0>0$ such that for any $\eta\in (0,\eta_0)$, there holds
\begin{align}\label{e:positivity-estimate-1}
&\int_{0}^{+\infty}\int_0^{2\pi}(\tilde{u}_\theta^2+\tilde{v}_\theta^2)d\theta ds
\nonumber\\ &\leq \varepsilon^2\eta\Big(\int_{0}^{+\infty}\int_0^{2\pi} \big(\tilde{u}_s^2+\tilde{v}_s^2\big)d\theta ds\Big)+\varepsilon^2\eta\Big(\int_{0}^{+\infty}\int_0^{2\pi}\big(u_{0s}^2+u_{0}^2\big)d\theta ds\Big)\nonumber\\
&\quad+\eta\|e'\|^2_{2}+\eta\|f'\|^2_{2}
+C\left|\int_{0}^{+\infty}\int_0^{2\pi}
\big(e^{-s}F_u\tilde{u}_\theta+ e^{-s}F_v\tilde{v}_\theta\big)
d\theta ds\right|.
\end{align}
\end{Lemma}
\begin{proof}
Multiplying the first equation of \eqref{rewrite-linear-equation-of-error} by $\tilde{u}_\theta$ and the second equation by $\tilde{v}_\theta$, integrating in $\bar{\Omega}$ and adding the resultant together, we arrive at
\begin{align}\label{p:positivity-estimate}
&\iint-\varepsilon^2\big(u_{ss}
+u_{\theta\theta}+2v_{\theta}-u\big)\tilde{u}_\theta-\varepsilon^2\big(v_{ss}+v_{\theta\theta}-2u_{\theta}-v\big)\tilde{v}_\theta
\nonumber\\
&+\iint\big(e^{-s}S_u\tilde{u}_\theta+e^{-s}S_v\tilde{v}_\theta\big)=\iint\big(e^{-s}F_u\tilde{u}_\theta+e^{-s}F_v\tilde{v}_\theta\big),
\end{align}
where we have used  the divergence-free condition
$\tilde{v}_s-\tilde{v}-\tilde{u}_\theta=0$.

We first deal with the terms related to $S_u, S_v$, which provide positive terms. We write
\begin{align}
&\iint\big(e^{-s}S_u\tilde{u}_\theta+e^{-s}S_v\tilde{v}_\theta\big)\nonumber\\
&=\underbrace{\iint e^{-s}\big((\frac{b}{r}+ar)u_\theta+v(-\frac{b}{r}+ar)+v(\frac{b}{r}+ar)\big)\tilde{u}_\theta
+\iint e^{-s}\big((ar+\frac{b}{r})v_\theta-(2ar+\frac{2b}{r})u\big)\tilde{v}_\theta }_{I_4}\nonumber\\
&\quad+\underbrace{\iint\big(({u}^a-u_e(r))u_\theta+{v}^aru_r+uu^a_\theta+vr(u^a-u_e(r))_r+{v}^au+v(u^a-u_e(r))\big) e^{-s}\tilde{u}_\theta}_{I_5}\nonumber\\
&\quad+\underbrace{\iint\big((u^a-u_e(r))v_\theta+{v}^arv_r+u{v}^a_\theta+vr{v}^a_r-2u(u^a-u_e(r))\big)e^{-s}\tilde{v}_\theta}_{I_6}.\label{suv}
\end{align}

{\bf 1)Estimate of $I_4$.} This term will be a positive term. First, note that $\inf\limits_{r\in(0,1]} (ar^2+b)>0$  and
\begin{align*}
	2&\iint ar^2\tilde{v}\tilde{u}_\theta =-2\int_{0}^{1}\int_0^{2\pi} ar\tilde{v}\tilde{u}_\theta d\theta dr\\
	&=2\int_{0}^{1}\int_0^{2\pi} ar\tilde{v}(r\tilde{v})_r d\theta dr=a\big((r\tilde{v})^2|_{r=1}-(r\tilde{v})^2|_{r=0}\big)=0.
\end{align*}
Then using the Cauchy--Schwarz inequality, Wirtinger's inequality \eqref{Wirtinger-inequality-1}($T=2\pi$)
and direct computation gives that
\begin{align}
I_4&=\iint b\tilde{u}^2_\theta+\iint b\tilde{v}^2_\theta-2\iint bu\tilde{v}_\theta+\iint ar^2\tilde{u}^2_\theta+\iint ar^2\tilde{v}^2_\theta
+4\iint ar^2\tilde{v}\tilde{u}_\theta\nonumber\\
&=\iint(b+ar^2)\tilde{u}^2_\theta+\iint(b+ar^2)\tilde{v}^2_\theta+2\iint(b+ar^2)\tilde{v}\tilde{u}_\theta\nonumber\\
&\geq\iint (b+ar^2)\tilde{u}^2_\theta+(b+ar^2)\tilde{v}^2_\theta-\frac{1}{2}(b+ar^2)\tilde{u}^2_\theta-2(b+ar^2)\tilde{v}^2\nonumber\\
&\geq\frac{1}{2}\iint (b+ar^2)\tilde{u}^2_\theta+\frac{1}{2}\iint (b+ar^2)\tilde{v}^2_\theta\nonumber\\
&\gtrsim\iint\tilde{u}^2_\theta+\tilde{v}^2_\theta.\label{I4}
\end{align}

{\bf 2)Estimate of $I_5$.}
Using (\ref{estimate-on-Euler-parts}), (\ref{estimate-on-Prandtl-parts}) and Lemma \ref{hardy-for-c-d-e-f}, we obtain
\begin{align}
\iint e^{-s}({u}^a-u_e(r))u_\theta\tilde{u}_\theta&=\iint e^{-s}({u}^a-u_e(r))(c\sin\theta+d\cos\theta+\tilde{u})_\theta\tilde{u}_\theta\nonumber\\
&\lesssim \eta\|e'\|^2_{2}+\eta\|f'\|^2_{2}+\eta\|\tilde{u}_\theta\|^2_{2}\label{I5-1}
\end{align}
and
\begin{align}
\iint e^{-s}v^aru_r\tilde{u}_\theta&=-\iint e^{-s}v^au_s\tilde{u}_\theta\nonumber\\&=-\iint e^{-s}v^a(u_0+\tilde{u})_s\tilde{u}_\theta-\iint e^{-s}v^a(u_1)_s\tilde{u}_\theta\nonumber\\&=-\iint e^{-s}v^a(u_0+\tilde{u})_s\tilde{u}_\theta-\iint e^{-s}v^au_1\tilde{u}_\theta+\iint e^{-s}v_s^au_1\tilde{u}_\theta+\iint e^{-s}v^au_1\tilde{u}_{\theta s}\nonumber
\\&=-\iint e^{-s}v^a(u_0+\tilde{u})_s\tilde{u}_\theta-\iint e^{-s}v^au_1\tilde{u}_\theta+\iint e^{-s}v_s^au_1\tilde{u}_\theta\nonumber\\&\quad-\iint e^{-s}v_\theta^au_1\tilde{u}_{s}-\iint e^{-s}v_\theta^a(u_1)_\theta\tilde{u}_{s}\nonumber\\
&\lesssim\varepsilon^2\eta\|u_{0s}\|^2_{2}+\varepsilon^2\eta\|\tilde{u}_s\|^2_{2}+\eta\|\tilde{u}_\theta\|^2_{2}+\eta\|e'\|^2_{2}+\eta\|f'\|^2_{2}.\label{I5-2}
\end{align}

According to the definition of $u^a$, the Hardy type inequality  \eqref{weight-Hardy} and Lemma \ref{hardy-for-c-d-e-f}, one has
\begin{align}
\iint e^{-s}uu^a_\theta\tilde{u}_\theta&=\iint e^{-s}u_0u^a_\theta\tilde{u}_\theta+\iint e^{-s}(u_1+\tilde{u})u^a_\theta\tilde{u}_\theta\nonumber\\
&=\iint e^{-s}u_0\partial_\theta{u^a_e}\tilde{u}_\theta +\varepsilon\iint\frac{Y}{(r-1)}u_0\chi\partial_\theta{u^a_p}\tilde{u}_\theta d\theta dr\nonumber\\&\quad+\iint e^{-s}(u_1+\tilde{u})u^a_\theta\tilde{u}_\theta
\nonumber\\&\lesssim\varepsilon^2\eta\|{u}_0\|^2_{2}+\eta\|\tilde{u}_\theta\|^2_{2}+\varepsilon\eta\Big(\int_0^1\int_0^{2\pi} \frac{({\chi} {u_0})^2}{(r-1)^2}d\theta dr\Big)^{\frac{1}{2}}\Big(\int_0^1\int_0^{2\pi} \tilde{u}^2_\theta d\theta dr\Big)^{\frac{1}{2}}\nonumber\\&\quad+\eta\|e'\|^2_{2}+\eta\|f'\|^2_{2}\nonumber\\
&\lesssim\varepsilon^2\eta\|{u}_0\|^2_{2}+\varepsilon^2\eta\|{u}_{0s}\|^2_{2}+\eta\|\tilde{u}_\theta\|^2_{2}+\eta\|e'\|^2_{2}+\eta\|f'\|^2_{2}.\label{I5-3}
\end{align}

Moreover, we have
\begin{align}
\iint e^{-s}vr(u^a-u_e(r))_r\tilde{u}_\theta&=\iint e^{-s}(v_1+\tilde{v})r\big((u^a_e-u_e(r))_r+\chi'u^a_p\big)\tilde{u}_\theta\nonumber\\
&\ \ +\iint e^{-s} Y\partial_Yu^a_p\frac{vr\chi}{r-1}\tilde{u}_\theta
\nonumber\\
&\lesssim\eta\|e'\|^2_{2}+\eta\|f'\|^2_{2}+\eta\|\tilde{v}_\theta\|^2_{2}+\eta\|\tilde{u}_\theta\|^2_{2}\nonumber\\
&\quad+\eta\Big(\int_0^1\int_0^{2\pi}\frac{({\chi}rv)^2}{(r-1)^2}d\theta dr\Big)^{\frac{1}{2}}
\Big(\int_0^1\int_0^{2\pi} \tilde{u}^2_\theta d\theta dr\Big)^{\frac{1}{2}}\nonumber\\
&\lesssim\eta\|e'\|^2_{2}+\eta\|f'\|^2_{2}+\eta\|\tilde{v}_\theta\|^2_{2}+\eta\|\tilde{u}_\theta\|^2_{2}\label{I5-4}
\end{align}
and
\begin{align}
\iint e^{-s}v^au\tilde{u}_\theta&=\iint e^{-s}v^a(u_0+u_1+\tilde{u})\tilde{u}_\theta\nonumber\\&\lesssim\varepsilon^2\eta\|u_{0}\|^2_{2}+\eta\|e'\|^2_{2}+\eta\|f'\|^2_{2}+\eta\|\tilde{u}_\theta\|^2_{2}\label{I5-5}\end{align}
and
\begin{align}
\iint e^{-s}v(u^a-u_e(r))\tilde{u}_\theta&=\iint e^{-s}(v_1+\tilde{v})(u^a-u_e(r))\tilde{u}_\theta
\nonumber\\&\lesssim\eta\|e'\|^2_{2}+\eta\|f'\|^2_{2}+\eta\|\tilde{u}_\theta\|^2_{2}+\eta\|\tilde{v}_\theta\|^2_{2}.\label{I5-6}
\end{align}

Combining \eqref{I5-1}-\eqref{I5-6} we immediately get
\begin{align}
I_5&\lesssim\varepsilon^2\eta\|{u}_0\|^2_{2}+\varepsilon^2\eta\|{u}_{0s}\|^2_{2}+\varepsilon^2\eta\|\tilde{u}_s\|^2_{2}+\eta\|\tilde{u}_\theta\|^2_{2}+\eta\|\tilde{v}_\theta\|^2_{2}+\eta\|e'\|^2_{L^2}+\eta\|f'\|^2_{2}.\label{I5}
\end{align}

{\bf 3)Estimate of $I_6$.}
First, from (\ref{estimate-on-Euler-parts}), \eqref{estimate-on-Prandtl-parts}, \eqref{decompose-u}, \eqref{decompose-u2} and Lemma \ref{hardy-for-c-d-e-f}, we obtain
\begin{align}
\iint e^{-s}(u^a-u_e(r))v_\theta \tilde{v}_\theta&=\iint e^{-s}(u^a-u_e(r))(v_{1\theta}+\tilde{v}_\theta) \tilde{v}_\theta\nonumber\\&\lesssim\eta\|e'\|^2_{2}+\eta\|f'\|^2_{2}+\eta\|\tilde{v}_\theta\|^2_{2},\label{I6-1}
\end{align}
and
\begin{align}
\iint e^{-s}v^arv_r\tilde{v}_\theta=-\iint e^{-s}v^a(v_1+\tilde{v})_s\tilde{v}_\theta\lesssim\varepsilon^2\eta\|\tilde{v}_s\|^2_{2}+\eta\|e'\|^2_{2}+\eta\|f'\|^2_{2}+\eta\|\tilde{v}_\theta\|^2_{2},\label{I6-2}
\end{align}
and
\begin{align}
\iint e^{-s}uv^a_\theta\tilde{v}_\theta=\iint e^{-s}(u_0+u_1+\tilde{u})v^a_\theta\tilde{v}_\theta\lesssim\varepsilon^2\eta\|u_{0}\|^2_{2}+\eta\|e'\|^2_{2}+\eta\|f'\|^2_{2}+\eta\|\tilde{u}_\theta\|^2_{2}
+\eta\|\tilde{v}_\theta\|^2_{2},\label{I6-3}
\end{align}
and
\begin{align}
\iint e^{-s}vrv^a_r\tilde{v}_\theta=\iint e^{-s}(v_1+\tilde{v})rv^a_r\tilde{v}_\theta\lesssim\eta\|e'\|^2_{2}+\eta\|f'\|^2_{2}+\eta\|\tilde{v}_\theta\|^2_{2},\label{I6-4}
\end{align}
where we have used the fact $v^a_r=O(\eta)$.

Secondly, one has
\begin{align}
&-2\iint e^{-s}u(u^a-u_e(r))\tilde{v}_\theta\nonumber\\
&=-2\iint e^{-s}u_0\chi u^a_p\tilde{v}_\theta-2\iint e^{-s}u_0(u^a_e-u_e(r))\tilde{v}_\theta-2\iint e^{-s}(u^{(1)}+\tilde{u})[u^a-u_e(r)]\tilde{v}_\theta\nonumber\\
&\lesssim\varepsilon^2\eta\|u_{0}\|^2_{2}+\varepsilon^2\eta\|u_{0s}\|^2_{2}+\eta\|e'\|^2_{2}+\eta\|f'\|^2_{2}+\eta\|\tilde{u}_\theta\|^2_{2}
+\eta\|\tilde{v}_\theta\|^2_{2},\label{I6-5}
\end{align}
where we have used
\begin{align*}
-2\iint e^{-s}u_0\chi u^a_p\tilde{v}_\theta&=-2\varepsilon\int_0^1\int_0^{2\pi} \frac{Y}{r-1}u_0\chi u^a_p\tilde{v}_\theta d\theta dr\\
&\lesssim\varepsilon\eta\Big(\int_0^1\int_0^{2\pi}\frac{({\chi} u_0)^2}{(r-1)^2}d\theta dr\Big)^{\frac{1}{2}}\Big(\int_0^1\int_0^{2\pi}\tilde{v}^2_\theta d\theta dr\Big)^{\frac{1}{2}}\\
&\lesssim\varepsilon^2\eta\|u_{0}\|^2_{2}+\varepsilon^2\eta\|u_{0s}\|^2_{2}
+\eta\|\tilde{v}_\theta\|^2_{2}.
\end{align*}

With the help of \eqref{I6-1}-\eqref{I6-5} we arrive at
\begin{align}
I_6\lesssim\varepsilon^2\eta\|u_{0}\|^2_{2}+\varepsilon^2\eta\|u_{0s}\|^2_{2}+\varepsilon^2\eta\|\tilde{v}_s\|^2_{2}+\eta\|e'\|^2_{2}+\eta\|f'\|^2_{2}+\eta\|\tilde{u}_\theta\|^2_{2}
	+\eta\|\tilde{v}_\theta\|^2_{2}.\label{I6}
\end{align}

Next we can deduce that
\begin{align}\label{diffusive-term-2}
&\iint-\varepsilon^2\big(u_{ss}
+u_{\theta\theta}+2v_{\theta}-u\big)\tilde{u}_\theta-\varepsilon^2\big(v_{ss}+v_{\theta\theta}-2u_{\theta}-v\big)\tilde{v}_\theta\nonumber\\
&=\iint\varepsilon^2\tilde{u}_s\tilde{u}_{s\theta}-2\varepsilon^2\tilde{v}_{\theta}\tilde{u}_\theta+\varepsilon^2\tilde{v}_s\tilde{v}_{s\theta}
+2\varepsilon^2\tilde{u}_{\theta}\tilde{v}_\theta=0.
\end{align}

Finally, collecting \eqref{suv} and the estimates \eqref{I4}, \eqref{I5}, \eqref{I6}, (\ref{diffusive-term-2}) together, we obtain (\ref{e:positivity-estimate-1}) which completes the proof of this lemma.
\end{proof}

In what follows, we will establish the  positivity estimate for first Fourier mode   $(u_1,v_1)$ of the solution $(u,v)$.
\begin{Lemma}\label{positivity-estimate-for-first-Fourier-mode}
If $(u,v)$ satisfies  \eqref{rewrite-linear-equation-of-error} and the hypothesis {\bf (H)}, then there exists $\eta_0>0$ such that for any $\eta\in (0,\eta_0)$, there holds
\begin{align}\label{e:positivity-estimate-for-first-Fourier-mode}
\int_{0}^{+\infty}\int_0^{2\pi}(e'^2+f'^2)d\theta ds
\lesssim& \varepsilon^2\eta\int_{0}^{+\infty}\int_0^{2\pi} \big(\tilde{u}_s^2+\tilde{v}_s^2\big)d\theta ds+\varepsilon^2\eta\int_{0}^{+\infty}\int_0^{2\pi}\big(u_{0s}^2+u_{0}^2\big)d\theta ds\nonumber\\
&+\eta\|\tilde{u}_\theta\|^2_{2}+\eta\|\tilde{v}_\theta\|^2_{2}
+\left|\int_{0}^{+\infty}\int_0^{2\pi}
\big(e^{-s}F_u{u}_{1\theta}+e^{-s}F_v{v}_{1\theta}\big)d\theta ds\right|.
\end{align}
\end{Lemma}

\begin{proof}
Multiplying the first equation of \eqref{rewrite-linear-equation-of-error} by ${u}_{1\theta}$ and the second equation by ${v}_{1\theta}$, integrating in $\mathcal C$ and adding the resultant  together, we arrive at
\begin{align}\label{p:positivity-estimate-first-Fourier-mode}
&\iint-\varepsilon^2\big(u_{ss}
+u_{\theta\theta}+2v_{\theta}-u\big){u}_{1\theta}-\varepsilon^2\big(v_{ss}+v_{\theta\theta}-2u_{\theta}-v\big){v}_{1\theta}\nonumber\\
&\qquad+\iint e^{-s}p_\theta{u}_{1\theta}-e^{-s}p_s{v}_{1\theta}+\iint e^{-s}S_u{u}_{1\theta}+e^{-s}S_v{v}_{1\theta}\nonumber\\
&=\iint\big(e^{-s}F_u{u}_{1\theta}+e^{-s}F_v{v}_{1\theta}\big).
\end{align}

{\bf Convective term.} We first deal with the terms related to $S_u$ and $S_v$.
We deduce that
\begin{align}
&\iint e^{-s}S_u{u}_{1\theta}+e^{-s}S_v{v}_{1\theta}\nonumber\\
&=\underbrace{\iint e^{-s}\big((ar+\frac{b}{r})u_\theta+v(-\frac{b}{r}+ar)+v(ar+\frac{b}{r})\big){u}_{1\theta}
+e^{-s}\big((ar+\frac{b}{r})v_\theta-(2ar+\frac{2b}{r})u\big){v}_{1\theta} }_{I_7}\nonumber\\
&\quad+\underbrace{\iint\big(({u}^a-u_e(r))u_\theta+{v}^aru_r+uu^a_\theta+vr(u^a-u_e(r))_r+{v}^au+v(u^a-u_e(r))\big) e^{-s}{u}_{1\theta}}_{I_8}\nonumber\\
&\quad+\underbrace{\iint\big((u^a-u_e(r))v_\theta+{v}^arv_r+u{v}^a_\theta+vr{v}^a_r-2u(u^a-u_e(r))\big)e^{-s}{v}_{1\theta}}_{I_9}.\label{suv-0}
\end{align}

{\bf 1)Estimate of $I_7$.}
Direct computation gives that
\begin{align}
I_7&=\int_0^\infty ds\int_0^{2\pi}\big(b{{u}^2_{1\theta}}+ b{v}^2_{1\theta}+2 b{u}_{1\theta}{v}_1+ ar^2{u}^2_{1\theta}+ ar^2{v}^2_{1\theta}+2 ar^2{u}_{1\theta}{v}_1\big)d\theta\nonumber\\
&=\pi\int_0^\infty \big((b+ar^2)(c^2+d^2)+(b+ar^2)(e^2+f^2)+2(b+ar^2)(cf-de)\big)ds\nonumber\\
&=\pi\int_0^\infty(b+ar^2)\big((c+f)^2+(d-e)^2\big)ds\nonumber\\
&=\pi\int_0^\infty (b+ar^2)(f'^2+e'^2)ds\nonumber\\
&\gtrsim \int_0^\infty (f'^2+e'^2)ds,\label{I7}
\end{align}
where we have used
\begin{align*}
2&\iint ar^2 u{v}_{1\theta}=-2\int_0^1\int_0^{2\pi} ar v_1{u}_{1\theta} d\theta dr=2\int_0^1\int_0^{2\pi} ar v_1{(r{v}_1)}_rd\theta dr=0
\end{align*}
and  $\inf\limits_{r\in(0,1]}\big(ar^2+b\big)>0$ .
\par

{\bf 2)Estimate of $I_8$.}
Due to \eqref{estimate-on-Euler-parts}, \eqref{estimate-on-Prandtl-parts} and Lemma \ref{hardy-for-c-d-e-f}, we have
\begin{align}
\iint e^{-s}({u}^a-u_e(r))u_\theta{u}_{1\theta}=\iint e^{-s}({u}^a-u_e(r))(u_1+\tilde{u})_\theta{u}_{1\theta}\lesssim \eta\|e'\|^2_{2}+\eta\|f'\|^2_{2}+\eta\|\tilde{u}_\theta\|^2_{2}.\label{I8-1}
\end{align}

Utilizing \eqref{decompose-u} one has
\begin{align*}
-\iint e^{-s}v^aru_ru_{1\theta}&=\iint e^{-s}v^au_su_{1\theta}\\&=\iint e^{-s}v^a(u_0+\tilde{u})_su_{1\theta}
+\iint e^{-s}v^au_{1s}u_{1\theta}\\&=:I_{81}+I_{82}.
\end{align*}
It is easy to get that
\begin{align*}
&I_{81}\lesssim\varepsilon^2\eta\|u_{0s}\|^2_{2}+\varepsilon^2\eta\|\tilde{u}_s\|^2_{2}+\eta\|e'\|^2_{L^2}+\eta\|f'\|^2_{2}.
\end{align*}
And direct computation implies
\begin{align*}
I_{82}&=\iint e^{-s}v^a(cc'\sin\theta\cos\theta-dd'\sin\theta\cos\theta-c'd\sin^2\theta+d'c\cos^2\theta)\\
&=\frac{1}{2}\iint e^{-s}(v^a-\partial_sv^a)(c^2-d^2)\sin\theta\cos\theta+e^{-s}v^a(c'd+d'c)\cos2\theta\\
&=\frac{1}{2}\iint e^{-s}(v^a-\partial_sv^a)(c^2-d^2)\sin\theta\cos\theta+e^{-s}(v^a-\partial_sv^a)cd\cos2\theta\\
&\lesssim\eta\big(\|e'\|^2_{2}+\|f'\|^2_{2}\big),
\end{align*}
where we have used the facts that  $\partial_sv^a=O(\eta)$, \eqref{hardy-for-c-d-e-f-0} and \begin{align*}\iint e^{-s}v^a(d'c-c'd)=\int^{+\infty}_0e^{-s}(d'c-c'd)ds\int^{2\pi}_0v^ad\theta=0.\end{align*}
We then obtain
\begin{align}
	-\iint e^{-s}v^aru_ru_{1\theta}\lesssim\varepsilon^2\eta\|u_{0s}\|^2_{2}+\varepsilon^2\eta\|\tilde{u}_s\|^2_{2}+\eta\|e'\|^2_{L^2}+\eta\|f'\|^2_{2}.\label{I8-2}
\end{align}

Next it follows from \eqref{decompose-u}, $\partial_\theta u^a_e=O(\varepsilon\eta)$, \eqref{weight-Hardy} and Lemma \ref{hardy-for-c-d-e-f} that
\begin{align}
\iint e^{-s}uu^a_\theta u_{1\theta}&=\iint e^{-s}u_0u^a_\theta u_{1\theta}+\iint e^{-s}(u_1 +\tilde{u})u^a_\theta u_{1\theta}\nonumber\\
&=\iint e^{-s}u_0\partial_\theta u^a_eu_{1\theta}+\varepsilon\int_0^1\int_0^{2\pi} e^{-s}\frac{\chi u_0}{r-1}Y\partial_\theta u^a_p u_{1\theta}d\theta dr+\iint e^{-s}(u_1 +\tilde{u})u^a_\theta u_{1\theta}\nonumber\\
&\lesssim\varepsilon^2\eta\|u_0\|^2_{2}+\varepsilon^2\eta\|{u}_{0s}\|^2_{2}+\eta\|e'\|^2_{2}+\eta\|f'\|^2_{2}.\label{I8-3}
\end{align}

Similarly to the process of \eqref{I5-4}, we have
\begin{align}
\iint e^{-s}vr(u^a-u_e(r))_ru_{1\theta}\lesssim\eta\|e'\|^2_{2}+\eta\|f'\|^2_{2}+\eta\|\tilde{v}_\theta\|^2_{2}+\eta\|\tilde{u}_\theta\|^2_{2}.\label{I8-4}
\end{align}

Moreover, it is easy to obtain
\begin{align}
\iint e^{-s}v^auu_{1\theta}\lesssim\varepsilon^2\eta\|u_{0}\|^2_{2}+\eta\|e'\|^2_{2}+\eta\|f'\|^2_{2}+\eta\|\tilde{u}_\theta\|^2_{2}\label{I8-5}
\end{align}
and
\begin{align}
\iint e^{-s}v(u^a-u_e(r))u_{1\theta}
\lesssim\eta\|e'\|^2_{2}+\eta\|f'\|^2_{2}+\eta\|\tilde{v}_\theta\|^2_{2}.\label{I8-6}
\end{align}

From \eqref{I8-1}-\eqref{I8-6} we immediately have \begin{align}
I_8\lesssim\varepsilon^2\eta\|u_0\|^2_{2}+\varepsilon^2\eta\|{u}_{0s}\|^2_{2}+\eta\|e'\|^2_{2}+\eta\|f'\|^2_{2}+\eta\|\tilde{v}_\theta\|^2_{2}+\eta\|\tilde{u}_\theta\|^2_{2}.\label{I8}
\end{align}

{\bf 3)Estimate of $I_9$.}
First, by (\ref{estimate-on-Euler-parts}), \eqref{estimate-on-Prandtl-parts}, Lemma \ref{hardy-for-c-d-e-f} and the Poincar\'{e} inequality (\ref{poin}), we obtain
\begin{align}
&\iint e^{-s}(u^a-u_e(r))v_\theta v_{1\theta}\lesssim\eta\|e'\|^2_{2}+\eta\|f'\|^2_{2}+\eta\|\tilde{v}_\theta\|^2_{2}\label{I9-1}
\end{align}
and
\begin{align}
\iint e^{-s}v^arv_rv_{1\theta}=\iint e^{-s}v^a(v_1+\tilde{v})_sv_{1\theta}
\lesssim\eta\|e'\|^2_{2}+\eta\|f'\|^2_{2}+\varepsilon^2\eta\|\tilde{v}_s\|^2_{2}\label{I9-2}
\end{align}
and
\begin{align}
&\iint e^{-s}uv^a_\theta v_{1\theta}=\iint e^{-s}v^a_\theta(u_0+u_1+\tilde{u})v_{1\theta}\lesssim\varepsilon^2\eta\|u_{0}\|^2_{2}+\eta\|e'\|^2_{2}+\eta\|f'\|^2_{2}+\eta\|\tilde{u}_\theta\|^2_{2}\label{I9-3}
\end{align}
and
\begin{align}
\iint e^{-s}vrv^a_rv_{1\theta}\lesssim\eta\|e'\|^2_{2}+\eta\|f'\|^2_{2}+\eta\|\tilde{v}_\theta\|^2_{2},\label{I9-4}
\end{align}
where we have used  $v^a_s=O(\eta)$.

Secondly, there holds
\begin{align}
&-2\iint e^{-s}u(u^a-u_e(r))v_{1\theta}\nonumber\\&=-2\iint e^{-s}u_0\chi u^a_pv_{1\theta}-2\iint e^{-s}u_0(u^a_e-u_e(r)v_{1\theta}-2\iint e^{-s}(u^{(1)}+\tilde{u})(u^a-u_e(r))v_{1\theta}
\nonumber\\&=-2\varepsilon\iint \frac{u_0\chi}{r-1} Yu^a_pv_{1\theta} d\theta dr-2\iint e^{-s}u_0(u^a_e-u_e(r)v_{1\theta}-2\iint e^{-s}(u^{(1)}+\tilde{u})(u^a-u_e(r))v_{1\theta}
\nonumber\\
&\lesssim\varepsilon^2\eta\|u_{0}\|^2_{2}+\varepsilon^2\eta\|u_{0s}\|^2_{2}
+\eta\|e'\|^2_{2}+\eta\|f'\|^2_{2}+\eta\|\tilde{u}_\theta\|^2_{2}.\label{I9-5}
\end{align}

Hence we get \begin{align}
	I_9\lesssim\varepsilon^2\eta\|u_{0}\|^2_{2}+\varepsilon^2\eta\|u_{0s}\|^2_{2}
	+\eta\|e'\|^2_{2}+\eta\|f'\|^2_{2}+\eta\|\tilde{u}_\theta\|^2_{2}+\eta\|\tilde{v}_\theta\|^2_{2}.\label{I9}
\end{align}

{\bf Pressure term.} By integration by parts, we deduce that
\begin{align}\label{e:pressure-estimate-in-positivity-estimate-5}
\iint e^{-s}p_\theta u_{1\theta}
-\iint e^{-s}p_sv_{1\theta}=\iint e^{-s}p(-u_{1\theta\theta}-v_{1\theta}+v_{1s\theta})=0.
\end{align}

{\bf Diffusion term.} From integration by parts we conclude
\begin{align}\label{diffusive-term-5}
&\iint-\varepsilon^2\big(u_{ss}
+u_{\theta\theta}+2v_{\theta}-u\big)u_{1\theta}-\varepsilon^2\big(v_{ss}+v_{\theta\theta}-2u_{\theta}-v\big)v_{1\theta}\nonumber\\
&=\iint\varepsilon^2u_{1s}u_{1s\theta}-2\varepsilon^2u_{1\theta}v_{1\theta}+\varepsilon^2v_{1s}v_{1s\theta}
+2\varepsilon^2u_{1\theta}v_{1\theta}=0.
\end{align}

Finally, collecting \eqref{suv-0} together with the estimates \eqref{I7}, \eqref{I8}, \eqref{I9}, (\ref{e:pressure-estimate-in-positivity-estimate-5}), and (\ref{diffusive-term-5}), we obtain \eqref{e:positivity-estimate-for-first-Fourier-mode}, which completes the proof of the lemma.
\end{proof}

\subsubsection{Linear stability estimate.}

\begin{Proposition}\label{proposition-1}
 If $(u,v)$ satisfies \eqref{rewrite-linear-equation-of-error} and the hypothesis {\bf (H)}, then there exist $\varepsilon_0>0, \eta_0>0$ such that for any $\varepsilon\in (0,\varepsilon_0), \eta\in(0,\eta_0)$, there holds
\begin{align}\label{e:linear-stability-estimates-for-linear-equation}
&\int_0^\infty(e'^2+f'^2)ds+\int_{0}^{+\infty}\int_0^{2\pi}(\tilde{u}_\theta^2+\tilde{v}_\theta^2)d\theta ds+\varepsilon^2\iint ({\tilde{u}}_s^2+{\tilde{v}}_s^2)+\varepsilon^2\iint (u^2_0+u^2_{0s})
\nonumber\\
&\lesssim\iint\big(e^{-s}F_u\tilde{u}+ e^{-s}F_v\tilde{v}\big)+
\iint e^{-s}F_uu_0\nonumber\\
&~~~~+\iint\big(e^{-s}F_u\tilde{u}_\theta+ e^{-s}F_v\tilde{v}_\theta\big)+\iint\big(e^{-s}F_uu_{1\theta}+ e^{-s}F_vv_{1\theta}\big).
\end{align}
\end{Proposition}
\begin{proof}

Combining
$\big((\ref{e:basic-Energy-estimate})+(\ref{e:basic-Energy-estimate-2})\big)+C\times\big((\ref{e:positivity-estimate-1})+(\ref{e:positivity-estimate-for-first-Fourier-mode})\big)$ for large $C$, we obtain the desired estimate \eqref{e:linear-stability-estimates-for-linear-equation}.
\end{proof}

\section{ Existence of $L^\infty$ solution for the error equations}
\indent

In this section, we first obtain partial $H^2$ estimates of $(u,v)$ and then apply the contraction mapping theorem to prove the existence of solutions to the error equations (\ref{e:error-equation}).

\subsection{Partial $H^2$ estimate for Stokes system}
\indent

A similar result to the following lemma has been  proven in \cite[Lemma 4.2]{FGLT2}.

\begin{lemma}\label{infinity-norm-estimate}
If $(u,v)$ satisfies  \eqref{rewrite-linear-equation-of-error} and the hypothesis {\bf (H)}, then
\begin{align}\label{s:Sobolev-embedding}
\|u\|_{\infty}&\leq\|u_0\|_2+\|u_{0s}\|_2+\|\tilde{u}_\theta\|_2+\|\tilde{u}_{s\theta}\|_2+\|c-c_\infty\|_{2}+\|d-d_\infty\|_{2}\nonumber\\
&\quad+\|c'\|_{L^2}+\|d'\|_{L^2}+|c_\infty|+|d_\infty|,
\end{align}
and
\begin{align}\label{s:Sobolev-embedding-2}
\|v\|_{\infty}&\leq\|\tilde{v}_\theta\|_2+\|\tilde{v}_{s\theta}\|_2+\|e-e_\infty\|_{2}+\|f-f_\infty\|_{2}+\|e'\|_{2}+\|f'\|_{2}\nonumber\\
&\quad+|e_\infty|+|f_\infty|.
\end{align}
\end{lemma}
\begin{proof}
For general $w(\theta,s)\in L^2\big(\mathcal C\big)$. We perform a Fourier series expansion to obtain
\begin{align}\label{Fourier-series-expansion}
w(\theta,s)=\sum_{k\in Z}w_k(s)e^{ik\theta}, \ \forall s\in [0,+\infty).
\end{align}
Thus
\begin{align}\label{Fourier-series-expansion-0}
|w(\theta,s)|\leq\sum_{k\in Z}|w_k(s)|, \ \forall (\theta,s)\in [0,2\pi]\times[0,+\infty).
\end{align}

Notice that
\beno
\|w_k\|_{\infty}\leq \sqrt{2}\|w_k\|^{\frac12}_2\|w'_k\|^{\frac12}_2
\eeno
and then
\beno
\sum_{k\neq 0}\|w_k\|_\infty\leq \sqrt{2}\sum_{k\neq 0}\|w_k\|^{\frac12}_2\|w'_k\|^{\frac12}_2\leq C\Big(\sum_{k\neq 0}|k|^2\|w_k\|^2_2\Big)^{\frac14}\Big(\sum_{k\neq 0}|k|^2\|w'_k\|^2_2\Big)^{\frac14}.
\eeno
Since
\beno
w_\theta(\theta,s)=\sum_{k\in Z}ikw_k(s)e^{ik\theta}, \quad w_{\theta s}(\theta,s)=\sum_{k\in Z}ikw'_k(s)e^{ik\theta},
\eeno
then
\beno
\|w_\theta\|_2^2=\sum_{k\in Z}|k|^2\|w_k\|_2^2, \quad \|w_{\theta s}\|_2^2=\sum_{k\in Z}|k|^2\|w'_k\|_2^2.
\eeno
Hence we obtain
\begin{align}\label{Fourier-series-expansion-1}
\sum_{k\neq 0}\|w_k\|_\infty\leq C\|w_\theta\|_2+C\|w_{\theta s}\|_2.
\end{align}
Moreover, there holds
\begin{align}\label{Fourier-series-expansion-2}
\|w_0\|_{\infty}\leq \sqrt{2}\|w_0\|^{\frac12}_2\|w'_0\|^{\frac12}_2.
\end{align}

Substituting \eqref{Fourier-series-expansion-1} and \eqref{Fourier-series-expansion-2} into \eqref{Fourier-series-expansion-0}  we obtain
\beno
\|w\|_{\infty}\leq C (\|w_\theta\|_2+\|w_{s\theta}\|_2+
\|w_0\|_2+\|w_{0s}\|_2)
\eeno
which, together with $u=u_0+(u_1-u_{1\infty})+\tilde{u}+u_{1\infty}$, implies \eqref{s:Sobolev-embedding} and  \eqref{s:Sobolev-embedding-2} similarly.
\end{proof}

Testing \eqref{rewrite-linear-equation-of-error} directly by $(u_1,v_1)$ does not yield the desired estimates, i.e., $\|(c-c_\infty,d-d_\infty,e-e_\infty,f-f_\infty,c',d')\|_2$. Accordingly we have to turn to
other multipliers to proceed.

\begin{lemma}\label{second-order-der-one-fre}
If $(u,v)$ satisfies  \eqref{rewrite-linear-equation-of-error} and the hypothesis {\bf (H)}, then
\begin{align}\label{eq:H2-estimate}
&\varepsilon^2\int_0^\infty\big( c''^2+d''^2+e''^2+f''^2+(c'^2+d'^2+e'^2+f'^2)\big)ds\nonumber\\
&\leq\frac{1}{\varepsilon^4}(\|e'\|^2_{2}+\|f'\|^2_{2})+\frac{1}{\varepsilon^4}(\|\tilde{u}_\theta\|^2_{2}+\|\tilde{v}_\theta\|^2_{2})
+\frac{1}{\varepsilon^4}(\|u_0\|^2_{2}+\|u_{0s}\|^2_{2}+\|\tilde{u}_{s}\|^2_{2}+
\|\tilde{v}_s\|^2_{2})\nonumber\\
&\ \ \ +\frac{1}{\varepsilon^4}\iint\big(e^{-2s}F^2_u+e^{-2s}F^2_v\big).
\end{align}
\end{lemma}

\begin{proof}
Multiplying the first equation of \eqref{rewrite-linear-equation-of-error} by $-{u}_{1ss}$ and the second equation by $-{v}_{1ss}$, integrating in $\bar{\Omega}$ and adding the resultant together, we obtain
\begin{align*}
&\iint-\varepsilon^2\big(u_{ss}
+u_{\theta\theta}+2v_{\theta}-u\big)(-{u}_{1ss})-\varepsilon^2\big(v_{ss}+v_{\theta\theta}-2u_{\theta}-v\big)(-{v}_{1ss})\\&\quad+\iint\big(e^{-s}p_\theta(-{u}_{1ss})
-e^{-s}p_s(-{v}_{1ss})\big)
+\iint\big(e^{-s}S_u(-{u}_{1ss})+e^{-s}S_v(-{v}_{1ss})\big)\\&=\iint\big(e^{-s}F_u(-{u}_{1ss})+e^{-s}F_v(-{v}_{1ss})\big).
\end{align*}
{\bf Diffusion term:} From integration by parts, we obtain
\begin{align}\label{cj1}
&\iint\varepsilon^2\big(u_{ss}
+u_{\theta\theta}+2v_{\theta}-u\big){u}_{1ss}+\varepsilon^2\big(v_{ss}+v_{\theta\theta}-2u_{\theta}-v\big){v}_{1ss}\nonumber\\
&=\iint\varepsilon^2\big(u_{1ss}\big)^2+\varepsilon^2\big(v_{1ss}\big)^2+2\varepsilon^2\big(u_{1s}\big)^2+2\varepsilon^2\big(v_{1s}\big)^2
-2\varepsilon^2v_{1s\theta}{u}_{1s}+2\varepsilon^2u_{1s\theta}{v}_{1s}\nonumber\\
&=\pi\varepsilon^2\int_0^\infty\big( c''^2+d''^2+e''^2+f''^2+2(c'^2+d'^2+e'^2+f'^2)+4(f'c'-e'd')\big)ds,
\end{align}
where we have used the fact $u_{1\theta\theta}-u_1=-2u_1$ and $v_{1\theta\theta}-v_1=-2v_1$.

{\bf Convective term:}
\begin{align*}
&\iint\big(e^{-s}S_u(-{u}_{1ss})+e^{-s}S_v(-{v}_{1ss})\big)\\
&=\underbrace{\iint e^{-s}\big((ar+\frac{b}{r})u_\theta+v(ar-\frac{b}{r})+v(ar+\frac{b}{r})\big)(-{u}_{1ss})
+e^{-s}\big((ar+\frac{b}{r})v_\theta-(2ar+\frac{2b}{r})u\big)(-{v}_{1ss})}_{J_1}\\
&\quad+\underbrace{\iint\big(({u}^a-u_e(r))u_\theta+{v}^aru_r+uu^a_\theta+vr(u^a-u_e(r))_r+{v}^au+v(u^a-u_e(r))\big) e^{-s}(-{u}_{1ss})}_{J_2}\\
&\quad+\underbrace{\iint\big((u^a-u_e(r))v_\theta+{v}^arv_r+u{v}^a_\theta+vr{v}^a_r-2u(u^a-u_e(r))\big)e^{-s}(-{v}_{1ss})}_{J_3}.
\end{align*}

It follows from the Young's inequality and Lemma \ref{hardy-for-c-d-e-f} that
\begin{align*}
J_1&=-\iint2bu_{1s}v_{1s}-\iint ar^2u_\theta{u}_{1ss}-2\iint ar^2v{u}_{1ss}-\iint ar^2v_\theta{v}_{1ss}+2\iint ar^2u{v}_{1ss}
\\&\leq\frac{\varepsilon^2}{100}(\|c'\|^2_{2}+\|d'\|^2_{2})+\frac{1}{\varepsilon^2}(\|e'\|^2_{2}+\|f'\|^2_{2})+
\frac{\varepsilon^2}{100}(\|c''\|^2_{2}+\|d''\|^2_{2}+\|e''\|^2_{2}+\|f''\|^2_{2}).
\end{align*}

It is easy to obtain that
\begin{align*}
-\iint e^{-s}({u}^a-u_e(r))u_\theta{u}_{1ss}\lesssim\frac{\varepsilon^2}{100}(\|c''\|^2_{2}+\|d''\|^2_{2})+\frac{1}{\varepsilon^2}(\|e'\|^2_{2}+\|f'\|^2_{2}+\|\tilde{u}_\theta\|^2_{2}),
\end{align*}
and
\begin{align*}
-\iint e^{-s}uu^a_\theta{u}_{1ss}\leq\frac{\varepsilon^2}{100}(\|c''\|^2_{2}+\|d''\|^2_{2})+\frac{1}{\varepsilon^2}(\|e'\|^2_{2}+\|f'\|^2_{2}+\|\tilde{u}_\theta\|^2_{2}+
\|u_0\|^2_{2}),
\end{align*}
and
\begin{align*}
-\iint e^{-s}vr(u^a-u_e(r))_r{u}_{1ss}\leq\frac{\varepsilon^2}{100}(\|c''\|^2_{2}+\|d''\|^2_{2})+\frac{1}{\varepsilon^4}(\|e'\|^2_{2}+\|f'\|^2_{2}+\|\tilde{v}_\theta\|^2_{2}),
\end{align*}
and
\begin{align*}
&-\iint e^{-s}v^au{u}_{1ss}\leq\frac{\varepsilon^2}{100}(\|c''\|^2_{2}+\|d''\|^2_{2})+(\|e'\|^2_{2}+\|f'\|^2_{2}+\|\tilde{u}_\theta\|^2_{2}+
\|u_0\|^2_{2}),
\end{align*}
and
\begin{align*}
&-\iint e^{-s}vu^a{u}_{1ss}\leq\frac{\varepsilon^2}{100}(\|c''\|^2_{2}+\|d''\|^2_{2})+\frac{1}{\varepsilon^2}(\|e'\|^2_{2}+\|f'\|^2_{2}+\|\tilde{v}_\theta\|^2_{2}).
\end{align*}

Moreover, based on integration by parts and  Lemma \ref{hardy-for-c-d-e-f} one has
\begin{align*}
&-\iint e^{-s}{v}^aru_r{u}_{1ss}\\&=\iint e^{-s}{v}^au_s{u}_{1ss}\\
&=\iint e^{-s}{v}^a(u_{0s}+\tilde{u}_{s}){u}_{1ss}+\iint e^{-s}{v}^au_{1s}{u}_{1ss}\\
&=\iint e^{-s}{v}^a(u_{0s}+\tilde{u}_{s}){u}_{1ss}+\frac{1}{2}\iint e^{-s}{v}^a\big((u_{1s}u_1)_s-{u}_{1ss}u_1\big)_s
\\
&=\iint e^{-s}{v}^a(u_{0s}+\tilde{u}_{s}){u}_{1ss}+\frac{1}{2}\iint \big(e^{-s}{v}^a\big)_{ss}u_{1s}u_1+\frac{1}{2}\iint \big(e^{-s}{v}^a\big)_{s}{u}_{1ss}u_1
\\
&\leq\frac{\varepsilon^2}{100}(\|c''\|^2_{2}+\|d''\|^2_{2})+\frac{\varepsilon^2}{100}(\|c'\|^2_{2}+\|d'\|^2_{2})+C\big(\|u_{0s}\|^2_{2}+\|\tilde{u}_{s}\|^2_{2}\big)\\&\quad+\frac{1}{\varepsilon^4}(\|e'\|^2_{2}+\|f'\|^2_{2}).
\end{align*}

Thus we get
\begin{align*}
J_2&\leq\frac{\varepsilon^2}{100}(\|c''\|^2_{2}+\|d''\|^2_{2})+\frac{\varepsilon^2}{100}(\|c'\|^2_{2}+\|d'\|^2_{2})
\\&\quad+\frac{1}{\varepsilon^4}(\|u_{0}\|^2_{2}+\|u_{0s}\|^2_{2}+\|\tilde{u}_\theta\|^2_{2}+\|\tilde{v}_\theta\|^2_{2}+\|\tilde{u}_{s}\|^2_{2}+\|e'\|^2_{2}+\|f'\|^2_{2}).
\end{align*}
Similarly, one can get
\begin{align*}
J_3&\leq\frac{\varepsilon^2}{100}(\|c''\|^2_{2}+\|d''\|^2_{2})
+\frac{1}{\varepsilon^2}(\|u_{0}\|^2_{2}+\|\tilde{u}_\theta\|^2_{2}+\|\tilde{v}_\theta\|^2_{2}+\|\tilde{v}_s\|^2_{2}+\|e'\|^2_{2}+\|f'\|^2_{2}).
\end{align*}

Accordingly there holds
\begin{align}
&\iint\big(e^{-s}S_u(-{u}_{1ss})+e^{-s}S_v(-{v}_{1ss})\big)\nonumber\\
&\leq\frac{\varepsilon^2}{100}(\|c''\|^2_{2}+\|d''\|^2_{2})+\frac{\varepsilon^2}{100}(\|c'\|^2_{2}+\|d'\|^2_{2})
\nonumber\\&\quad+\frac{1}{\varepsilon^4}(\|u_{0}\|^2_{2}+\|u_{0s}\|^2_{2}+\|\tilde{u}_\theta\|^2_{2}+\|\tilde{v}_\theta\|^2_{2}+\|\tilde{u}_{s}\|^2_{2}+\|\tilde{v}_{s}\|^2_{2}+\|e'\|^2_{2}+\|f'\|^2_{2}).\label{con-est-25}
\end{align}

{\bf Pressure term:} By integrating by parts and the divergence-free condition we deduce that
\begin{align}\label{e:pressure-estimate-in-positivity-estimate-14}
&\iint e^{-s}p_\theta(-{u}_{1ss})-e^{-s}p_s(-{v}_{1ss})\nonumber\\&=\iint e^{-s}p\big(u_{1ss\theta}+v_{1ss}-v_{1sss}\big)-\int^{2\pi}_0 p(\theta,0)v_{1ss}(\theta,0)d\theta\nonumber\\
&=-\int^{2\pi}_0 p(\theta,0)u_{1s\theta}(\theta,0)d\theta\nonumber\\
&=\int^{2\pi}_0\int_0^{\ln\frac{4}{3}}\partial_s\big(p(\theta,s)\chi_2(s) u_{1s\theta}(\theta,s)\big)d\theta ds,
\end{align}
where we have used $v_{1s}(\theta,0)=v_{1ss}(\theta,0)=0$, $\chi_2(s)$ is a cut-off function defined in $[0,1]$ such that $\chi_2(0)=1$ and $\text{supp}\chi_2\subset[0,\ln\frac{4}{3}]$.

Next, we go back to Euler coordinates and continue to deal  with \eqref{e:pressure-estimate-in-positivity-estimate-14}.
Let
\begin{align*}
\vec{e}_\theta=
\left(\begin{array}{c}
-\sin \theta \\[5pt]
\cos \theta
\end{array}\right), \ \
\vec{e}_r=
\left(\begin{array}{c}
\cos \theta \\[5pt]
\sin \theta
\end{array}\right)
\end{align*}
and
\begin{align*}
\left(\begin{array}{c}
\bar{u}(x,y) \\[5pt]
\bar{v}(x,y)
\end{array}\right)
&=u(\theta, r)\vec{e}_\theta+v(\theta,r)\vec{e}_r, \ \
\left(\begin{array}{c}
\bar{u}^a(x,y) \\[5pt]
\bar{v}^a(x,y)
\end{array}\right)
=u^a(\theta, r)\vec{e}_\theta+v^a(\theta,r)\vec{e}_r,\\
\bar{p}(x,y)&=p(\theta,r), \ \
\left(\begin{array}{c}
\bar{F}_u(x,y) \\[5pt]
\bar{F}_v(x,y)
\end{array}\right)
=\frac{F_u(\theta, r)}{r}\vec{e}_\theta+\frac{F_v(\theta, r)}{r}\vec{e}_r.
\end{align*}
Then the error equations (\ref{linear-equation-of-error}) can be written in Euler coordinates as follows
\begin{align}\label{error-equation-in-Euler-coordinates}
\left\{
\begin{array}{lll}
\bar{u}^a\partial_x \bar{u}+\bar{v}^a\partial_y \bar{u} + \bar{u}\partial_x \bar{u}^a +\bar{v}\partial_y \bar{u}^a
+\partial_x\bar{p}-\varepsilon^2 \Delta \bar{u}= \bar{F}_u,\\[5pt]
\bar{u}^a\partial_x \bar{v}+\bar{v}^a\partial_y \bar{v} + \bar{u}\partial_x \bar{v}^a +\bar{v}\partial_y \bar{v}^a
+\partial_y\bar{p}-\varepsilon^2 \Delta \bar{v}= \bar{F}_v,\\[5pt]
\partial_x\bar{u}+\partial_y \bar{v}=0,\\[5pt]
(\bar{u},\bar{v})|_{\partial B_1}=(0,0).
\end{array}
\right.
\end{align}

It is noted that
\begin{align}
&\int^{2\pi}_0\int_0^{\ln\frac{4}{3}}\partial_s\big(p(\theta,s)\chi_2 u_{1s\theta}(\theta,s)\big)d\theta ds\nonumber\\&=\iint
\partial_sp\chi_2 u_{1s\theta}+p(\theta,s)\chi'_2 u_{1s\theta}+p(\theta,s)\chi_2 u_{1ss\theta}\nonumber\\
&=\iint\partial_sp\chi_2 u_{1s\theta}-\partial_\theta p(\theta,s)\chi'_2 u_{1s}-\partial_\theta p(\theta,s)\chi_2 u_{1ss}\nonumber\\
&\lesssim\|(c',d')\|^2_{L^2((0,\ln\frac{4}{3}))}+\frac{\varepsilon^2}{100}\|(c'',d'')\|^2_{L^2((0,\ln\frac{4}{3}))}+\frac{1}{\varepsilon^2}\|(p_\theta,p_s)\|^2_{L^2((0,2\pi)\times(0,\ln\frac{4}{3}))}
\nonumber\\
&\lesssim\|(\nabla\bar{u},\nabla\bar{v})\|^2_{L^2(V)}+\frac{\varepsilon^2}{100}\|(c'',d'')\|^2_{L^2((0,\ln\frac{4}{3}))}+\frac{1}{\varepsilon^2}\|\nabla\bar{p}\|^2_{L^2(V)},\label{e:pressure-estimate-in-positivity-estimate-15}
\end{align}
where $V=\{\frac{3}{4}\leq\sqrt{x^2+y^2}\leq1\}$.

In terms of the estimate for the  Stokes equations (e.g., \cite{GP} and \cite{Tasi}) and the trace theorem, we get
\begin{align}\label{stokes-estimate-3}
  \|\varepsilon^2\nabla^2(\bar{u},\bar{v})\|_{L^2(V)}+\|\nabla\bar{p}\|_{L^2(V)}
  \lesssim\| \mathbf{F}\|_{L^2(V')}+\varepsilon^{2}\|(\bar{u},\bar{v})\|_{L^2(V')},
\end{align}
where $V\subset V'=\{\frac{2}{3}\leq\sqrt{x^2+y^2}\leq1\}$ and
both annuli stay a fixed positive distance away from the singular
origin.  Hence the standard local Stokes estimate and the trace theorem
apply without any singular-core correction; the boundary terms produced by
the first-mode integrations by parts are controlled by the quantities in the
$Y$-norm and by the source norms on $V'$.
\begin{equation}\label{def:F}
 \mathbf{F}=(\bar{u}^a\partial_x\bar{u}+\bar{v}^a\partial_y \bar{u} + \bar{u}\partial_x \bar{u}^a +\bar{v}\partial_y \bar{u}^a-\bar{F}_u,
\bar{u}^a\partial_x \bar{v}+\bar{v}^a\partial_y \bar{v} + \bar{u}\partial_x \bar{v}^a +\bar{v}\partial_y \bar{v}^a-\bar{F}_v).
\end{equation}

To proceed we need to estimate $\|\nabla(\bar{u},\bar{v})\|_{L^2(V')}$. For that we write(e.g., \cite{GP})
\begin{align*}
 \mathbf{F}=&\big(\partial_x(\bar{u}^a\bar{u}+\bar{v}\bar{u}^a)+\partial_y(\bar{v}^a\bar{u}
 +\bar{v}\bar{u}^a)+\text{div}\phi_1,\\& \partial_x(\bar{u}^a\bar{v}+\bar{u}\bar{v}^a)+\partial_y(\bar{v}^a\bar{v}+\bar{v}\bar{v}^a)
 +\text{div}\phi_2\big),\ (x,y)\in V'',
\end{align*}
where  $V'\subset V''=\{\frac{1}{2}\leq\sqrt{x^2+y^2}\leq1\}$ and
\begin{align}
\|\phi_1\|_{H^1(V'')}\lesssim\|\bar{F}_u\|_{L^2(V'')},\ \ \  \|\phi_2\|_{H^1(V'')}\lesssim\|\bar{F}_v\|_{L^2(V'')}.\label{phi-est}
\end{align}

It follows from the  estimate for the Stokes equations(e.g., \cite{GP}  and  \cite{Tasi}) and the trace theorem, we have
\begin{align}
\|\varepsilon^2\nabla(\bar{u},\bar{v})\|_{L^2(V')}&\lesssim\|\bar{u}^a\bar{u}+\bar{v}\bar{u}^a+\bar{v}^a\bar{u}
 +\bar{v}\bar{u}^a+\phi_1\|_{L^2(V'')}
 \nonumber\\ &\quad +\|\bar{u}^a\bar{v}+\bar{u}\bar{v}^a+\bar{v}^a\bar{v}+\bar{v}\bar{v}^a+\phi_2\|_{L^2(V'')}
+\|\varepsilon^2(\bar{u},\bar{v})\|_{L^2(V'')}.\label{stokes-estimate-4}
\end{align}

\par

  It follows from \eqref{weight-Hardy}, \eqref{stokes-estimate-3} and \eqref{def:F} and that
\begin{align}
\|\nabla\bar{p}\|^2_{L^2(V)}&\lesssim \| \mathbf{F}\|^2_{L^2(V')}+\|\varepsilon^2(\bar{u},\bar{v})\|^2_{L^2(V')}\nonumber\\&\lesssim\|(c',d',e',f',u_0)\|^2_{L^2((0,\ln 2))}+\|(u_{0s},\tilde{u}_\theta,\tilde{v}_\theta,\tilde{u}_s,\tilde{v}_s)\|^2_{L^2((0,2\pi)\times(0,\ln 2))}\nonumber\\&\quad+\|\big(F_u,F_v\big)\|^2_{L^2((0,2\pi)\times(0,\ln 2))}.\label{pre-est-der}
\end{align}

Moreover, by using \eqref{phi-est}, \eqref{stokes-estimate-4} and direct computations one has
\begin{align*}
&\varepsilon^2\|(c',d')\|^2_{L^2((0,\ln 2))}\nonumber\\&\lesssim\varepsilon^2\|\nabla(\bar{u},\bar{v})\|_{L^2(V')}\nonumber
\\&\lesssim\|\bar{u}^a\bar{u}+\bar{v}\bar{u}^a+\bar{v}^a\bar{u}
 +\bar{v}\bar{u}^a+\phi_1\|_{L^2(V'')}+\|\bar{u}^a\bar{v}+\bar{u}\bar{v}^a+\bar{v}^a\bar{v}+\bar{v}\bar{v}^a+\phi_2\|_{L^2(V'')}\nonumber\\&\quad+\|\varepsilon^2(\bar{u},\bar{v})\|_{L^2(V'')}
\nonumber\\&\lesssim\|(e',f',u_0)\|^2_{L^2((0,\ln 2))}+\|(\tilde{u}_\theta,\tilde{v}_\theta)\|^2_{L^2((0,2\pi)\times(0,\ln 2))}+\|(F_u,F_v)\|^2_{L^2((0,2\pi)\times(0,\ln 2))}
\end{align*}
which, together with \eqref{pre-est-der}, implies
\begin{align*}
\|\nabla\bar{p}\|^2_{L^2(V)}&\lesssim \varepsilon^{-2}\|(e',f',u_0)\|^2_{L^2((0,\ln 2))}+ \varepsilon^{-2}\|(u_{0s},\tilde{u}_\theta,\tilde{v}_\theta,\tilde{u}_s,\tilde{v}_s)\|^2_{L^2((0,2\pi)\times(0,\ln 2))}\nonumber\\&\quad+\varepsilon^{-2}\|e^{-s}F_u,e^{-s}F_v\|^2_{L^2((0,2\pi)\times(0,\ln 2))}.
\end{align*}

Therefore we arrive at
\begin{align*}
&\int^{2\pi}_0\int_0^{\ln\frac{4}{3}}\partial_s\big(p(\theta,s)\chi u_{1s\theta}(\theta,s)\big)d\theta ds\nonumber\\&\lesssim\varepsilon^{-4}\|(e',f',u_0)\|^2_{L^2((0,\ln 2))}+ \varepsilon^{-4}\|(u_{0s},\tilde{u}_\theta,\tilde{v}_\theta,\tilde{u}_s,\tilde{v}_s)\|^2_{L^2((0,2\pi)\times(0,\ln 2))}\nonumber\\&\quad+\varepsilon^{-4}\|e^{-s}F_u,e^{-s}F_v\|^2_{L^2((0,2\pi)\times(0,\ln 2))}+\frac{\varepsilon^2}{100}\|(c'',d'')\|^2_{L^2((0,\ln 2))}
\end{align*}
which, together with \eqref{e:pressure-estimate-in-positivity-estimate-14}, leads to
\begin{align}\label{e:pressure-estimate-in-positivity-estimate-24}
&\iint e^{-s}p_\theta(-{u}_{1ss})-e^{-s}p_s(-{v}_{1ss})\nonumber\\&\lesssim\varepsilon^{-4}\|(e',f',u_0)\|^2_{L^2((0,\infty))}+ \varepsilon^{-4}\|(u_{0s},\tilde{u}_\theta,\tilde{v}_\theta,\tilde{u}_s,\tilde{v}_s)\|^2_{L^2((0,2\pi)\times(0,\infty))}\nonumber\\&\quad+\varepsilon^{-4}\|e^{-s}F_u,e^{-s}F_v\|^2_{L^2((0,2\pi)\times(0,\ln 2))}+\frac{\varepsilon^2}{100}\|(c'',d'')\|^2_{L^2((0,\infty))}.
\end{align}

Collecting the estimates (\ref{cj1}), \eqref{con-est-25} and (\ref{e:pressure-estimate-in-positivity-estimate-24}) together, we obtain (\ref{eq:H2-estimate})  and then complete the proof of this Lemma.
\end{proof}

In order to estimate $L^\infty$, we also need to estimate the terms $\|\tilde{u}_{s\theta}\|_2$ and $\|\tilde{v}_{s\theta}\|_2$.
\begin{lemma}\label{second-order-der-high-fre}
If $(u,v)$ satisfies  \eqref{rewrite-linear-equation-of-error} and the hypothesis {\bf (H)}, then
\begin{align}\label{eq:H2-estimate-high}
&\varepsilon^2\int_{0}^{+\infty}\int_0^{2\pi}\big(\tilde{u}^2_{s\theta}+{\tilde{u}}^2_{\theta\theta}+\tilde{v}^2_{s\theta}+{\tilde{v}}^2_{\theta\theta}\big)d\theta ds
\nonumber\\
&\leq C\|u_{0}\|^2_{2}+C\varepsilon^2\big(\|u_{0s}\|^2_{2}+\|\tilde{u}_s\|^2_{2}+\|\tilde{v}_s\|^2_{2}
+\|c'\|^2_{2}+\|d'\|^2_{2}\big)\nonumber\\&\quad+C\varepsilon^{-2}\big(\|\tilde{u}_{\theta}\|^2_{2}+\|\tilde{v}_{\theta}\|^2_{2}+\|e'\|^2_{2}+\|f'\|^2_{2}\big)\nonumber\\&\quad+\Big|\int_{0}^{+\infty}\int_0^{2\pi}\big(e^{-s}F_u\tilde{u}_{\theta\theta}+ e^{-s}F_v\tilde{v}_{\theta\theta})\big)d\theta ds\Big|.
\end{align}
\end{lemma}
\begin{proof}
Multiplying the first equation of \eqref{rewrite-linear-equation-of-error} by $-\tilde{u}_{\theta\theta}$ and the second equation by $-\tilde{v}_{\theta\theta}$, integrating in $\bar{\Omega}$ and adding the resultant together, we obtain
\begin{align*}
&\iint\varepsilon^2\big(u_{ss}
+u_{\theta\theta}+2v_{\theta}-u\big){\tilde{u}}_{\theta\theta}+\varepsilon^2\big(v_{ss}+v_{\theta\theta}-2u_{\theta}-v\big)\tilde{v}_{\theta\theta}-\iint\big(e^{-s}S_u\tilde{u}_{\theta\theta}+e^{-s}S_v\tilde{v}_{\theta\theta}\big)\\&=-\iint\big(e^{-s}F_u\tilde{u}_{\theta\theta}+ e^{-s}F_v\tilde{v}_{\theta\theta}\big),
\end{align*}
where we have used  the divergence-free condition
$\tilde{v}_s-\tilde{v}-\tilde{u}_\theta=0$.

{\bf Diffusion term:}
\begin{align}\label{cj2}
&\iint\varepsilon^2\big(u_{ss}
+u_{\theta\theta}+2v_{\theta}-u\big){\tilde{u}}_{\theta\theta}+\varepsilon^2\big(v_{ss}+v_{\theta\theta}-2u_{\theta}-v\big)\tilde{v}_{\theta\theta}
\nonumber\\
&=\varepsilon^2\iint\tilde{u}^2_{s\theta}+{\tilde{u}}^2_{\theta\theta}+{\tilde{u}}^2_{\theta}+\tilde{v}^2_{s\theta}+{\tilde{v}}^2_{\theta\theta}+{\tilde{v}}^2_{\theta}+2\varepsilon^2\iint
\tilde{v}_{\theta}{\tilde{u}}_{\theta\theta}-\tilde{u}_{\theta}\tilde{v}_{\theta\theta}\nonumber\\
&\geq\frac{\varepsilon^2}{2}\iint\big(\tilde{u}^2_{s\theta}+{\tilde{u}}^2_{\theta\theta}+\tilde{v}^2_{s\theta}+{\tilde{v}}^2_{\theta\theta}\big)-C\varepsilon^2\iint\big({\tilde{u}}^2_{\theta}+{\tilde{v}}^2_{\theta}\big).
\end{align}
{\bf Convective term:}
\begin{align}
&\iint\big(e^{-s}S_u({\tilde{u}}_{\theta\theta})+e^{-s}S_v(\tilde{v}_{\theta\theta})\big)\nonumber\\
&=\underbrace{\iint e^{-s}\big((ar+\frac{b}{r})u_\theta+v(ar-\frac{b}{r})+v(ar+\frac{b}{r})\big){\tilde{u}}_{\theta\theta}
+\iint e^{-s}\big((ar+\frac{b}{r})v_\theta-(2ar+\frac{2b}{r})u\big)\tilde{v}_{\theta\theta}}_{K_1}\nonumber\\
&\quad+\underbrace{\iint\big(({u}^a-u_e(r))u_\theta+{v}^aru_r+uu^a_\theta+vr(u^a-u_e(r))_r+{v}^au+v(u^a-u_e(r))\big) e^{-s}{\tilde{u}}_{\theta\theta}}_{K_2}\nonumber\\
&\quad+\underbrace{\iint\big((u^a-u_e(r))v_\theta+{v}^arv_r+u{v}^a_\theta+vr{v}^a_r-2u(u^a-u_e(r))\big)e^{-s}\tilde{v}_{\theta\theta}}_{K_3}.\label{Convective-term-low}
\end{align}

 By integration by part, we obtain
\begin{align}
K_1=\iint2ar^2\tilde{v}\tilde{u}_{\theta\theta}-2\iint (ar^2+b)\tilde{u}{\tilde{v}}_{\theta\theta}
\lesssim \|\tilde{u}_\theta\|^2_{2}+\|\tilde{v}_\theta\|^2_{2}.\label{Convective-term-low-1}
\end{align}

It follows from the (\ref{estimate-on-Euler-parts}), \eqref{estimate-on-Prandtl-parts}, \eqref{decompose-u}, \eqref{decompose-u2} and Lemma \ref{hardy-for-c-d-e-f} that
\begin{align*}
\iint e^{-s}({u}^a-u_e(r))u_\theta{\tilde{u}}_{\theta\theta}&=\iint e^{-s}({u}^a-u_e(r))(u_1+\tilde{u})_\theta{\tilde{u}}_{\theta\theta}\\
&\leq C\varepsilon^{-2}\big(\|e'\|^2_{L^2}+\|f'\|^2_{2}\big)+\frac{\varepsilon^2}{100}\|\tilde{u}_{\theta\theta}\|^2_{2},
\end{align*}
and
\begin{align*}
\iint e^{-s}v^aru_r{\tilde{u}}_{\theta\theta}&=\iint e^{-s}v^a_{\theta}u_s{\tilde{u}}_{\theta}+
\iint e^{-s}v^au_{s\theta}{\tilde{u}}_{\theta}\\
&\leq C\big(\varepsilon^2\|u_{0s}\|^2_{2}+\varepsilon^2\|c'\|^2_{2}+\varepsilon^2\|d'\|^2_{2}+\varepsilon^2\|\tilde{u}_s\|^2_{2}+\|\tilde{u}_{\theta}\|^2_{2}\big)\\
&\ \ \ +\frac{\varepsilon^2}{100}\|\tilde{u}_{s\theta}\|^2_{2},
\end{align*}
and
\begin{align*}
\iint e^{-s}uu^a_\theta{\tilde{u}}_{\theta\theta}&=-\iint e^{-s}u_{\theta}u^a_\theta{\tilde{u}}_{\theta}-
\iint e^{-s}uu^a_{\theta\theta}{\tilde{u}}_{\theta}\\
&\lesssim\|u_{0}\|^2_{2}+\|e'\|^2_{2}+\|f'\|^2_{2}+\|\tilde{u}_{\theta}\|^2_{2},
\end{align*}
and
\begin{align*}
\iint e^{-s}vr(u^a-u_e(r))_r{\tilde{u}}_{\theta\theta}\leq C\varepsilon^{-2}\big(\|\tilde{v}_{\theta}\|^2_{2}+\|e'\|^2_{2}+\|f'\|^2_{2}+\|\tilde{u}_{\theta}\|^2_{2}\big)+\frac{\varepsilon^2}{100}\|\tilde{u}_{\theta\theta}\|^2_{2},
\end{align*}
and
\begin{align*}
\iint e^{-s}v^au{\tilde{u}}_{\theta\theta}&=-\iint e^{-s}v^a_{\theta}u{\tilde{u}}_{\theta}-
\iint e^{-s}v^a u_\theta{\tilde{u}}_{\theta}\nonumber\\
&\leq C\big(\|u_{0}\|^2_{2}+\|e'\|^2_{2}+\|f'\|^2_{2}+\|\tilde{u}_{\theta}\|^2_{2}\big)+\frac{\varepsilon^2}{100}\|\tilde{u}_{\theta\theta}\|^2_{2},
\end{align*}
and
\begin{align*}
\iint e^{-s}v(u^a-u_e(r)){\tilde{u}}_{\theta\theta}&=-\iint e^{-s}v_\theta(u^a-u_e(r)){\tilde{u}}_{\theta}-
\iint e^{-s}v(u^a-u_e(r))_{\theta}{\tilde{u}}_{\theta}\\
&\lesssim \|e'\|^2_{2}+\|f'\|^2_{2}+\|\tilde{u}_{\theta}\|^2_{2}+\|\tilde{v}_{\theta}\|^2_{2}.
\end{align*}

Accordingly one has
\begin{align}
K_2&\leq C\big(\|u_{0}\|^2_{2}+\varepsilon^2\|u_{0s}\|^2_{2}+\varepsilon^2\|c'\|^2_{2}+\varepsilon^2\|d'\|^2_{2}+\varepsilon^2\|\tilde{u}_s\|^2_{2}\big)\nonumber\\
&\quad+C\varepsilon^{-2}\big(\|\tilde{v}_{\theta}\|^2_{2}+\|e'\|^2_{2}+\|f'\|^2_{2}+\|\tilde{u}_{\theta}\|^2_{2}\big)+\frac{\varepsilon^2}{100}\|\tilde{v}_{\theta\theta}\|^2_{2}.\label{Convective-term-low-2}
\end{align}

Similarly we have
\begin{align}
K_3&\leq C\big(\|u_{0}\|^2_{2}+\|\tilde{u}_{\theta}\|^2_{2}+\varepsilon^2\|\tilde{v}_s\|^2_{2}\big)+C\varepsilon^{-2}\big(\|e'\|^2_{2}+\|f'\|^2_{2}+\|\tilde{v}_{\theta}\|^2_{2}\big)
\nonumber\\ &\quad+\frac{\varepsilon^2}{100}\|\tilde{v}_{\theta\theta}\|^2_{2}+\frac{\varepsilon^2}{100}\|\tilde{v}_{\theta s}\|^2_{2}.\label{Convective-term-low-3}
\end{align}

Substituting \eqref{Convective-term-low-1}-\eqref{Convective-term-low-3} into \eqref{Convective-term-low} we have
\begin{align}
&\iint e^{-s}\big(S_u{\tilde{u}}_{\theta\theta}+S_v\tilde{v}_{\theta\theta}\big)\nonumber\\ &\leq C \|u_{0}\|^2_{2}+\varepsilon^2\|u_{0s}\|^2_{2}+\varepsilon^2\|c'\|^2_{2}+\varepsilon^2\|d'\|^2_{2}+\varepsilon^2\|\tilde{u}_s\|^2_{2}+\varepsilon^2\|\tilde{v}_s\|^2_{2}\nonumber\\
&\quad+\varepsilon^{-2}\big(\|\tilde{v}_{\theta}\|^2_{2}+\|e'\|^2_{2}+\|f'\|^2_{2}+\|\tilde{u}_{\theta}\|^2_{2}\big)+\frac{\varepsilon^2}{100}\|\tilde{v}_{\theta\theta}\|^2_{2}+\frac{\varepsilon^2}{100}\|\tilde{v}_{\theta s}\|^2_{2}.\label{Convective-term-low-4}
\end{align}

 Combining the estimates \eqref{cj2} and \eqref{Convective-term-low-4} together, we obtain (\ref{eq:H2-estimate-high}) which completes the proof of this lemma.
\end{proof}
We eliminate the pressure of the system \eqref{rewrite-linear-equation-of-error} and then get
\begin{align}\label{c-d-e-f}
&-\varepsilon^2e^{s}\big(u_{sss}
+u_{s\theta\theta}+2v_{s\theta}-u_{s}+u_{ss}+u_{\theta\theta}+2v_{\theta}-u\nonumber\\
&\quad+v_{ss\theta}+v_{\theta\theta\theta}-2u_{\theta\theta}-v_{\theta}\big)
+(S_u)_s+(S_v)_\theta\nonumber\\
&=(F_u)_s+(F_v)_\theta.
\end{align}
\begin{lemma}\label{L2-one-fre}
If $(u,v)$ satisfies  \eqref{rewrite-linear-equation-of-error}, the hypothesis {\bf (H)} and
\begin{align}\label{apri-est}e^{\frac{s}{2}}\big(c-c_\infty,d-d_\infty,e-e_\infty,f-f_\infty,c',d',e',f'\big)\in L^2((0,+\infty)),\end{align}
then
\begin{align}\label{eq:H2-estimate-first}
&\int_0^\infty e^s\big((c-c_\infty)^2+(d-d_\infty)^2+(e-e_\infty)^2+(f-f_\infty)^2\big)ds
\nonumber\\&\quad+\int_0^\infty e^s(c'^2+d'^2+e'^2+f'^2)ds+c^2_{\infty}+d^2_{\infty}
\nonumber\\
&\lesssim\varepsilon^2(\|c''\|^2_{2}+\|d''\|^2_{2}+\|e''\|^2_{2}+\|f''\|^2_{2})
+\|u_{0s}\|^2_{2}+\|\tilde{u}_s\|^2_{2}+
\|\tilde{v}_s\|^2_{L^2}\nonumber\\&\quad+\frac{1}{\varepsilon^2}\big(\|u_{0}\|^2_{2}+\|\tilde{u}_{\theta}\|^2_{2}\big)+
\frac{1}{\varepsilon^4}\big(\|e'\|^2_{2}+\|f'\|^2_{2}+\|\tilde{v}_{\theta}\|^2_{2}\big)\nonumber\\&\quad+\Big|\int_{0}^{+\infty}\int_0^{2\pi} (\partial_sF_u+\partial_\theta F_v)u_{1s\theta}d\theta ds\Big|+\|e^{-s}F_u\|_{2}+\|e^{-s}F_v\|_{2}.
\end{align}
\end{lemma}
\begin{proof}
Multiplying  \eqref{c-d-e-f} by $u_{1s\theta}$, we obtain
\begin{align}
&\iint-\varepsilon^2e^{s}\big(u_{sss}
+u_{s\theta\theta}+2v_{s\theta}-u_{s}+u_{ss}-u+v_{ss\theta}+v_{\theta\theta\theta}-u_{\theta\theta}+v_{\theta}\big)u_{1s\theta}\nonumber\\
&\ \ \ +\iint\big((S_u)_s+(S_v)_\theta\big)u_{1s\theta}=\iint \big(\partial_sF_u+\partial_\theta F_v\big)u_{1s\theta}.\label{dif-con-term}
\end{align}

{\bf Diffusion term:} Due to direct computations we find
\begin{align}\label{cj10}
&-\varepsilon^2\iint e^{s}\big(u_{sss}
+u_{s\theta\theta}+2v_{s\theta}-u_{s}+u_{ss}-u+v_{ss\theta}+v_{\theta\theta\theta}-u_{\theta\theta}+v_{\theta}\big)u_{1s\theta}
\nonumber\\&=\varepsilon^2\int_0^{2\pi} u_{1ss}(\theta,0)u_{1s\theta}(\theta,0)d\theta-3\varepsilon^2\iint e^{s}v_{s\theta}u_{1s\theta}
\nonumber\\&=-\varepsilon^2\iint \partial_s\big(\chi_2u_{1ss}(\theta,s)u_{1s\theta}(\theta,s)\big)-3\varepsilon^2\iint e^{s}v_{s\theta}u_{1s\theta}.
\end{align}

According to  \eqref{c-d-e-f}, one has
  \begin{align}\label{cj10-10}
&-\varepsilon^2\iint \partial_s\big(\chi_2u_{1ss}(\theta,s)u_{1s\theta}(\theta,s)\big)
\nonumber\\&=-\varepsilon^2\iint \chi_2'u_{1ss}u_{1s\theta}+\chi_2u_{1sss}u_{1s\theta}
\nonumber\\&\geq -C\varepsilon^2(\|c''\|^2_{2}+\|d''\|^2_{2}+\|e''\|^2_{2}+\|f''\|^2_{2})-\frac{1}{100}\int_0^\infty e^s(c'^2+d'^2)ds-\frac{C}{\varepsilon^4}(\|e'\|^2_{L^2}+\|f'\|^2_{L^2}+\|\tilde{v}_\theta\|^2_{L^2})\nonumber\\& \quad-C\eta\big(\|{u}_{0s}\|^2_{2}+\|\tilde{u}_s\|^2_{2}+\|\tilde{v}_s\|^2_{2}+\|\tilde{u}_\theta\|^2_{2}\big)-C\big(\|e^{-s}F_u\|_{2}+\|e^{-s}F_v\|_{2}\big) .
\end{align}
Moreover,
\begin{align*}
-3\varepsilon^2\iint e^{s}v_{s\theta}u_{1s\theta}\gtrsim -\varepsilon^2\int_0^\infty e^s(c'^2+d'^2+e'^2+f'^2)ds.
\end{align*}

Therefore we have
\begin{align}\label{cj10-11}
&-\varepsilon^2\iint e^{s}\big(u_{sss}
+u_{s\theta\theta}+2v_{s\theta}-u_{s}+u_{ss}-u+v_{ss\theta}+v_{\theta\theta\theta}-u_{\theta\theta}+v_{\theta}\big)u_{1s\theta}
\nonumber\\&\geq -C\varepsilon^2(\|c''\|^2_{2}+\|d''\|^2_{2}+\|e''\|^2_{2}+\|f''\|^2_{2})-\frac{1}{100}\int_0^\infty e^s(c'^2+d'^2+e'^2+f'^2)ds\nonumber\\&\quad-\frac{C}{\varepsilon^4}(\|e'\|^2_{2}+\|f'\|^2_{2}+\|\tilde{v}_\theta\|^2_{2})-C\eta\big(\|{u}_{0s}\|^2_{2}+\|\tilde{u}_s\|^2_{2}+\|\tilde{v}_s\|^2_{2}+\|\tilde{u}_\theta\|^2_{2}\big)
\nonumber\\& \quad-C\big(\|e^{-s}F_u\|_{2}+\|e^{-s}F_v\|_{2}\big) .
\end{align}

{\bf Convective term:}
We divide the convective term into four parts as follows:
\begin{align}
&\iint\big((S_u)_s+(S_v)_\theta\big)u_{1s\theta}\nonumber\\
&=\underbrace{b\iint e^s(u_{s\theta}+v_{\theta\theta}-u_{\theta})u_{1s\theta}}_{K_4}+\underbrace{a\iint(-3e^{-s}u_\theta
+e^{-s}u_{\theta s}-2e^{-s}v+2e^{-s}v_s+e^{-s}v_{\theta\theta})u_{1s\theta}}_{K_{5}}\nonumber\\
&\quad+\underbrace{\iint\big(({u}^a-u_e(r))u_\theta+{v}^aru_r+u(u^a)_\theta+vr(u^a-u_e(r))_r+{v}^au+v(u^a-u_e(r))\big)_s u_{1s\theta}}_{K_6}\nonumber\\
&\quad+\underbrace{\iint\big((u^a-u_e(r))v_\theta+{v}^arv_r+u{v}^a_\theta+vr{v}^a_r-2u(u^a-u_e(r))\big)_\theta u_{1s\theta}}_{K_7}.\label{convective-term-11}
\end{align}
{\bf 1)Estimate of $K_4$ and $K_5$.}
First,
\begin{align*}
\iint e^{s}u_{s\theta}u_{1s\theta}=\pi\int_0^\infty e^s(c'^2+d'^2)ds.
\end{align*}

Next, it follows from \eqref{divergence-free-condition-1}, \eqref{apri-est} and integration by parts that
\begin{align}
&\int_0^\infty e^s\big((c-c_\infty)^2+(d-d_\infty)^2\big)ds\nonumber\\
&=\int_0^\infty e^s\bigg(\big(f'-(f-f_\infty)\big)^2+\big(e'-(e-e_\infty)\big)^2\bigg)ds\nonumber\\
&=\int_0^\infty e^s\big(e'^2+f'^2+2(e-e_\infty)^2+2(f-f_\infty)^2\big)ds+e^2_{\infty}+f^2_{\infty}
\nonumber\\
&=\int_0^\infty e^s\big(e'^2+f'^2+2(e-e_\infty)^2+2(f-f_\infty)^2\big)ds+c^2_{\infty}+d^2_{\infty}.\label{apri-est-1}
\end{align}
Moreover, one has
\begin{align*}
-\iint e^s(u_1-u_{1\infty})_{\theta}u_{1s\theta}&=-\iint e^s(u_1-u_{1\infty})_{\theta}(u_1-u_{1\infty})_{s\theta}\\
&=\iint\frac{e^s}{2}\big((u_1-u_{1\infty})_{\theta}\big)^2d\theta ds+\frac{1}{2}\int_0^{\pi}\big((u_{1\infty})_\theta\big)^2 d\theta\\
&=\frac{\pi}{2}\int_0^\infty e^s\big((c-c_\infty)^2+(d-d_\infty)^2\big)ds+\frac{\pi}{2}(c^2_{\infty}+d^2_{\infty})\\
&=\frac{\pi}{2}\int_0^\infty e^s\big(e'^2+f'^2+2(e-e_\infty)^2+2(f-f_\infty)^2\big)ds+\pi(c^2_{\infty}+d^2_{\infty}),
\end{align*}
where we have used \eqref{apri-est-1}. And
we have
\begin{align*}
-\iint e^s(v_1-v_{1\infty})u_{1s\theta}&=\iint e^s(v_1-v_{1\infty})_\theta u_{1s}
\nonumber\\
&\geq -\pi\int_0^\infty \big(|e^s(e-e_\infty)d'|+|(f-f_\infty)c'|\big) ds\nonumber\\
&\geq -\frac{\pi}{2}\int_0^\infty e^s
\big(c'^2+d'^2+(e-e_\infty)^2+(f-f_\infty)^2\big)ds.
\end{align*}
Thus
\begin{align}
K_4\geq\frac{b\pi}{2}\int_0^\infty e^s\big(c'^2+d'^2+e'^2+f'^2+(e-e_\infty)^2+(f-f_\infty)^2\big)ds+b\pi(c^2_{\infty}+d^2_{\infty}),\label{K4-est}
\end{align}
where we have used
\begin{align*}
v_{1\theta\theta}-u_{1\theta}=-v_1-u_{1\theta}=-(v_1-v_{1\infty})-(u_1-u_{1\infty})_{\theta}.
\end{align*}

Furthermore,
\begin{align}K_{5}&=a\iint(-3e^{-s}u_\theta
+e^{-s}u_{\theta s}-2e^{-s}v+2e^{-s}v_s+e^{-s}v_{\theta\theta})u_{1s\theta}
\nonumber\\
&=a\iint(-3e^{-s}u_{1\theta}
-2e^{-s}v_1+2e^{-s}v_{1s}+e^{-s}v_{1\theta\theta})u_{1s\theta}\nonumber\\
&\quad-\frac{a}{2}\iint e^{-s}(u_{1\theta})^2-a\iint e^{-s}u_{1\theta}u_{1ss\theta}
\nonumber\\&\geq-\frac{C}{\varepsilon^2}(\|e'\|^2_{2}+\|f'\|^2_{2})-{\varepsilon^2}(\|e''\|^2_{2}+\|f''\|^2_{2}).\label{K5-est}
\end{align}

Finally, we then have
\begin{align}
K_4+K_5&\gtrsim\int_0^\infty e^s\big(c'^2+d'^2+e'^2+f'^2+(c-c_\infty)^2 +(d-d_\infty)^2+(e-e_\infty)^2+(f-f_\infty)^2\big)ds\nonumber\\
&\quad+c^2_{\infty}+d^2_{\infty}-\frac{C}{\varepsilon^2}(\|e'\|^2_{2}+\|f'\|^2_{2})-{\varepsilon^2}(\|e''\|^2_{2}+\|f''\|^2_{2}).\label{K45-est}
\end{align}
{\bf 2)Estimate of $K_6$.}
First, using integration by parts, (\ref{estimate-on-Euler-parts}), \eqref{estimate-on-Prandtl-parts}, \eqref{decompose-u}, \eqref{decompose-u2} and Lemma \ref{hardy-for-c-d-e-f}, we obtain
\begin{align*}
&\iint \big(({u}^a-u_e(r))u_\theta\big)_s u_{1s\theta}\nonumber\\
&=\iint ({u}^a-u_e(r))_su_\theta u_{1s\theta}+\iint({u}^a-u_e(r))u_{s\theta}u_{1s\theta}
\nonumber\\
&=\iint ({u}^a-u_e(r))_s(u_{1\theta}+\tilde{u}_{\theta}) u_{1s\theta}+\iint({u}^a-u_e(r))(u_{1s\theta}+\tilde{u}_{s\theta})u_{1s\theta}
\\&\lesssim\eta\|e^{\frac{s}{2}}(c',d')\|^2_{2}+\frac{\eta}{\varepsilon^2}\|(e',f',\tilde{u}_\theta)\|^2_{2}
-\iint \big(({u}^a-u_e(r))\tilde{u}_{\theta}\big) u_{1ss\theta}-\iint \big(({u}^a-u_e(r))_s\tilde{u}_{\theta}\big) u_{1s\theta}\\
&\lesssim\eta\|e^{\frac{s}{2}}(c',d')\|^2_{2}+\frac{1}{\varepsilon^2}\|(e',f',\tilde{u}_\theta)\|^2_{2}+\varepsilon^2(\|c''\|^2_{2}+\|d''\|^2_{2}).
\end{align*}

\par
Similarly, one has
\begin{align*}
\iint (v^aru_r)_su_{1s\theta}&=\iint v^au_su_{1ss\theta}\nonumber\\&=\iint v^a(u_{0s}+u_{1s}+\tilde{u}_{s})u_{1ss\theta}
\nonumber\\
&\lesssim\varepsilon^2(\|c''\|^2_{2}+\|d''\|^2_{2})+\eta^2\|c'\|^2_{2}+\eta^2\|d'\|^2_{2}+\eta^2\|\tilde{u}_s\|^2_{2}+\eta^2\|u_{0s}\|^2_{2}
\end{align*}
and
\begin{align*}
\iint(uu^a_\theta)_su_{1s\theta}&=\iint u_su^a_\theta u_{1s\theta}+\iint u(u^a_\theta)_su_{1s\theta}
\nonumber\\&=\iint (u_{0s}+u_{1s}+\tilde{u}_{s})u^a_\theta u_{1s\theta}+\iint (u_{0}+u_{1}+\tilde{u})(u^a_\theta)_su_{1s\theta}\\
&\lesssim\eta\big(\|e^{\frac{s}{2}}c'\|^2_{2}+\|e^{\frac{s}{2}}d'\|^2_{2}+\|u_{0s}\|^2_{2}+\|\tilde{u}_{s}\|^2_{2}
\big)\nonumber\\&\quad+\frac{1}{\varepsilon^2}\big(\|u_{0}\|^2_{2}+\|e'\|^2_{2}+\|f'\|^2_{2}
+\|\tilde{u}_{\theta}\|^2_{2}\big)
\end{align*}
and
\begin{align*}
\iint (v^au)_su_{1s\theta}&=\iint (v^a)_suu_{1s\theta}+\iint v^au_su_{1s\theta}
\nonumber\\&=\iint (v^a)_s(u_{0}+u_{1}+\tilde{u})u_{1s\theta}+\iint v^a(u_{0s}+u_{1s}+\tilde{u}_{s})u_{1s\theta}
\nonumber \\
&\lesssim\eta\big(\|e^{\frac{s}{2}}c'\|^2_{2}+\|e^{\frac{s}{2}}d'\|^2_{2}\big)+\varepsilon^2\big(\|u_{0s}\|^2_{2}+\|\tilde{u}_{s}\|^2_{2}
\big)\nonumber\\&\quad+\|u_{0}\|^2_{2}+\|e'\|^2_{2}+\|f'\|^2_{2}
+\|\tilde{u}_{\theta}\|^2_{2}
\end{align*}
and
\begin{align*}
\iint \big(v(u^a-u_e(r))\big)_su_{1s\theta}&=\iint v_s(u^a-u_e(r))u_{1s\theta}+\iint v(u^a-u_e(r))_su_{1s\theta}
\nonumber\\&=\iint (v_{1s}+\tilde{v}_{s})(u^a-u_e(r))u_{1s\theta}+\iint (v_{1}+\tilde{v})(u^a-u_e(r))_su_{1s\theta}
\nonumber \\
&\lesssim\eta\big(\|e^{\frac{s}{2}}c'\|^2_{2}+\|e^{\frac{s}{2}}d'\|^2_{2}\big)+\eta\|\tilde{v}_{s}\|^2_{2}
\big)\nonumber\\&\quad+\varepsilon^{-2}\big(\|e'\|^2_{2}+\|f'\|^2_{2}
+\|\tilde{v}_{\theta}\|^2_{2}\big).
\end{align*}

Secondly, from the divergence-free condition $u_\theta+(rv)_r=0$, \eqref{weight-Hardy} and Lemma \ref{hardy-for-c-d-e-f}, we have
\begin{align*}
&\iint \big(vr(u^a-u_e(r))_r\big)_su_{1s\theta}\nonumber\\&=-\iint vr(u^a-u_e(r))_ru_{1ss\theta}\nonumber\\
&=-\iint (v_1+\tilde{v})r\big((u^a_e-u_e(r))_r+\chi'u^a_p\big)u_{1ss\theta}-\int_0^1\int_0^{2\pi} \frac{\chi vr}{r-1}(Y\partial_Yu^a_p)u_{1ss\theta}d\theta dr\\
&\lesssim\varepsilon^2(\|c''\|^2_{2}+\|d''\|^2_{2})+\frac{1}{\varepsilon^2}(\|\tilde{v}_{\theta}\|^2_{2}+\|\tilde{u}_{\theta}\|^2_{2})+\frac{1}{\varepsilon^2}(\|e'\|^2_{2}+\|f'\|^2_{2}).
\end{align*}

Finally,  we arrive at
\begin{align}
K_6&\gtrsim -\varepsilon^2(\|c''\|^2_{2}+\|d''\|^2_{2})-\frac{1}{\varepsilon^2}\big(\|\|u_{0}\|^2_{2}+\|e'\|^2_{2}+\|f'\|^2_{2}+\|\tilde{u}_{\theta}\|^2_{2}+\|\tilde{v}_{\theta}\|^2_{2}\big)
\nonumber \\
&\quad-\eta\big(\|e^{\frac{s}{2}}c'\|^2_{2}+\|e^{\frac{s}{2}}d'\|^2_{2}+\|\tilde{u}_{s}\|^2_{2}+\|\tilde{v}_{s}\|^2_{2}\big).\label{K6-est}
\end{align}

{\bf 3)Estimate of $K_7$.}
It is easy to deduce that
\begin{align}
K_7&=-\iint\big((u^a-u_e(r))v_\theta-{v}^av_s+u{v}^a_\theta+vr{v}^a_r-2u(u^a-u_e(r))\big) u_{1s\theta\theta}
\nonumber \\
&=-\iint\bigg((u^a-u_e(r))(v_{1\theta}+\tilde{v}_{\theta})-{v}^a(v_{1s}+\tilde{v}_{s})+(u_0+u_{1}+\tilde{u}){v}^a_\theta\nonumber \\
&\qquad+(v_{1}+\tilde{v})r{v}^a_r-2(u_0+u_{1}+\tilde{u})(u^a-u_e(r))\bigg) u_{1s\theta\theta}.
\nonumber\\&\gtrsim -\eta\big(\|e^{\frac{s}{2}}(c',d')\|^2_{2}+\|u_0\|^2_{2}+\|e'\|^2_{2}+\|f'\|^2_{2}+\|\tilde{v}_{\theta}\|^2_{2}+\|\tilde{v}_{s}\|^2_{2}\big).\label{K7-est}
\end{align}

Substituting \eqref{K45-est}, \eqref{K6-est} and \eqref{K7-est} into \eqref{convective-term-11} leads to
\begin{align}
&\iint\big((S_u)_s+(S_v)_\theta\big)u_{1s\theta}\nonumber\\
&\gtrsim\int_0^\infty e^s\big(c'^2+d'^2+e'^2+f'^2+(c-c_\infty)^2 +(d-d_\infty)^2+(e-e_\infty)^2+(f-f_\infty)^2\big)ds\nonumber\\
&\quad+c^2_{\infty}+d^2_{\infty}-\frac{C}{\varepsilon^2}(\|u_{0}\|^2_{2}+\|e'\|^2_{2}+\|f'\|^2_{2}+\|\tilde{u}_{\theta}\|^2_{2}+\|\tilde{v}_{\theta}\|^2_{2})
\nonumber \\
&\quad-\eta\big(\|e^{\frac{s}{2}}c'\|^2_{2}+\|e^{\frac{s}{2}}d'\|^2_{2}+\|\tilde{u}_{s}\|^2_{2}+\|\tilde{v}_{s}\|^2_{2}\big)-{\varepsilon^2}(\|e''\|^2_{2}+\|f''\|^2_{2}).\label{convective-term-111}
\end{align}

Thanks to \eqref{dif-con-term}, \eqref{cj10-11} and \eqref{convective-term-111}, we obtain the desired inequality \eqref{eq:H2-estimate-first}. Thus we complete the proof of the lemma.
\end{proof}

\subsection{Existence of error equations}
\indent

We apply the contraction mapping theorem to prove the existence of the error equations (\ref{e:error-equation}).

\begin{Proposition}\label{existence-and-error-estimate-of-error-equation}
There exist $\varepsilon_0>0,\eta_0>0$ such that for any $\varepsilon\in (0,\varepsilon_0), \eta\in (0,\eta_0)$, the error equations (\ref{e:error-equation}) have a unique solution $(u,v)$ which satisfies
$$\|(u,v)\|_\infty\leq \varepsilon^{\frac{17}{2}}.$$
\end{Proposition}
\begin{proof}
{

We consider the following linear system
\begin{align}\label{linear-equation-error-1}
\left\{
\begin{array}{lll}
-\varepsilon^2\big(\mathcal{U}_{ss}
+\mathcal{U}_{\theta\theta}+2\mathcal{V}_{\theta}-\mathcal{U}\big)+e^{-s}\mathcal{P}_\theta+e^{-s}S_\mathcal{U}=e^{-s}\mathcal{R}_u,\\[5pt]
-\varepsilon^2\big(\mathcal{V}_{ss}+\mathcal{V}_{\theta\theta}-2\mathcal{U}_{\theta}-\mathcal{V}\big)-e^{-s}\mathcal{P}_s+e^{-s}S_\mathcal{V}=e^{-s}\mathcal{R}_v,\\[5pt]
\mathcal{U}_\theta-\mathcal{V}_s+\mathcal{V}=0,  \\[5pt]
\mathcal{U}(\theta,s)=\mathcal{U}(\theta+2\pi,s),\  \mathcal{V}(\theta,s)=\mathcal{V}(\theta+2\pi,s), \\[5pt]
\mathcal{U}(\theta,0)=0,\ \mathcal{V}(\theta,0)=0,
 \end{array}
\right.
\end{align}
where $(u,v)$ satisfies
\begin{align}\label{linear-equation-error-1-new}
\left\{
\begin{array}{lll}
u_\theta-v_s+v=0,  \\[5pt]
u(\theta,s)=u(\theta+2\pi,s),\  v(\theta,s)=v(\theta+2\pi,s), \\[5pt]
u(\theta,0)=0,\ v(\theta,0)=0
 \end{array}
\right.
\end{align}
and
\begin{align*}
S_\mathcal{U}&=u^a\mathcal{U}_\theta+v^ar\mathcal{U}_r+\mathcal{U}u^a_\theta+\mathcal{V}ru^a_r+v^a\mathcal{U}+\mathcal{V}u^a,\\[5pt]
S_\mathcal{V}&=u^a\mathcal{V}_\theta+v^ar\mathcal{V}_r+\mathcal{U}v^a_\theta+\mathcal{V}rv^a_r-2\mathcal{U}u^a,\\[5pt]
\mathcal{R}_u&=R_u^a-uu_\theta-vru_r-vu,\\[5pt]
 \mathcal{R}_v&=R_v^a-uv_\theta-vrv_r+u^2.
\end{align*}

We also decompose $(\mathcal{U},\mathcal{V})$ as follows:
\begin{align*}
&\mathcal{U}(\theta,s)=\mathcal{U}_0+\underbrace{\mathcal{C}(s)\sin\theta+\mathcal{D}(s)\cos\theta}_{\mathcal{U}_1}+\tilde{\mathcal{U}},\\
&\mathcal{V}(\theta,s)=\underbrace{\mathcal{E}(s)\sin\theta+\mathcal{F}(s)\cos\theta}_{\mathcal{V}_1}+\tilde{\mathcal{V}}.
\end{align*}

Applying \eqref{eq:H2-estimate} and \eqref{eq:H2-estimate-first}  to \eqref{linear-equation-error-1} one has
\begin{align}\label{l-infinity-est-1}
&\int_0^\infty \big(\mathcal{C}''^2+\mathcal{D}''^2+\mathcal{E}''^2+\mathcal{F}''^2+(\mathcal{C}'^2+\mathcal{D}'^2+\mathcal{E}'^2+\mathcal{F}'^2)\big)ds\nonumber\\
&\lesssim\frac{1}{\varepsilon^6}(\|(\mathcal{E}',\mathcal{F}')\|^2_{2})+\frac{1}{\varepsilon^6}\|(\tilde{\mathcal{U}}_\theta
,\tilde{\mathcal{V}}_\theta)\|^2_{2}+\frac{1}{\varepsilon^6}\|(\tilde{\mathcal{U}}_s
,\tilde{\mathcal{V}}_s,\mathcal{U}_0,\mathcal{U}_{0s})\|^2_{2}
\nonumber\\&\quad+\frac{1}{\varepsilon^6}(\|e^{-s}\mathcal{R}_u\|^2_{2}+\|e^{-s}\mathcal{R}_v\|^2_{2})
\end{align}
and
\begin{align}\label{l-infinity-est-3}
&\int_0^\infty e^s\big(\mathcal{C}'^2+\mathcal{D}'^2+\mathcal{E}'^2+\mathcal{F}'^2\big)ds+\mathcal{C}^2_\infty+\mathcal{D}^2_\infty
\nonumber\\
&\ \ +\int_0^\infty  e^s\big((\mathcal{E}-\mathcal{E}_\infty)^2+(\mathcal{F}-\mathcal{F}_\infty)^2+(\mathcal{C}-\mathcal{C}_\infty)^2
+(\mathcal{D}-\mathcal{D}_\infty)^2\big)ds\nonumber\\
&\lesssim\varepsilon^2|(\mathcal{C}'',\mathcal{D}'',\mathcal{E}'',\mathcal{F}'')\|^2_{2}+\|(\mathcal{U}_{0s},\tilde{\mathcal{U}}_s
,\tilde{\mathcal{V}}_s)\|^2_{2}+\frac{1}{\varepsilon^2}\|(\mathcal{U}_0,\tilde{\mathcal{U}}_\theta
)\|^2_{2}+\frac{1}{\varepsilon^4}\|(\mathcal{E}',\mathcal{F}',\tilde{\mathcal{V}}_\theta)\|^2_{2}
\nonumber\\&\quad+\|e^{-s}\partial_sR^a_u\|^2_{2}+\iint(\mathcal{R}_u-R^a_u)_s\mathcal{U}_{1s\theta}+\|e^{-s}\mathcal{R}_v\|^2_{2}+\|e^{-s}\mathcal{R}_u\|^2_{2}.
\end{align}

By applying the linear stability estimate (\ref{e:linear-stability-estimates-for-linear-equation}) to \eqref{linear-equation-error-1}, we deduce that there exist $\varepsilon_0>0, \eta_0>0$ such that for any $\varepsilon\in (0,\varepsilon_0), \eta\in(0,\eta_0)$,
\begin{align}\label{l-infinity-est-4}
&\int_0^\infty(\mathcal{E}'^2+\mathcal{F}'^2)ds+\iint(\tilde{\mathcal{U}}_\theta^2+\tilde{\mathcal{V}}_\theta^2)+\varepsilon^2\iint\big({\tilde{\mathcal{U}}}_s^2+{\tilde{\mathcal{V}}}_s^2\big)+\varepsilon^2\int_0^\infty (\mathcal{U}^2_0+\mathcal{U}^2_{0s})ds
\nonumber\\
&\leq\iint\big(e^{-s}\mathcal{R}_u\tilde{\mathcal{U}}+ e^{-s}\mathcal{R}_v\tilde{\mathcal{V}}\big)+
\iint e^{-s}\mathcal{R}_u\mathcal{U}_0\nonumber\\
&~~~~+\iint\big(e^{-s}\mathcal{R}_u\tilde{\mathcal{U}}_\theta+ e^{-s}\mathcal{R}_v\tilde{\mathcal{V}}_\theta\big)+\iint\big(e^{-s}\mathcal{R}_u\mathcal{U}_{1\theta}+ e^{-s}\mathcal{R}_v\mathcal{V}_{1\theta}\big)\nonumber\\
&\leq\frac{1}{100}\big(\|\tilde{\mathcal{U}}_\theta\|^2_{2}+\|\tilde{\mathcal{V}}_\theta\|^2_{2}+\|\mathcal{E}'\|^2_{2}+\|\mathcal{F}'\|^2_{2}\big)+\frac{\varepsilon^2}{100}\|\mathcal{U}_0\|^2_{2}
\nonumber\\
&\quad+\frac{1}{\varepsilon^2}\|e^{-\frac{s}{2}}\mathcal{R}_u\|^2_{2}+\|e^{-\frac{s}{2}}\mathcal{R}_v\|^2_{2}
\end{align}
which immediately implies
\begin{align}\label{l-infinity-est-4-1}
&\iint(\mathcal{E}'^2+\mathcal{F}'^2)+\iint(\tilde{\mathcal{U}}_\theta^2+\tilde{\mathcal{V}}_\theta^2)+\varepsilon^2\iint\big({\tilde{\mathcal{U}}}_s^2+{\tilde{\mathcal{V}}}_s^2\big)+\varepsilon^2\int \mathcal{U}^2_0+\mathcal{U}^2_{0s}
\nonumber\\
&\lesssim\frac{1}{\varepsilon^2}\|e^{-\frac{s}{2}}\mathcal{R}_u\|^2_{2}+\|e^{-\frac{s}{2}}\mathcal{R}_v\|^2_{2}.
\end{align}

Thanks to  \eqref{l-infinity-est-1}, \eqref{l-infinity-est-3} and \eqref{l-infinity-est-4-1},  we have
\begin{align}\label{key-2}
&\big\|\big(\varepsilon\mathcal{{U}}_0,\varepsilon\mathcal{{U}}_{0s},\varepsilon\mathcal{C}'',\varepsilon\mathcal{D}''
,\varepsilon\mathcal{E}'',\varepsilon\mathcal{F}'',\varepsilon\mathcal{\tilde{U}}_s,\varepsilon\mathcal{\tilde{V}}_s,\tilde{\mathcal{U}}_\theta,\tilde{\mathcal{V}}_\theta,
\nonumber\\& \ \ \ \ e^{\frac{s}{2}}(\mathcal{C}-\mathcal{C}_\infty,\mathcal{D}-
\mathcal{D}_\infty,\mathcal{E}-\mathcal{E}_\infty,\mathcal{F}-\mathcal{F}_\infty,
\mathcal{C'},\mathcal{D'},\mathcal{E'},\mathcal{F'})\big)\big\|^2_{2}+|\mathcal{E}_\infty|^2+|\mathcal{F}_\infty|^2
\nonumber\\&\leq\frac{1}{\varepsilon^8}\|e^{-\frac{s}{2}}\mathcal{R}_u\|^2_{2}+\frac{1}{\varepsilon^6}\|e^{-\frac{s}{2}}\mathcal{R}_v\|^2_{2}
+\|e^{-s}\partial_sR^a_u\|^2_{2}+\iint(\mathcal{R}_u-R^a_u)_s\mathcal{U}_{1s\theta}.
\end{align}

Applying \eqref{eq:H2-estimate-high} to \eqref{linear-equation-error-1} one has
\begin{align}\label{l-infinity-est-2}
&\iint\tilde{\mathcal{U}}^2_{s\theta}+{\tilde{\mathcal{U}}}^2_{\theta\theta}+\tilde{\mathcal{V}}^2_{s\theta}+{\tilde{\mathcal{V}}}^2_{\theta\theta}
\nonumber\\
&\lesssim\|(\mathcal{U}_{0s},\mathcal{C}',\mathcal{D}',\tilde{\mathcal{U}}_s
,\tilde{\mathcal{V}}_s)\|^2_{2}+\frac{1}{\varepsilon^2}\|\mathcal{U}_0\|^2_{2}
+\frac{1}{\varepsilon^4}\|(\mathcal{E}',\mathcal{F}',\tilde{\mathcal{U}}_\theta,\tilde{\mathcal{V}}_\theta)\|^2_{2}
\nonumber\\&\quad+\frac{1}{\varepsilon^4}(\|e^{-s}\mathcal{R}_u\|^2_{2}+\|e^{-s}\mathcal{R}_v\|^2_{2}).
\end{align}
Thanks to \eqref{l-infinity-est-1}, \eqref{l-infinity-est-3}, \eqref{l-infinity-est-4-1}, \eqref{l-infinity-est-2} and  Lemma \ref{infinity-norm-estimate}, we then derive
\begin{align}\label{infinity-100}
\|(\mathcal{U},\mathcal{V})\|^2_{L^\infty}&\lesssim\frac{1}{\varepsilon^8}\|e^{-\frac{s}{2}}\mathcal{R}_u\|^2_{2}+\frac{1}{\varepsilon^6}\|e^{-\frac{s}{2}}\mathcal{R}_v\|^2_{2}+\|e^{-s}\partial_sR^a_u\|^2_{2}+\iint(\mathcal{R}_u-R^a_u)_s\mathcal{U}_{1s\theta}.
\end{align}

To proceed we need to  deal with the nonlinear terms including $\mathcal{R}_u$ and $\mathcal{R}_v$ in \eqref{key-2} and \eqref{infinity-100}. We only show $\|e^{-\frac{s}{2}}(\mathcal{R}_u-R^a_u)\|^2_{2}$ and $\iint(\mathcal{R}_u-R^a_u)_s\mathcal{U}_{1s\theta}$.
It is easy to obtain that
\begin{align*}
&\iint e^{-s}(uu_\theta)^2\leq\|u\|^2_{\infty}\|(e',f',\tilde{u}_\theta)\|^2_{2},
\end{align*}
and
\begin{align*}
&\iint e^{-s}(vu_s)^2\leq\|v\|^2_{\infty}\|(c',d',\tilde{u}_s)\|^2_{2},
\end{align*}
and
\begin{align*}
&\iint e^{-s}(vu)^2\leq\|u\|^2_{\infty}\|(e',f',\tilde{v}_\theta)\|^2_{2},
\end{align*}
and
\begin{align*}
&\iint (uu_\theta+vru_r+vu)_s\mathcal{U}_{1s\theta}
\nonumber\\&=\iint vu_s\mathcal{U}_{1ss\theta}+(vu_s+uv_s)\mathcal{U}_{1s\theta}
-\frac{1}{2}\partial_s(u^2)\mathcal{U}_{1s\theta\theta}\\
&\leq\frac{1}{100}\|\mathcal{C}',\mathcal{D}'\|^2_{2}+\|(v,u)\|^2_{\infty}\|(u_{0s},c',d',e',f', \tilde{u}_s,\tilde{v}_s)\|^2_{2}
\nonumber\\&\quad+\frac{1}{\varepsilon^2}\|v\|^2_{\infty}\|(u_{0s},c',d', \tilde{u}_s)\|^2_{2}
+\frac{\varepsilon^2}{100}\|(\mathcal{C}'',\mathcal{D}'')\|^2_{2}.
\end{align*}
Thus we have
\begin{align}\label{l-infinity-est-6}
&\|e^{-\frac{s}{2}}(\mathcal{R}_u-R^a_u)\|^2_{2}\leq\|u\|^2_{\infty}\|(e',f',\tilde{u}_\theta,\tilde{v}_\theta)\|^2_{2}
+\|v\|^2_{L^\infty}\|(c',d',\tilde{u}_s)\|^2_{2}
\end{align}
and
\begin{align}\label{l-infinity-est-7}
\iint(\mathcal{R}_u-R^a_u)_s\mathcal{U}_{1s\theta}
&\leq\frac{1}{100}\|\mathcal{C}',\mathcal{D}'\|^2_{2}+\|(v,u)\|^2_{\infty}\|(u_{0s},c',d',e',f', \tilde{u}_s,\tilde{v}_s)\|^2_{2}
\nonumber\\&\quad+\frac{1}{\varepsilon^2}\|v\|^2_{L^\infty}\|(u_{0s},c',d', \tilde{u}_s)\|^2_{2}
+\frac{\varepsilon^2}{100}\|(\mathcal{C}'',\mathcal{D}'')\|^2_{2}.
\end{align}

Similarly, we get
\begin{align}\label{l-infinity-est-8}
\|e^{-\frac{s}{2}}(R_v-R^a_v)\|^2_{2}\lesssim\|u\|^2_{\infty}\|( u_0,e',f',\tilde{u}_\theta,\tilde{v}_\theta)\|^2_{2}+
\|v\|^2_{\infty}\|(e',f',\tilde{v}_s)\|^2_{2}.
\end{align}

 Now we define
\begin{align*}
\|(u,v)\|^2_Y=&\|\big(\varepsilon u_0,\varepsilon
u_{0s},\varepsilon\tilde{u}_s,\varepsilon\tilde{v}_s,\tilde{u}_\theta,\tilde{v}_\theta,e^{\frac{s}{2}}(c',d',e',f',c-c_\infty,d-d_\infty,e-e_\infty,f-f_\infty)\big)\|^2_{2}
\nonumber\\&+ \|\big(\varepsilon c'',\varepsilon d'',\varepsilon e'',\varepsilon f'', \varepsilon\tilde{u}_{s\theta},\varepsilon\tilde{u}_{\theta\theta},\varepsilon\tilde{v}_{s\theta},\varepsilon\tilde{v}_{\theta\theta}\big)\|^2_{2}
\nonumber\\&+|e_\infty|^2+|f_\infty|^2+\|(u,v)\|^2_{\infty}.
\end{align*}

Let $Y_0$ be the set of smooth $2\pi$-periodic pairs satisfying
\eqref{linear-equation-error-1-new}, the homogeneous boundary condition at
$s=0$, and the first-mode limiting relations used in Hypothesis {\bf (H)}.
We define $Y$ to be the completion of $Y_0$ under the above norm.  The norm
contains the quantities appearing in Hypothesis {\bf (H)}.
For each input $(u,v)$ in a closed ball of $Y$, let
$T(u,v)=(\mathcal U,\mathcal V)$ denote the solution of the linear system
\eqref{linear-equation-error-1} with nonlinear sources
$\mathcal R_u=R_u^a-uu_\theta-vru_r-vu$ and
$\mathcal R_v=R_v^a-uv_\theta-vrv_r+u^2$.

 By \eqref{key-2}, \eqref{infinity-100}, \eqref{l-infinity-est-6}, \eqref{l-infinity-est-7} and \eqref{l-infinity-est-8}, we deduce that for small $\varepsilon, \eta$, there holds
\begin{align}\label{inductive-estimate}
\|(\mathcal{U},\mathcal{V})\|^2_Y
\leq C\bigg(\frac{1}{\varepsilon^8}\|e^{-\frac{s}{2}}{R}^a_u\|^2_2+\frac{1}{\varepsilon^6}\|e^{-\frac{s}{2}}{R}^a_v\|^2_2+\|e^{-s}\partial_sR^a_u\|^2_{2}+\frac{1}{\varepsilon^{10}}
\|(u,v)\|^4_Y\bigg).
\end{align}
Thus, due to \eqref{remainder-estimate} and the differentiated residual bound in Section~2, one has}
\begin{equation*}
\|(e^{-\frac{s}{2}}R^a_u,e^{-\frac{s}{2}}R^a_v)\|^2_2\leq C_0\varepsilon^{26}, \|e^{-s}\partial_sR^a_u\|^2_{2}\leq C_0\varepsilon^{24}
\end{equation*}
and then
there exists $\varepsilon_1>0,$ $\eta_0>0$ such that for any $\varepsilon\in(0,\varepsilon_1),$ $\eta\in(0,\eta_0),$ the operator
\begin{equation*}
 T:(u,v)\mapsto(\mathcal{U},\mathcal{V})
\end{equation*}
maps the ball $\{\|(u,v)\|^2_Y\leq \varepsilon^{17}\}$ in $Y$ into itself.

Moreover, for every two pairs
$(u, v)$ and $(\bar{u}, \bar{v})$ in the ball, we have
\begin{align}\label{contraction-estimate}
\|(\mathcal{U}-\overline{\mathcal{U}}, \mathcal{V}-\overline{\mathcal{V}})\|^2_Y\leq C\varepsilon^{-10}(\|(u,v)\|^2_Y+\|(\bar{u},\bar{v})\|^2_Y)\|(u-\bar{u}, v-\bar{v})\|^2_Y,
\end{align}
where
\begin{equation*}
  (\bar{u},\bar{v})\mapsto(\overline{\mathcal{U}},\overline{\mathcal{V}}).
\end{equation*}

In fact, set
\begin{align*}
U=\mathcal{U}-\overline{\mathcal{U}},\ V=\mathcal{V}-\overline{\mathcal{V}}, \ P=\mathcal{P}-\overline{\mathcal{P}},
\end{align*}
then there holds
\begin{align}\label{linear}
\left\{
\begin{array}{lll}
-\varepsilon^2\big({U}_{ss}
+{U}_{\theta\theta}+2{V}_{\theta}-{U}\big)+e^{-s}{P}_\theta+e^{-s}\big(S_{\mathcal{U}}-S_{\overline{\mathcal{U}}}\big)=e^{-s}\big(R_{u}-R_{\bar{u}}\big),\\[5pt]
-\varepsilon^2\big({V}_{ss}+{V}_{\theta\theta}-2{U}_{\theta}-{V}\big)-e^{-s}{P}_s+e^{-s}\big(S_{\mathcal{V}}-S_{\overline{\mathcal{V}}}\big)=e^{-s}\big(R_{v}-R_{\bar{v}}\big),\\[5pt]
{U}_\theta-{V}_s+{V}=0,  \\[5pt]
{U}(\theta,s)={U}(\theta+2\pi,s),\  {V}(\theta,s)={V}(\theta+2\pi,s), \\[5pt]
{U}(\theta,0)=0,\ {V}(\theta,0)=0.
 \end{array}
\right.
\end{align}
Thus the desired conclusion (\ref{contraction-estimate}) can be derived by following the argument leading to \eqref{inductive-estimate}.

Hence, for sufficiently small $\varepsilon$ and $\eta$, the operator
\begin{align*}
T:(u,v)\mapsto (\mathcal{U},\mathcal{V})
\end{align*}
maps the ball
$\{(u,v)\in Y:\ \|(u,v)\|_Y\leq \varepsilon^{\frac{17}{2}}\}$
into itself and is a contraction.  By the contraction mapping theorem,
$T$ has a unique fixed point $(u,v)$ in this ball.  This fixed point solves
the nonlinear error equations \eqref{e:error-equation} and satisfies
\begin{align*}
\|(u,v)\|_{Y}\leq \varepsilon^{\frac{17}{2}}.
\end{align*}
In particular, $\|(u,v)\|_{L^\infty}\leq \varepsilon^{\frac{17}{2}}$.
This completes the proof of this proposition.
\end{proof}

\section{$L^\infty$ estimate of the vorticity error and proofs of main results}
\indent

In this section we complete the proofs of the main results. First, we prove an $L^\infty$ estimate for the vorticity error.

The following trace estimate is obtained by a localized modification of the
standard Stokes regularity estimate in Tsai's book \cite{Tasi}.
\begin{Lemma}[Boundary vorticity trace estimates]\label{lem:boundary-vorticity-trace}
Let $V=\{3/4<r<1\}$ and $V'=\{2/3<r<1\}$.  For the solution of the transformed error system, let
$\omega=e^s(-u_s+u-v_\theta)$ and let $\omega_{\neq}$ denote its non-zero angular part.  The boundary contributions appearing in the first- and second-derivative vorticity identities satisfy
\[
\varepsilon^{-2}\left|\int_0^{2\pi}(\omega_{\neq})_s(\theta,0)\omega_{\neq}(\theta,0)\,d\theta\right|
+\varepsilon^{-2}\left|\int_0^{2\pi}(\omega_{\neq})_s(\theta,0)(\omega_{\neq})_{\theta\theta}(\theta,0)\,d\theta\right|
\]
\[
+\left|\int_0^{2\pi}(\omega_{\neq})_s(\theta,0)(\omega_{\neq})_\theta(\theta,0)\,d\theta\right|
\le C\varepsilon^3.
\]
The boundary terms produced when estimating $(u_0)_{ss}$, $\widetilde u_{ss}$, and $\widetilde v_{ss}$ are controlled in the same way.
\end{Lemma}
\begin{proof}
The estimates are local near $\partial B_1$.  On $V'$ the point vortex is smooth and bounded, so no singular-core correction is needed.  Rewriting the error equations in Euler coordinates and applying the standard local Stokes estimate in $V'$, followed by the trace theorem on $\partial B_1$, gives the boundary traces of $\omega_{\neq}$, $(\omega_{\neq})_s$, $(\omega_{\neq})_\theta$, and $(\omega_{\neq})_{\theta\theta}$ in terms of the $Y$-norm and the residual bounds in \eqref{eq:log-residual-bounds}.  Since $\|(u,v)\|_Y\le \varepsilon^{17/2}$ and the residual terms have order stated in \eqref{eq:log-residual-bounds}, the displayed estimate follows from Cauchy--Schwarz and Young's inequality.
\end{proof}
\begin{Proposition}\label{vorticity-behavior-000}
There exist $\varepsilon_0>0,\eta_0>0$ such that for any $\varepsilon\in (0,\varepsilon_0), \eta\in (0,\eta_0)$, the vorticity error satisfies
\begin{align} \label{vorticity-behavior-00}
\|\omega\|_\infty\leq C\varepsilon^{\frac{3}{2}}.
\end{align}
\end{Proposition}

\begin{proof}
 In order to verify \eqref{vorticity-behavior-00}, we first claim
\begin{align}\label{L-infinity-near-origin-point}
\|(u_0)_{ss},\tilde{u}_{ss},\tilde{v}_{ss}\|_{L^2}\leq C\varepsilon^{9}.
\end{align}
In fact,
it follows from the divergence-free condition and $(u,v)\in Y$ that $\|\tilde{v}_{ss}\|_{L^2}\leq C\varepsilon^{\frac{15}{2}}$.
Moreover, multiplying \eqref{c-d-e-f} by  $e^{-s}\partial_su_0$ and integrating on $[0,2\pi]\times[0,+\infty)$ implies $\|(u_0)_{ss}\|_{L^2}\leq C\varepsilon^{9}$.
 Similarly, one has  $\|\tilde{u}_{ss}\|_{L^2}\leq C\varepsilon^{9}$. The boundary terms, e.g. $u_0''(0)u_0'(0)$ and $\tilde{u}_{ss}(0)\tilde{u}_{s}(0)$, are controlled by Lemma~\ref{lem:boundary-vorticity-trace}, which is based on the trace theorem and the local Stokes estimate near $\partial B_1$. Therefore we prove \eqref{L-infinity-near-origin-point} with the help of $\|(u,v)\|_Y\leq \varepsilon^{\frac{17}{2}}$.
 \par

Next, we split the argument into two steps, treating the zero-mode part $\omega_{0}$ and the non-zero-mode part $\omega_{\neq}$ of $\omega$, respectively.

 {\bf Step 1}. In this step we focus on the zero-mode part of $\omega$ and will prove
  \begin{align}
\|\omega_{0}(s)\|_{\infty}\leq C\varepsilon^{\frac{11}{2}}, \label{zero-vorticity-estimate}
\end{align}
where
\begin{align*}
\omega_{0}(s)=\frac{1}{2\pi}\int_0^{2\pi}\omega(\theta,s)d\theta.
\end{align*}

Note that
\begin{align}
  \omega=\frac{1}{r}\partial_r(ru)-\frac{1}{r}\partial_\theta v=e^s(-u_s+u-v_\theta),\label{vorticity-def}
\end{align}
we have
\begin{align*}
\omega_{0}=e^s(-(u_0)_s+u_0).
\end{align*}
According to the first equations of \eqref{e:error-equation} and \eqref{rewrite-linear-equation-of-error}, we find
\begin{align}
 \varepsilon^2(\omega_{0})_s=(R_u)_0-(S_u)_{0}\label{ODE-vorticity}
\end{align}
which gives
\begin{align}
\varepsilon^2\omega_{0}(s)=\int^s_0\big((R_u)_0-(S_u)_{0}\big)ds'+\varepsilon^2\omega_0(0),\label{formula-zero-fre-vorticity}
\end{align}
where
\begin{align*}
(R_u)_0(s)=\frac{1}{2\pi}\int_0^{2\pi}R_u(\theta,s)d\theta,\ (S_u)_0(s)=\frac{1}{2\pi}\int_0^{2\pi}S_u(\theta,s)d\theta.
\end{align*}

Using \eqref{remainder-decay}, \eqref{decompose-u}, \eqref{decompose-u2} and  $\|(u,v)\|_Y\leq \varepsilon^{\frac{17}{2}}$, we obtain
\begin{align}
\int^s_0(R_u)_0ds'&=\frac{1}{2\pi}\int^s_0\int^{2\pi}_0(R_u^a-uu_\theta+vu_s-vu)d\theta ds'
\nonumber\\&=\frac{1}{2\pi}\int^s_0\int^{2\pi}_0(R_u^a+vu_{s}-uv_s)d\theta ds'
\nonumber\\&=\frac{1}{2\pi}\int^s_0\int^{2\pi}_0(R_u^a+v_1u_{1s}-u_1v_{1s}+\tilde{v}\tilde{u}_{s}-\tilde{u}\tilde{v}_{s})d\theta ds'
\nonumber\\&\leq\frac{C}{2\pi}\bigg(\varepsilon^{16}\int^s_0e^{-s'}ds'+\|v_1\|_{\infty}\|e^{\frac{s}{2}}u_{1s}\|_{L^2}+\|\tilde{v}_\theta\|_{L^2}\|\tilde{u}_{s}\|_{L^2}+\|\tilde{u}_\theta\|_{L^2}\|\tilde{v}_{s}\|_{L^2}\bigg)
\nonumber\\&\leq C\varepsilon^{15}.\label{est:zero-fre-Ru}
\end{align}

The following observation on the structure of $u^a$ and $v^a$ will be used:
\begin{align}
|u^a-(u^a)_0-(u^a)_1|, \ |v^a-(v^a)_1|\leq Cr\sim e^{-s},\label{construction-of-app-sol}
\end{align}
where $(\cdot)_0$ and $(\cdot)_1$ denote the zero-mode part and first Fourier mode respectively.

Using an argument similar to that for \eqref{est:zero-fre-Ru}, and combining \eqref{decompose-u}, \eqref{decompose-u2}, \eqref{construction-of-app-sol} and  $\|(u,v)\|_Y\leq \varepsilon^{\frac{11}{2}}$, we have
\begin{align}
\int^s_0(S_u)_0ds'&=\frac{1}{2\pi}\int^s_0\int^{2\pi}_0(-v^au_s+e^{-s}vu^a_r+v^au+vu^a)d\theta ds'
\nonumber\\&=\frac{1}{2\pi}\int^s_0\int^{2\pi}_0\big(-v^a(u_1+(u-u_1))_s+e^{-s}v\big(u^a-\frac{b}{r}\big)_r+v^au+vu^a\big)d\theta ds'
\nonumber\\&\leq C\varepsilon^{\frac{17}{2}},\label{est:zero-fre-Su}
\end{align}
where we have used
\begin{align*}
\int^s_0\int^{2\pi}_0v^au+vu^ad\theta ds'&=\int^s_0\int^{2\pi}_0(v^a_s-u^a_\theta)u+(v_s-u_\theta)u^ad\theta ds'\\
&=\int^s_0\int^{2\pi}_0(v^a_s u+v_su^a)d\theta ds'
\\&=-\int^s_0\int^{2\pi}_0e^{-s}v^a_rud\theta ds'+\int^s_0\int^{2\pi}_0(v_1+v-v_1)_su^ad\theta ds'
\\&\leq C\varepsilon^{\frac{17}{2}}.
\end{align*}

It follows from the trace theorem, the standard Stokes estimate near the boundary $\partial B_1$, and  \eqref{L-infinity-near-origin-point} that
\begin{align}
\|\omega_{0}(0)\|_{\infty}\leq C\varepsilon^{\frac{11}{2}}. \label{zero-vorticity-estimate-1}
\end{align}
which, together with \eqref{formula-zero-fre-vorticity}-\eqref{est:zero-fre-Su}, leads to \eqref{zero-vorticity-estimate}.

 {\bf Step 2}. In this step we focus on the non-zero mode part of $\omega$ and prove
  \begin{align}
\|\omega_{\neq}(\theta,s):=\omega(\theta,s)-\omega_{0}(s)\|_{\infty}\leq C\varepsilon^{\frac{3}{2}}. \label{nonzero-vorticity-estimate}
\end{align}

Recall that the singular component of the leading swirl is
$b/r=be^s$ in logarithmic coordinates.  We therefore isolate this part before
estimating the nonzero angular modes.  Thanks to  the equation
\eqref{c-d-e-f}  we have
\begin{align}\label{c-d-e-f-2}
&-\varepsilon^2e^{s}\big(u_{sss}
+u_{s\theta\theta}+2v_{s\theta}-u_{s}+u_{ss}+u_{\theta\theta}+2v_{\theta}-u\nonumber\\
&\quad+v_{ss\theta}+v_{\theta\theta\theta}-2u_{\theta\theta}-v_{\theta}\big)
+b(e^su_\theta)_s+b(e^sv_\theta)_\theta-2b(e^su)_\theta\nonumber\\
&=(R_u)_s+(R_v)_\theta
-(S_u-be^su_\theta)_s-(S_v-be^sv_\theta+2be^su)_\theta
=:\Lambda_b
\end{align}
and then we get the equation for the regular vorticity error
\begin{align}\label{10.01}
\varepsilon^2(\omega_{ss}+\omega_{\theta\theta})-b\omega_\theta=\Lambda_b.
\end{align}

It follows from \eqref{remainder-decay}, \eqref{decompose-u},
\eqref{decompose-u2}, $(u,v)\in Y$, \eqref{L-infinity-near-origin-point} and
Sobolev embedding that $\|\Lambda_b\|_{L^2}\leq
C\varepsilon^{\frac{15}{2}}$.
The derivative $L^\infty$ factors used below, such as $\|u_\theta\|_\infty$ and $\|v_\theta\|_\infty$, are controlled by the one-dimensional Sobolev estimates for the first Fourier mode and the two-dimensional Sobolev estimates for the higher-mode component contained in the $Y$-norm.  Thus the estimate of $\Lambda_b$ uses only quantities already controlled by $\|(u,v)\|_Y$ and the residual bounds \eqref{eq:log-residual-bounds}.
We only show the estimate of the last term in $\Lambda_b$; the remaining terms are similar. By  \eqref{decompose-u}, \eqref{decompose-u2}, \eqref{relation-inf} and $\|(u,v)\|_Y\leq \varepsilon^{\frac{17}{2}}$, we get
 \begin{align*}
\iint\big((v_\theta-u)^2+((v_\theta-u)_\theta)^2\big)
&=\iint\Big((v_{1\theta}-(v_{1\infty})_{\theta})+\tilde{v}+u_0+\tilde{u}-(u_1-u_{1\infty})\Big)^2
\nonumber\\&\quad+\iint\Big((v_{1\theta}-(v_{1\infty})_{\theta})_\theta+\tilde{v}_\theta+u_{0\theta}+\tilde{u}_\theta-(u_1-u_{1\infty})_\theta\Big)^2\\
&\leq C\frac{1}{\varepsilon^2}\|(u,v)\|_{Y}^2\nonumber\\
&\leq  C\varepsilon^{15}.
\end{align*}
Consequently, there holds
\begin{align*}
\iint (S_v-be^sv_\theta+2be^su)^2_\theta&=\iint \Big(\big(-v^av_s+vrv^a_r\big)_\theta+\big((u^a-\frac{b}{r})(v_\theta-u)+u(v^a_\theta-(u^a-\frac{b}{r}))\big)_\theta\Big)^2
\nonumber\\&\lesssim\|v_s\|_{L^2}^2+\|v_{s\theta}\|_{L^2}^2+\|v\|_{\infty}^2+\|v_\theta\|_{\infty}^2+\|v_\theta-u\|_{L^2}^2+\|(v_\theta-u)_\theta\|_{L^2}^2+\|u\|_{\infty}^2+\|u_\theta\|_{\infty}^2
\nonumber\\
&\leq C\frac{1}{\varepsilon^2}\|(u,v)\|_{Y}^2\nonumber\\
&\leq C\varepsilon^{15},
\end{align*}
where we have used
$|v^a_\theta-(u^a-\frac{b}{r})|\leq C r$ and $|\big(v^a_\theta-(u^a-\frac{b}{r})\big)_\theta|\leq C r$.

Multiplying \eqref{10.01} by $-\omega_{\neq}$ and integrating by parts, we obtain
\begin{align}
&\varepsilon^2\iint \big((\omega_{\neq})_{s}^2+((\omega_{\neq})_{\theta})^2\big)
=-\iint \Lambda_b \omega_{\neq}
-\varepsilon^2\int^{2\pi}_0 (\omega_{\neq})_{s}(\theta,0)\omega_{\neq}(\theta,0)d\theta.\label{1derivative-vorticity}
\end{align}
Furthermore, multiplying \eqref{10.01} by $(\omega_{\neq})_{ss}$ and integrating by parts, we obtain
\begin{align}
{\varepsilon^2\iint \big((\omega_{\neq})^2_{ss}+(\omega_{\neq})^2_{\theta s}\big)}
&={\varepsilon^2\int^{2\pi}_0 (\omega_{\neq})_{s}(\theta,0)(\omega_{\neq})_{\theta\theta}(\theta,0)d\theta
-\int^{2\pi}_0 (\omega_{\neq})_{s}(\theta,0)(\omega_{\neq})_{\theta}(\theta,0)d\theta}\nonumber\\
&\quad{+\iint \Lambda_b (\omega_{\neq})_{ss}.}\label{2derivative-vorticity}
\end{align}
Thanks to the Poincar\'{e} inequality, Young inequality, \eqref{1derivative-vorticity} and \eqref{2derivative-vorticity} we immediately get

\begin{align}
&\iint \big((\omega_{\neq})^2+(\omega_{\neq})_{s}^2
+(\omega_{\neq})_{\theta}^2+(\omega_{\neq})^2_{ss}
+(\omega_{\neq})^2_{\theta s}\big)\nonumber\\
&\leq C\varepsilon^{-4}\iint\Lambda_b^2
+C\varepsilon^{-2}\left|\int^{2\pi}_0
(\omega_{\neq})_{s}(\theta,0)\omega_{\neq}(\theta,0)d\theta\right|\nonumber\\
&\quad+C\varepsilon^{-2}\left|\int^{2\pi}_0
(\omega_{\neq})_{s}(\theta,0)(\omega_{\neq})_{\theta\theta}(\theta,0)d\theta\right|
+C\left|\int^{2\pi}_0
(\omega_{\neq})_{s}(\theta,0)(\omega_{\neq})_{\theta}(\theta,0)d\theta\right|\nonumber\\
&\leq C\varepsilon^{3},\label{H2derivative-vorticity-1}
\end{align}

where the boundary terms are estimated by Lemma~\ref{lem:boundary-vorticity-trace}.

On the other hand, owing to \eqref{10.01}, we have
\begin{align}
&\iint (\omega_{\neq})^2_{\theta \theta}
\leq C\varepsilon^{3}.\label{H2derivative-vorticity-2}
\end{align}
Finally, \eqref{H2derivative-vorticity-1}, \eqref{H2derivative-vorticity-2} and the Sobolev embedding give \eqref{nonzero-vorticity-estimate}; together with \eqref{zero-vorticity-estimate}, this yields the desired conclusion.
\end{proof}

Next, we prove Theorem \ref{main-theorem}.
\begin{proof}
Combining Proposition \ref{existence-and-error-estimate-of-error-equation}, Proposition \ref{vorticity-behavior-000} and the approximate solution (\ref{approximate-solution}) of the Navier--Stokes equations (\ref{NS-curvilnear}), we obtain Theorem \ref{main-theorem}.
\end{proof}

Finally, we prove Theorem \ref{main-theorem-0}.
\begin{proof}
It follows from the approximate solution and the error estimate that
\[
\mathbf{u}^\varepsilon
=\frac{b}{r}\mathbf t
+\left(u^\varepsilon-\frac{b}{r}\right)\mathbf t
+v^\varepsilon\mathbf n,
\qquad
u^\varepsilon-\frac{b}{r},\ v^\varepsilon\in L^\infty(B_1).
\]
Hence $\mathbf u^\varepsilon\in L^q(B_1)$ for every $1\leq q<2$.
Fix $1<q<2$ and set $q'=q/(q-1)>2$. Let
$\boldsymbol{\varphi}=\varphi_1\mathbf t+\varphi_2\mathbf n
\in W_0^{1,q'}(B_1)\cap W^{2,q'}(B_1)$ be divergence-free with
$\boldsymbol{\varphi}(0)=0$.  We verify the formulation in
Definition~\ref{def:very-weak-whole-disk}.

First, the nonlinear term is interpreted as a punctured-domain limit.  Write
\[
\mathbf u^\varepsilon=\frac{b}{r}\mathbf t+\mathbf u^\varepsilon_{\rm reg},
\qquad \mathbf u^\varepsilon_{\rm reg}\in L^\infty(B_1).
\]
Since $q'>2$, the embedding
$W^{2,q'}(B_1)\hookrightarrow C^{1,\alpha}(B_1)$ holds for some
$\alpha\in(0,1)$.  Hence
\[
\nabla\boldsymbol\varphi(x)=\nabla\boldsymbol\varphi(0)+O(|x|^\alpha),
\qquad
\boldsymbol\varphi(x)=O(|x|),
\]
near the origin.  All terms in
$\mathbf u^\varepsilon\otimes\mathbf u^\varepsilon:\nabla\boldsymbol\varphi$
containing at least one factor of $\mathbf u^\varepsilon_{\rm reg}$ are
integrable near $r=0$.  It remains only to consider the purely singular part.  In polar coordinates,
\[
\iint_{B_\rho\setminus B_\delta}
\frac{b^2}{r^2}\mathbf t\otimes\mathbf t:
\nabla\boldsymbol\varphi\,dxdy
=b^2\int_\delta^\rho\frac{1}{r}
\int_0^{2\pi}
\mathbf t\otimes\mathbf t:
\nabla\boldsymbol\varphi(r,\theta)\,d\theta\,dr.
\]
If $\nabla\boldsymbol\varphi(r,\theta)$ is replaced by the constant matrix
$\nabla\boldsymbol\varphi(0)$, then the angular integral equals
$\pi\operatorname{div}\boldsymbol\varphi(0)=0$.  The remaining part is
bounded by $Cr^\alpha$, and therefore the radial integrand is bounded by
$Cr^{\alpha-1}$, which is integrable at $r=0$.  Thus
\[
\lim_{\delta\to0^+}\iint_{B_1\setminus B_\delta}
\mathbf u^\varepsilon\otimes\mathbf u^\varepsilon:
\nabla\boldsymbol\varphi\,dxdy
\]
exists.

The inner convective boundary flux also vanishes.  Indeed, the singular
part is tangential, so $\mathbf u^\varepsilon\cdot\mathbf n$ is the bounded
regular normal component on $\partial B_\delta$, while
$|\boldsymbol\varphi(\theta,\delta)|\le C\delta$.  Consequently
\begin{align}
\left|
\int_{\partial B_\delta}
(\mathbf u^\varepsilon\cdot\boldsymbol\varphi)
(\mathbf u^\varepsilon\cdot\mathbf n)\,dl
\right|\le C\delta\to0.
\label{weak-solution-0-4}
\end{align}
For $\delta>0$, integrating the smooth equation on
$B_1\setminus B_\delta$ and using $\nabla\cdot\boldsymbol{\varphi}=0$ gives
\begin{align}
&-\varepsilon^2\iint_{B_1\setminus B_\delta}
\mathbf u^\varepsilon\cdot\Delta\boldsymbol{\varphi}\,dxdy
-\iint_{B_1\setminus B_\delta}
\mathbf u^\varepsilon\otimes\mathbf u^\varepsilon:
\nabla\boldsymbol{\varphi}\,dxdy\nonumber\\
&=-\varepsilon^2\int_{\partial B_1}
\mathbf g\cdot\frac{\partial\boldsymbol{\varphi}}{\partial\mathbf n}\,dl
+\mathcal B_\delta,\label{weak-solution-0-6}
\end{align}
where $\mathcal B_\delta$ denotes the inner-boundary contribution at
$\partial B_\delta$.
The diffusion part of $\mathcal B_\delta$ is generated only by the singular
field $(b/r)\mathbf t$; direct computation gives
\[
\lim_{\delta\to0^+}\mathcal B_\delta^{\rm diff}
=\varepsilon^2\big(-4\pi b\,\partial_r\overline{\varphi_1}(0)\big)
=\varepsilon^2\langle\mathbf F,\boldsymbol{\varphi}\rangle .
\]
The convective inner-boundary flux vanishes because
$\boldsymbol{\varphi}(0)=0$; the remaining terms contain the bounded regular
part and are $o(1)$ as $\delta\to0$.

It remains to check the pressure flux.  The volume pressure term has already
disappeared because the test field is divergence-free.  The singular outer
pressure satisfies
\[
p_e'(r)=\frac{(ar+b/r)^2}{r},\qquad
p_e(r)=\frac{a^2}{2}r^2+2ab\log r-\frac{b^2}{2r^2}+C.
\]
Moreover, the polar divergence-free condition
\[
\partial_\theta\varphi_1+\partial_r(r\varphi_2)=0
\]
implies, after integration in $\theta$, that
\[
\partial_r\left(r\int_0^{2\pi}\varphi_2(\theta,r)\,d\theta\right)=0.
\]
Since $\boldsymbol{\varphi}=0$ on $\partial B_1$, we have
\[
\int_0^{2\pi}\varphi_2(\theta,r)\,d\theta=0,\qquad 0<r\leq1.
\]
Consequently the singular pressure flux vanishes identically:
\[
\delta\,p_e(\delta)\int_0^{2\pi}\varphi_2(\theta,\delta)\,d\theta=0.
\]
The difference $p^\varepsilon-p_e$ is regular at the origin, so its flux is
$o(1)$.  Letting $\delta\to0^+$ in \eqref{weak-solution-0-6} yields
\begin{align}
-\varepsilon^2\iint_{B_1}
\mathbf u^\varepsilon\cdot\Delta\boldsymbol{\varphi}\,dxdy
-\lim_{\delta\to0^+}\iint_{B_1\setminus B_\delta}
\mathbf u^\varepsilon\otimes\mathbf u^\varepsilon:
\nabla\boldsymbol{\varphi}\,dxdy
&=\varepsilon^2\langle\mathbf F,\boldsymbol{\varphi}\rangle
-\varepsilon^2\int_{\partial B_1}
\mathbf g\cdot
\frac{\partial\boldsymbol{\varphi}}{\partial\mathbf n}\,dl.
\label{weak-solution-0-0}
\end{align}

For the divergence equation, let $\psi\in W^{1,q'}(B_1)$.  Since the
singular component is tangential and the regular normal component is
bounded near the origin,
\begin{align}
\iint_{B_1}\mathbf u^\varepsilon\cdot\nabla\psi\,dxdy
&=\lim_{\delta\to0^+}
\iint_{B_1\setminus B_\delta}
\mathbf u^\varepsilon\cdot\nabla\psi\,dxdy
=\lim_{\delta\to0^+}
\int_{\partial B_\delta}v^\varepsilon\psi\,dl=0.
\label{weak-solution-0-00}
\end{align}
Thus $\mathbf u^\varepsilon$ is a very weak solution on the whole disk.
Finally, Proposition~\ref{existence-and-error-estimate-of-error-equation}
and Proposition~\ref{vorticity-behavior-000} show that the regular velocity
and regular vorticity extend across the origin, so
\eqref{velocity-behavior-01}--\eqref{vorticity-behavior-0} yield
\eqref{velocity-behavior-1}--\eqref{vorticity-behavior}.
\end{proof}

\appendix

\section{Construction of corrector $h(\theta,r)$}

In this section, we give a construction of the corrector $h(\theta,r)$. First, we prove a simple lemma.

\begin{Lemma}\label{corector-equation}
Assume that $K(\theta,r)$ is a $2\pi$-periodic smooth function which satisfies
\beno
\int_0^{2\pi}K(\theta,r)d\theta=0, \ \forall r\in (0,1]; \ K(\theta, 1)=0,
\eeno
then there exists a $2\pi$-periodic function $h(\theta,r)$ such that
\begin{align}\label{eq:corrector-h}
&\partial_\theta h(\theta,r)=K(\theta,r); \quad h(\theta, 1)=0;\nonumber\\
&\int_0^{2\pi}h(\theta,r)d\theta=0, \quad \left(\int_0^1\int_0^{2\pi}\left|\partial_\theta^j\partial_r^kh\right|^2d\theta dr\right)^{\frac12}\leq C\left(\int_0^1\int_0^{2\pi}\left|\partial_\theta^j\partial_r^kK\right|^2d\theta dr\right)^{\frac12}.
\end{align}
\end{Lemma}
\begin{proof}
Let
\beno
K(\theta,r)=\sum_{n\neq 0}K_n(r)e^{in\theta}, \ \ K_n(1)=0.
\eeno
Set
\beno
h(\theta,r)=\sum_{n\neq 0}\frac{K_n(r)}{in}e^{in\theta}.
\eeno
It is easy to justify that $h(\theta,r)$ satisfies (\ref{eq:corrector-h}) which completes the proof.
\end{proof}

Next, we construct the corrector $h(\theta,r)$ consistently with the
thirteenth-order approximate solution.
\begin{Proposition}\label{prop:corrector-h}
The function $h$ may be chosen so that
\[
u^a_\theta+rv^a_r+v^a=0,\qquad h(\theta,1)=0,\qquad
\int_0^{2\pi}h(\theta,r)\,d\theta=0,
\]
and, for all fixed $j,k\geq0$,
\[
\left(\int_0^1\int_0^{2\pi}\left|\partial_\theta^j\partial_r^k h\right|^2d\theta dr\right)^{\frac12}\leq C(j,k)\varepsilon^{-k}.
\]
\end{Proposition}

\begin{proof}
Direct computation gives
\begin{align*}
u^a_\theta+rv^a_r+v^a
=\varepsilon^{13}\partial_\theta h(\theta,r)+K(\theta,r),
\end{align*}
where
\begin{align*}
K(\theta,r)=&\varepsilon^{13}\chi(r)
\big[Y\partial_Yv_p^{(13)}(\theta,Y)+v_p^{(13)}(\theta,Y)\big]
+r\chi'(r)\Big(\sum_{i=1}^{13}\varepsilon^i v_p^{(i)}(\theta,Y)\Big).
\end{align*}
Since $\chi'(r)=0$ on $[0,\frac12]\cup[\frac34,1]$, the cut-off
contribution is supported where $Y=(r-1)/\varepsilon\leq-1/(4\varepsilon)$.
The exponential decay of the Prandtl profiles therefore implies that, for
each fixed $j,k$,
\[
\left(\int_0^1\int_0^{2\pi}\left|\partial_\theta^j\partial_r^kK\right|^2d\theta dr\right)^{\frac12}\leq C(j,k)\varepsilon^{13-k}.
\]
Moreover, $K(\theta,1)=0$ and, using
$\int_0^{2\pi}v_p^{(i)}(\theta,Y)\,d\theta=0$ for all $1\leq i\leq13$,
\[
\int_0^{2\pi}K(\theta,r)\,d\theta=0,\qquad 0<r\leq1.
\]
By Lemma \ref{corector-equation} we choose $h$ so that
\[
\varepsilon^{13}\partial_\theta h(\theta,r)+K(\theta,r)=0.
\]
The preceding bound on $K$ gives
\[\left(\int_0^1\int_0^{2\pi}\left|\partial_\theta^j\partial_r^kh\right|^2d\theta dr\right)^{\frac12}\leq C(j,k)\varepsilon^{-k}\]
and hence
$(u^a,v^a)$ is exactly divergence-free.
\end{proof}

\par
\par\noindent
\section*{{Funding and Acknowledgments}}
\par\noindent
Z. Chen was supported by the National Natural Science Foundation of China (No. 12501294).
M. Fei is supported partly by NSF of China under Grant No.12271004, No.12471222 and NSF of Anhui
Province of China under Grant No. 2308085J10. Z. Lin is supported partly by NSF of China under Grant No. 12494544. He is a member of LMNS (Laboratory of Mathematics for Nonlinear Science) of Fudan University. He thanks LMNS for the support.
J. Zhao did not receive any funding for this study.
\par\noindent
 \section*{{Statements and Declarations}}
 \par\noindent
\textbf{Data availability statement}
There is no associated data to the manuscript.
\par\noindent
\textbf{Conflict of interest}
The authors declare no conflict of interest.

\end{document}